\documentclass[a4paper,11pt,leqno]{article}
\author{Mo\"ise J\'ER\'EMIE}
\date{}
\usepackage{amsmath}
\usepackage{amsfonts}
\usepackage{amssymb}
\usepackage{pifont}
\usepackage{amsthm}
\usepackage[nottoc,notlof,notlot]{tocbibind}
\usepackage[english,frenchb]{babel}
\usepackage[colorlinks=true]{hyperref}
\hypersetup{urlcolor=blue,linkcolor=black,citecolor=black,colorlinks=true}
\usepackage[latin1]{inputenc}
\usepackage{titlesec}
\newtheorem{teo}{Th\'eor\`eme}[section]
\newtheorem{prop}{Proposition}[section]
\newtheorem{cor}{Corollaire}[section]
\newtheorem{rmq}{Remarque}[section]
\newtheorem{expl}{Exemple}[section]
\newtheorem{lemme}{Lemme}[section]

\newtheorem{fait}{Fait}
\newcommand{\ud}{\mathrm{d}}
\numberwithin{equation}{section}

\titleformat{\chapter}[display]
  {\normalfont\huge\bfseries\centering}{\chaptertitlename\ \thechapter}{20pt}{\Huge}

\titleformat{\section}
  {\normalfont\Large\bfseries\centering}{\thesection}{1em}{}
\titleformat{\subsection}
  {\normalfont\large\bfseries\centering}{\thesubsection}{1em}{}

\begin{document}
\begin{center}
{\bf Processus empiriques des rapports de $m$-espacements uniformes disjoints}\\ ({\bf Non-overlapping uniform $m$-spacings-ratio empirical processes})\\
\vspace{0,3cm}Mo\"ise JEREMIE\\
\vspace{0,3cm}L.S.T.A., Universit\'e de Paris VI, Jussieu\\
\vspace{0,5cm}
\end{center}

\begin{center}
{\bf R\'esum\'e}
\end{center}

\noindent {\scriptsize Nous considérons un processus empirique basé sur des paires sélectionnées de rapports de $m$-es\-pa\-ce\-ments
disjoints, engendrés par des échantillons in\-dé\-pen\-dants de tailles arbitraires. Comme résultat principal, nous montrons
que lorsque les deux échantillons sont uniformément distribués sur des intervalles de même longueur,
ce processus empirique converge vers un pont Brownien centré sur une moyenne de la forme $\ \displaystyle (B\circ H_{m})_{C}(v)=B(H_{m}(v))-2(2m+1)C\ \frac{(2m-1)!}{m((m-1)!)^2}(v(1-v))^m\int_{0}^{1}B(H_{m}(s))\ud s,\ $ pour $\ \displaystyle 0\leq v\leq 1,\ $ où $\ \displaystyle B(.)\ $ désigne un pont Brownien, $\ \displaystyle H_{m},\ $ la f.r. de la loi béta de paramètres $m$ et $m,\ $ et $C,\ $ une constante.}

\begin{center}
{\bf Abstract}
\end{center}

\noindent {\scriptsize We consider an empirical process based upon ratio of selected pair of the
non-overlapping $m$-spacings generated by independent samples of arbitrary sizes.
As a main result, we show that when both samples are uniformly distributed on
intervals of equal lengths, this empirical process converges to a mean-centered
Brownian bridge of the form $\ \displaystyle (B\circ H_{m})_{C}(v)=B(H_{m}(v))-2(2m+1)C\ \frac{(2m-1)!}{m((m-1)!)^2}(v(1-v))^m\int_{0}^{1}B(H_{m}(s))\ud s,\ $ for $\ \displaystyle 0\leq v\leq 1,\ $ where $\ \displaystyle B(.)\ $ denotes a Brownian bridge, $\ \displaystyle H_{m},\ $ the distribution function of the Beta distribution with parameters $m$ and $m,\ $ and $C,\ $ a constant.}\vskip5pt

\vskip5pt

\noindent {\scriptsize\textit{Keys words and phrases:} processus empiriques, processus Gaussiens, principes d'invariance faibles, lois limites, espacements, statistiques d'ordre, lois faibles.}

\section{Introduction et résultats}\label{ma-sectiona}

\noindent Deheuvels et Derzko \cite{MR2728435} ont étudié un processus
empirique basé sur les rapports d'espacements de deux échantillons indépendants.
Dans le cas d'une loi uniforme sur $(0,1)$, ils ont
montré que ce processus
empirique convergeait vers une limite gaussienne de la forme
\[B(t)-6Ct(1-t)\int_{0}^1B(s)ds,\quad \hbox{pour}\quad 0\leq t\leq 1,\]
où $\{B(t):0\leq t\leq 1\}$ est un pont Brownien, $C$, une
constante dépendant du comportement limite des tailles des deux
échantillons. Dans cet article, nous généralisons leurs
résultats aux $m$-espacements. Nous adopterons les
notations suivantes. Sur $(\Omega,\mathcal{A},\mathbb{P})$, $\{X_n:n\geq 1\}$ et $\{Y_n:n\geq
1\}$ deux suites indépendantes, de variables aléatoires
indépendantes de loi uniforme sur $(0,1),$ et de fonctions de
répartition $F_{X}(t)=P(X_{1}\leq t)=F_{Y}(t)=P(Y_{1}\leq t)=t$,
pour $0\leq t\leq 1$. Pour tout couple d'entiers, $n_1\geq 1$ et
$n_2\geq 1,\ $ soient respectivement $\ \displaystyle
X_{1,n_1}<\ldots<X_{n_1,n_1}\ $ et $\ \displaystyle
Y_{1,n_2}<\ldots<Y_{n_2,n_2}\ $ les statistiques d'ordre de $\
\displaystyle X_1,\ldots,X_{n_1}$ et $\ \displaystyle
Y_1,\ldots,Y_{n_2}$. Nous posons $\ \displaystyle
X_{0,n_1}=0,\ $  $\ \displaystyle X_{n_1+1,n_1}=1\ $ et $\
\displaystyle Y_{0,n_2}=0,$  $\displaystyle Y_{n_2+1,n_2}=1,\ $ pour
$\ \displaystyle n_1,n_2\geq 0$. Puisque ces statistiques d'ordre
sont distinctes presque sûrement, nous travaillerons, sans perte de
généralité, sur l'événement de probabilité $1$ correspondant. Pour
tout $1\leq m\leq (n_1\wedge n_2)+1,$ nous
considérons les $m$-espacements engendrés par chacune de ces suites,
et définis par
\[D_{i,n_1;X}^{(m)}=X_{i+m,n_1}-X_{i,n_1}\quad\hbox{et}\quad D_{j,n_2;Y}^{(m)}=Y_{j+m,n_2}-Y_{j,n_2},\]
pour $0\leq i\leq n_1-m+1$, et $0\leq j\leq n_2-m+1$. L'objet de notre étude est l'ensemble des rapports d'espacements de la forme
\[R_{i,j,n_1,n_2}=\frac{(n_1-m+2)D_{i,n_1;X}^{(m)}}{(n_1-m+2)D_{i,n_1;X}^{(m)}+(n_2-m+2)D_{j,n_2;Y}^{(m)}}.\]
La loi de $R_{i,j,n_1,n_2}$ est indépendante de $i,j$, et satisfait, lorsque $n_1\wedge n_2\rightarrow\infty$,
\[R_{i,j,n_1,n_2}\stackrel{d}{\rightarrow}\beta_{m,m},\]
(voir~(\ref{eq:iaaaaaaaw}) ci-dessous pour la définition de $\beta_{m,m}).$ Nous allons considérer un processus empirique basé sur une collection de ces rapports, choisis de telle sorte que les espacements qui y interviennent soient disjoints au sens que les intervalles ouverts correspondants sont d'intersection vide. Cette propriété est toujours satisfaite pour $m=1$, cas étudié par
Deheuvels et Derzko \cite{MR2435038,MR2728435}. A cet effet, nous introduisons des indices $0\leq i_{0,n_1}<\ldots<i_{N,n_1}
\leq n_1-m+1$, et $0\leq j_{0,n_2}<\ldots<j_{N,n_2}\leq n_2-m+1$, tels que $i_{k+1,n_1}-i_{k,n_1}\geq m$ et
$j_{k+1,n_2}-j_{k,n_2}\geq m,$ pour $0\leq k\leq N\leq N_1\land N_2,\ $ où
\begin{equation}
N_1:=\left\lfloor\frac{n_1+1}{m}\right\rfloor-1\quad\hbox{ et }\quad N_2:=\left\lfloor\frac{n_2+1}{m}\right\rfloor-1,\label{eq:iaaaaaaaaf}
\end{equation}
et $\ \displaystyle \left\lfloor u\right\rfloor\ $ désigne la partie entière inférieure de $u\in\mathbb{R}.\ $\vskip5pt
\noindent Nous considérons donc les paires de $m$-espacements disjoints définies par
\begin{flalign}
&S_{k,N;X}^{(m)}=D_{i_{k,n_{1}},n_1;X}^{(m)}=X_{i_{k,n_{1}}+m,n_1}-X_{i_{k,n_{1}},n_1}\ \mbox{ et }\label{eq:iaaaaaaar}\\
&S_{k,N;Y}^{(m)}=D_{j_{k,n_{2}},n_2;Y}^{(m)}=Y_{j_{k,n_{2}}+m,n_2}-Y_{j_{k,n_{2}},n_2},\ \mbox{ pour }\ k=0,\ldots,N.\nonumber
\end{flalign}
Par convention, nous posons
\begin{equation}
P=P(N)=N_{1}-N\geq 0\ \mbox{ et }\ Q=Q(N)=N_{2}-N\geq 0,\label{eq:iaaaaaaaag}
\end{equation}
et supposons que $\ \displaystyle N,\ P=P(N)\ \mbox{ et }\ Q=Q(N)\ $ vérifient les conditions limites, lorsque
$\ \displaystyle N\to\infty,$
\begin{equation}
\frac{P}{N+1}\to c\in[0,\infty],\ \mbox{ et }\ \frac{Q}{N+1}\to d\in[0,\infty],\label{eq:iaaaaaaas}
\end{equation}
\noindent où $\ \displaystyle c\ \mbox{ et }\ d \ $ sont des
constantes (éventuellement infinies).
Pour $0\leq k\leq N\leq N_1\wedge N_2$, les rapports de
$m$-espacements disjoints sont définis par
\begin{flalign}
R_{k;n_1,n_2}=R_{i_{k,n_1},j_{k,n_2},n_1,n_2}&=\frac{(N+P+1)\ S_{k,N;X}^{(m)}}{(N+P+1)\ S_{k,N;X}^{(m)}+(N+Q+1)\ S_{k,N;Y}^{(m)}}\label{eq:iaaaaaaat}\\
&=\frac{\left(N_{1}+1\right)\ S_{k,N;X}^{(m)}}{\left(N_{1}+1\right)\ S_{k,N;X}^{(m)}+\left(N_{2}+1\right)\ S_{k,N;Y}^{(m)}},\nonumber
\end{flalign}
La loi de $\{R_{k;n_1,n_2}:0\leq k\leq N\}$ ne dépend pas de  $\{i_{k,n_1},j_{k,n_2}:0\leq k\leq N\}$. Nous notons la fonction de répartition $\ \displaystyle [\mbox{f.r.}]\ $ empirique
correspondante par
\begin{equation}
H_{N;n_1,n_2}(x)=\frac{1}{N+1}\sum_{k=0}^N{\rm 1\!I}_{\left\{R_{k;n_1,n_2}\leq x\right\}},
\quad\hbox{pour}\quad x\in \mathbb{R}.\label{eq:iaaaaaaau}
\end{equation}
Notons $H_{m}$, la f.r. de la loi $\beta_{m,m}$, définie, pour $\ \displaystyle m>0,\ $ par
\begin{equation}
H_{m}(x)=\frac{1}{\beta(m,m)}\int_{0}^xt^{m-1}(1-t)^{m-1}dt,\quad\hbox{pour}\quad 0\leq x\leq 1.\label{eq:iaaaaaaaw}
\end{equation}Notons aussi, $\ \displaystyle H_{m}^{-1},\ $ la fonction de quantile de la loi $\ \displaystyle \beta_{m,m},\ $ définie, pour $\ \displaystyle m>0,\ $ par
\begin{equation}
H_{m}^{-1}(t)=\inf\left\{x:H_{m}(x)\geq t\right\},\ \hbox{pour tout }\quad 0<t<1.\label{eq:iaaaaaaaw1}
\end{equation}
Nous abordons dans ce qui suit l'étude du processus empirique
\begin{equation}
\hspace{-1cm}\gamma_{N;n_1,n_2}(x)=\left(N+1\right)^{1/2}(H_{N;n_1,n_2}(x)-H_{m}(x)),\quad\hbox{pour}\quad 0\leq x\leq 1,\label{eq:iaaaaaaav}
\end{equation}
dont nous étudions le comportement sous les conditions~(\ref{eq:iaaaaaaas}).\vskip5pt
\noindent Comme résultat principal, nous montrons que lorsque les v.a. $\ \displaystyle \left\{X_{n}:\right.$ $\displaystyle\left.n\geq 1\right\}\ $ et $\ \displaystyle \left\{Y_{n}:n\geq 1\right\}\ $ suivent des lois uniformes sur des intervalles de $\ \displaystyle \mathbb{R},\ $
$\ \displaystyle\left\{\gamma_{N;n_1,n_2}(v):\right.$ $ \displaystyle\left.0\leq v\leq 1\right\}\ $ peut être approché asymptotiquement par {\it des processus Gaussiens centrés} (voir Deheuvels 2007\cite{MR2392790}) de la forme
{\small
\begin{flalign}
\hspace{-0,7cm}(B\circ H_{m})_{C}(v)&=B(H_{m}(v))-2(2m+1)C\ \frac{(2m-1)!}{m((m-1)!)^2}(v(1-v))^m\label{eq:mv}\\
\hspace{-0,7cm}&\quad\times\int_{0}^{1}B(H_{m}(s))\ud s,\ 0\leq v \leq 1.\nonumber
\end{flalign}
}\noindent Ici $\ \displaystyle \left\{B(v):0\leq v\leq 1\right\}\ $ désigne un pont Brownien et $\ \displaystyle C\in\mathbb{R},\ $ une constante. Cette catégorie de processus dont nous rappellerons les propriétés au paragraphe~\ref{ma-sectionb} ont été introduites par Deheuvels et Derzko~\cite{MR2728435}. Nous établirons, en particulier l'identité en loi
{\footnotesize
\begin{flalign}
&\left\{(B\circ H_{m})(t):0\leq t\leq 1\right\}\stackrel{d}{=}\left\{(B\circ H_{m})(t)\right.\label{eq:mw}\\
&\left.\hspace{1,5cm}-4(2m+1)\frac{(2m-1)!}{m((m-1)!)^2}(t(1-t))^m\int_{0}^{1}B(H_{m}(s))\ud s: 0\leq t\leq 1\right\}.\nonumber
\end{flalign}}
\noindent Dans le cas suivant, la structure limite de $\ \displaystyle \gamma_{N;n_1,n_2}^{(m)}(.)\ $ est simple.
\begin{expl}\label{exp:ex1} Si $N=o(N_{1})$ et $N=o(N_{2})$, lorsque $N_1\wedge N_2\rightarrow\infty$,
de sorte que $c=d=\infty$ dans~(\ref{eq:iaaaaaaas}),
il existe, pour un choix convenable de $\left(\Omega,\mathcal{A},\mathbb{P}\right)$,
une suite $\{B_{N}^{\prime}(.):N\geq1\}$ de
ponts Browniens, telle que, lorsque $N\to\infty$,
\begin{equation}
\sup_{0\leq v\leq 1}\left\vert\gamma_{N;n_1,n_2}(v)-B_{N}^{\prime}(H_{m}(v))\right\vert=o_{\mathbb{P}}(1).\label{eq:iaaaaaaaaad}
\end{equation}
\end{expl}
\noindent La démonstration de~(\ref{eq:iaaaaaaaaad}) sera précisée, plus loin, dans la remarque~\ref{rmque:remqu1}.\vskip5pt
\noindent L'exemple~\ref{exp:ex1} laisse entendre que le processus
empirique $\gamma_{N;n_1,n_2}(.)$ pourrait converger faiblement
(relativement à la topologie de Skorohod, Bil\-ling\-sley \cite{MR0233396}) vers un
pont Brownien. Tel n'est pas le cas en général com\-me l'ont établi Deheuvels et Derzko~\cite{MR2728435} pour $\ \displaystyle m=1.\ $
Le théorème~\ref{thm:teo11112} suivant généralise leur résultat à $\ \displaystyle m\geq 1.$
\begin{teo}\label{thm:teo11112} Lorsque $c=d=0$, pour un choix convenable de
$(\Omega,\mathcal{A},\mathbb{P})$, il existe deux suites $\ \displaystyle \left\{{\mathcal{B}}_{N}^{+}(v):0\leq v\leq 1\right\}\ $ et $\ \displaystyle \left\{{\mathcal{B}}_{N}^{-}(v):0\leq v\leq 1\right\},\ $ $\ \displaystyle N=1,2,\ldots,\ $ de ponts Browniens, liés par les relations réciproques,
{\small
\begin{flalign}
{\mathcal{B}}_{N}^{\pm}\circ H_{m}(v)&\stackrel{d}{=}{\mathcal{B}}_{N}^{\pm}\circ H_{m}(v)-8(2m+1)\ \frac{(2m-1)!}{2m\ ((m-1)!)^2}(v(1-v))^m\times\label{eq:mt}\\
&\quad\int_{0}^{1}{\mathcal{B}}_{N}^{\pm}\circ H_{m}(s)\ \ud s,\ \mbox{ pour }\ 0\leq v\leq 1,\nonumber
\end{flalign}
}\noindent et vérifiant les relations, lorsque $\ \displaystyle N\to\infty,\ $
\begin{flalign}
&\hspace{-1,2cm}\sup_{0\leq v\leq 1}\left\vert \gamma_{N;n_1,n_2}(v)-\left\{{\mathcal{B}}_{N}^{+}\circ H_{m}(v)-4(2m+1)\ \frac{(2m-1)!}{2m\ \left((m-1)!\right)^2}\times\right.\right.\label{eq:iaaaaaaaab}\\
&\hspace{-0,7cm}\left.\left.\left(1+\frac{1}{\sqrt{2m+1}}\right)\left(v(1-v)\right)^m\int_{0}^{1}{\mathcal{B}}_{N}^{+}\circ H_{m}(s)\ \ud s\right\}\right\vert=o_{\mathbb{P}}(1),\nonumber
\end{flalign}
\noindent et
\begin{flalign}
&\hspace{-1,2cm}\sup_{0\leq v\leq 1}\left\vert \gamma_{N;n_1,n_2}(v)-\left\{{\mathcal{B}}_{N}^{-}\circ H_{m}(v)-4(2m+1)\ \frac{(2m-1)!}{2m\ \left((m-1)!\right)^2}\times\right.\right.\label{eq:mu}\\
&\hspace{-0,7cm}\left.\left.\left(1-\frac{1}{\sqrt{2m+1}}\right)\left(v(1-v)\right)^m\int_{0}^{1}{\mathcal{B}}_{N}^{-}\circ H_{m}(s)\ \ud s\right\}\right\vert=o_{\mathbb{P}}(1).\nonumber
\end{flalign}
\end{teo}
\noindent La preuve du théorème~\ref{thm:teo11112} est exposée dans le paragraphe~\ref{ma-soussectionn}.\vskip5pt
\begin{rmq}\label{rmque:remqu2} Soit $\ \displaystyle (B\circ H_{m})_{C}(\cdot)\ $ comme dans~(\ref{eq:mv}). Il est facile de montrer (voir par exemple~(\ref{eq:ms}) ci-dessous) que
\begin{equation}
(B\circ H_{m})_{C_1}(\cdot)\stackrel{d}{=}(B\circ H_{m})_{C_2}(\cdot)\iff C_{1}(C_{1}-2)=C_{2}(C_{2}-2).\label{eq:mx}
\end{equation}
\noindent Ceci implique, en particulier que les versions des processus gaussiens\\*
$\displaystyle(B\circ H_{m})_{\displaystyle\Bigg\{1+\frac{1}{\sqrt{2m+1}}\Bigg\}}(\cdot)\ $ et $\ \displaystyle(B\circ H_{m})_{\displaystyle\Bigg\{1-\frac{1}{\sqrt{2m+1}}\Bigg\}}(\cdot)\ $ approchant $\ \displaystyle\gamma_{N;n_1,n_2}(\cdot)\ $ dans~(\ref{eq:iaaaaaaaab}) et~(\ref{eq:mu}), respectivement, sont identiquement distribués (voir le paragraphe~\ref{ma-sectionb}).
\end{rmq}
\noindent Le théorème~\ref{thm:teo11112} peut être généralisé comme suit.\vskip5pt
\begin{teo}\label{thm:teo11113} Sur un espace $(\Omega,\mathcal{A},\mathbb{P})$ de probabilité
convenable, il existe deux suites $\ \displaystyle \left\{{\mathcal{B}}_{N}^{+}(v):0\leq v\leq 1\right\}\ $ et $\ \displaystyle \left\{{\mathcal{B}}_{N}^{-}(v):0\leq v\leq 1\right\},\ $ $\ \displaystyle N=1,2,\ldots,\ $ de ponts Browniens, liés par les relations réciproques,
{\small
\begin{flalign}
{\mathcal{B}}_{N}^{\pm}\circ H_{m}(v)&\stackrel{d}{=}{\mathcal{B}}_{N}^{\pm}\circ H_{m}(v)-8(2m+1)\ \frac{(2m-1)!}{2m\ ((m-1)!)^2}(v(1-v))^m\times\label{eq:iaaaaaaaaae3}\\
&\quad\int_{0}^{1}{\mathcal{B}}_{N}^{\pm}\circ H_{m}(s)\ \ud s,\ \mbox{ pour }\ 0\leq v\leq 1,\nonumber
\end{flalign}
}\noindent et vérifiant les relations, lorsque $\ \displaystyle N\to\infty,\ $
{\footnotesize
\begin{flalign}
&\hspace{0cm}\sup_{0\leq v\leq 1}\left\vert \gamma_{N;n_1,n_2}(v)-\left\{{\mathcal{B}}_{N}^{+}\circ H_{m}(v)-4(2m+1)\ \frac{(2m-1)!}{2m\ \left((m-1)!\right)^2}\left(1+\sqrt{R_{N,m}}\right)\right.\right.\label{eq:my}\\
&\left.\left.\hspace{0,5cm}\times\left(v(1-v)\right)^m\int_{0}^{1}{\mathcal{B}}_{N}^{+}\circ H_{m}(s)\ \ud s\right\}\right\vert=O_{\mathbb{P}}\left(N^{-1/4}\ \left(\log N\right)^{1/2}\right),\nonumber
\end{flalign}
}\noindent et
{\footnotesize
\begin{flalign}
&\hspace{0cm}\sup_{0\leq v\leq 1}\left\vert \gamma_{N;n_1,n_2}(v)-\left\{{\mathcal{B}}_{N}^{-}\circ H_{m}(v)-4(2m+1)\ \frac{(2m-1)!}{2m\ \left((m-1)!\right)^2}\left(1-\sqrt{R_{N,m}}\right)\right.\right.\label{eq:mz}\\
&\left.\left.\hspace{0,5cm}\times\left(v(1-v)\right)^m\int_{0}^{1}{\mathcal{B}}_{N}^{-}\circ H_{m}(s)\ \ud s\right\}\right\vert=O_{\mathbb{P}}\left(N^{-1/4}\ \left(\log N\right)^{1/2}\right),\nonumber
\end{flalign}
}\begin{flalign}
\hspace{-1cm}\mbox{où }\ R_{N,m}&=1-\frac{m}{2m+1}\left[\frac{N+1}{N+1+P}+\frac{N+1}{N+1+Q}\right]\label{eq:maa}\\
\hspace{-1cm}&=1-\frac{m}{2m+1}\left[\frac{N+1}{N_1+1}+\frac{N+1}{N_2+1}\right].\nonumber
\end{flalign}
\end{teo}
\noindent Dans le théorème~\ref{thm:teo11114} suivant, nous considérons le cas où $\ \displaystyle \left\{X_{n}:n\geq 1\right\}\ $ et $\ \displaystyle \left\{Y_{n}:n\geq 1\right\}\ $ sont des suites indépendantes de v.a. i.i.d. de lois uniformes respectivement sur $\ \displaystyle(a,a+e)\ $ et sur $\ \displaystyle(b,b+h).\ $\vskip5pt
\begin{teo}\label{thm:teo11114} Si $\ \displaystyle X\stackrel{d}{=}\mathcal{U}(a,a+e)\ $ et $\ \displaystyle Y\stackrel{d}{=}\mathcal{U}(b,b+h).\ $ Alors, sur un espace $(\Omega,\mathcal{A},\mathbb{P})$ de probabilité
convenable, il existe deux suites $\ \displaystyle \left\{{\mathcal{B}}_{N}^{+}(v):\right.$ $\displaystyle\left.0\leq v\leq 1\right\}\ $ et $\ \displaystyle \left\{{\mathcal{B}}_{N}^{-}(v):0\leq v\leq 1\right\},\ $ $\ \displaystyle N=1,2,\ldots,\ $ de ponts Browniens, liés par les relations réciproques,
{\small
\begin{flalign}
&\hspace{0cm}{\mathcal{B}}_{N}^{\pm}\circ H_{m}\left(\frac{vh}{vh+(1-v)e}\right)\stackrel{d}{=}{\mathcal{B}}_{N}^{\pm}\circ H_{m}\left(\frac{vh}{vh+(1-v)e}\right)-8(2m+1)\times\label{eq:mab}\\
&\hspace{0,5cm}\frac{(2m-1)!}{2m\ ((m-1)!)^2}\left(\frac{vh(1-v)e}{\left(vh+(1-v)e\right)^2}\right)^m
\int_{0}^{1}{\mathcal{B}}_{N}^{\pm}\circ H_{m}(s)\ \ud s,\nonumber
\end{flalign}
}\noindent pour $\ \displaystyle0\leq v\leq 1,\ $ et vérifiant les relations, lorsque $\ \displaystyle N\to\infty,\ $
{\footnotesize
\begin{flalign}
&\hspace{0,5cm}\sup_{0\leq v\leq 1}\left\vert \gamma_{N;n_1,n_2}\left(\frac{vh}{vh+(1-v)e}\right)-\left\{{\mathcal{B}}_{N}^{+}\circ H_{m}\left(\frac{vh}{vh+(1-v)e}\right)-4(2m+1)\times\right.\right.\label{eq:mac}\\
&\hspace{1cm}\left.\left.\frac{(2m-1)!}{2m\ \left((m-1)!\right)^2}\left(1+\sqrt{R_{N,m}}\right)\left(\frac{vh(1-v)e}{\left(vh+(1-v)e\right)^2}\right)^m
\int_{0}^{1}{\mathcal{B}}_{N}^{+}\circ H_{m}(s)\ \ud s\right\}\right\vert\nonumber\\
&\hspace{1cm}=O_{\mathbb{P}}\big(N^{-1/4}\ \left(\log N\right)^{1/2}\big),\nonumber
\end{flalign}
}\noindent et
{\footnotesize
\begin{flalign}
&\hspace{0cm}\sup_{0\leq v\leq 1}\left\vert \gamma_{N;n_1,n_2}\left(\frac{vh}{vh+(1-v)e}\right)-\left\{{\mathcal{B}}_{N}^{-}\circ H_{m}\left(\frac{vh}{vh+(1-v)e}\right)-4(2m+1)\times\right.\right.\label{eq:mad}\\
&\hspace{0,5cm}\left.\left.\frac{(2m-1)!}{2m\ \left((m-1)!\right)^2}\left(1-\sqrt{R_{N,m}}\right)\left(\frac{vh(1-v)e}{\left(vh+(1-v)e\right)^2}\right)^m\int_{0}^{1}{\mathcal{B}}_{N}^{-}\circ H_{m}(s)\ \ud s\right\}\right\vert\nonumber\\
&\hspace{0,5cm}=O_{\mathbb{P}}\big(N^{-1/4}\ \left(\log N\right)^{1/2}\big),\nonumber
\end{flalign}
}\noindent où
\begin{flalign}
\hspace{-1cm}R_{N,m}&=1-\frac{m}{2m+1}\left[\frac{N+1}{N+1+P}+\frac{N+1}{N+1+Q}\right]\label{eq:mae}\\
\hspace{-1cm}&=1-\frac{m}{2m+1}\left[\frac{N+1}{N_1+1}+\frac{N+1}{N_2+1}\right].\nonumber
\end{flalign}
\end{teo}
\begin{rmq}\label{rmque:remqu1}1°) En posant, comme dans~(\ref{eq:iaaaaaaas}), $\ \displaystyle P/(N+1)\to c\ $ et $\ \displaystyle Q/(N+1)\to d,\ $
lorsque $\ N\to\infty$, nous obtenons que
\begin{displaymath}
R_{N,m}\to R_{\infty,m}:=1-\frac{m}{2m+1}\left[\frac{1}{1+c}+\frac{1}{1+d}\right],
\end{displaymath}
\noindent avec la convention $1/\infty=0.\ $ Sous cette hypothèse, nous pouvons remplacer $\ \displaystyle R_{N,m}\ $ par $\ \displaystyle R_{\infty,m}\ $ dans~(\ref{eq:my}) et~(\ref{eq:mz}) au prix du remplacement du terme $O_{\mathbb{P}}\big(N^{-1/4}(\log N)^{1/2}\big)$ par $o_{\mathbb{P}}(1)$. Lorsque $c=d=0$, ce dernier résultat coïncide avec le théorème~\ref{thm:teo11112}. De même, nous obtenons
~(\ref{eq:iaaaaaaaaad}), pour $\ c=d=\infty.\ $\vskip5pt
\noindent 2°) Nous avons toujours $\ \displaystyle 1-\frac{2m}{2m+1}<R_{N,m}<1,\ $ et de même, $\ \displaystyle R_{\infty,m}\ $ peut être pris pour chaque valeur possible dans $\ \displaystyle \left[1-\frac{2m}{2m+1},1\right].\ $ Le seul cas où le processus limite de $\ \gamma_{N;n_1,n_2}(.),\ $ lorsque $\ \displaystyle N\to\infty,\ $ est un pont Brownien
est obtenu pour $\ \displaystyle R_{\infty,m}=1,\ $ ou équivalent, lorsque $\ \displaystyle c=d=\infty.\ $\vskip5pt
\noindent 3°) Lorsque $m=1$, les théorèmes~\ref{thm:teo11112} et~\ref{thm:teo11113} précisent les théorèmes 1.1 et 1.2 de Deheuvels et Derzko\cite{MR2728435} comme suit.\vskip5pt
\begin{teo}[Deheuvels et Derzko 2011]\label{thm:teo8} Lorsque $c=d=0$, pour un choix convenable de
$(\Omega,\mathcal{A},\mathbb{P})$, il existe deux suites $\ \displaystyle \left\{{\mathcal{B}}_{N}^{+}(v):0\leq v\leq 1\right\}\ $ et $\ \displaystyle \left\{{\mathcal{B}}_{N}^{-}(v):0\leq v\leq 1\right\},\ $ $\ \displaystyle N=1,2,\ldots,\ $ de ponts Browniens, liés par les relations réciproques,
\begin{equation}
{\mathcal{B}}_{N}^{\pm}(v)\stackrel{d}{=}{\mathcal{B}}_{N}^{\pm}(v)-12v(1-v)\int_{0}^{1}{\mathcal{B}}_{N}^{\pm}(s)\ \ud s,\ \mbox{ pour }\ 0\leq v\leq 1,\label{eq:mt2}
\end{equation}
\noindent et vérifiant les relations, lorsque $\ \displaystyle N\to\infty,\ $
\begin{flalign}
&\hspace{-1cm}\sup_{0\leq v\leq 1}\left\vert \gamma_{N;n_1,n_2}(v)-\left\{{\mathcal{B}}_{N}^{+}(v)-2(3+\sqrt{3})v(1-v)\int_{0}^{1}{\mathcal{B}}_{N}^{+}(s)\ \ud s\right\}\right\vert\label{eq:iaaaaaaaab2}\\
&\hspace{-1cm}=o_{\mathbb{P}}(1),\nonumber
\end{flalign}
\noindent et
\begin{flalign}
&\hspace{-1cm}\sup_{0\leq v\leq 1}\left\vert \gamma_{N;n_1,n_2}(v)-\left\{{\mathcal{B}}_{N}^{-}(v)-2(3-\sqrt{3})v(1-v)\int_{0}^{1}{\mathcal{B}}_{N}^{-}(s)\ \ud s\right\}\right\vert\label{eq:mu2}\\
&\hspace{-1cm}=o_{\mathbb{P}}(1).\nonumber
\end{flalign}
\end{teo}
\begin{teo}[Deheuvels et Derzko 2011]\label{thm:teo9} Sur un espace $(\Omega,\mathcal{A},\mathbb{P})$ de probabilité
convenable, il existe deux suites $\ \displaystyle \left\{{\mathcal{B}}_{N}^{+}(v):0\leq v\leq 1\right\}\ $ et\\*
$\displaystyle \left\{{\mathcal{B}}_{N}^{-}(v):0\leq v\leq 1\right\},\ $ $\ \displaystyle N=1,2,\ldots,\ $ de ponts Browniens, liés par les relations réciproques,
\begin{equation}
{\mathcal{B}}_{N}^{\pm}(v)\stackrel{d}{=}{\mathcal{B}}_{N}^{\pm}(v)-12v(1-v)\int_{0}^{1}{\mathcal{B}}_{N}^{\pm}(s)\ \ud s,\ \mbox{ pour }\ 0\leq v\leq 1,\label{eq:iaaaaaaaaae4}
\end{equation}
\noindent et vérifiant les relations, lorsque $\ \displaystyle N\to\infty,\ $
{\small
\begin{flalign}
&\hspace{0cm}\sup_{0\leq v\leq 1}\left\vert \gamma_{N;n_1,n_2}(v)-\left\{{\mathcal{B}}_{N}^{+}(v)-6\left(1+\sqrt{R_{N,1}}\right)v(1-v)\int_{0}^{1}{\mathcal{B}}_{N}^{+}(s)\ \ud s\right\}\right\vert\label{eq:my2}\\
&\hspace{0cm}=O_{\mathbb{P}}\left(N^{-1/4}\ \left(\log N\right)^{1/2}\right),\nonumber
\end{flalign}
}\noindent et
{\small
\begin{flalign}
&\hspace{0cm}\sup_{0\leq v\leq 1}\left\vert \gamma_{N;n_1,n_2}(v)-\left\{{\mathcal{B}}_{N}^{-}(v)-6\left(1-\sqrt{R_{N,1}}\right)v(1-v)\int_{0}^{1}{\mathcal{B}}_{N}^{-}(s)\ \ud s\right\}\right\vert\label{eq:mz2}\\
&\hspace{0cm}=O_{\mathbb{P}}\left(N^{-1/4}\ \left(\log N\right)^{1/2}\right),\nonumber
\end{flalign}
{\small
}\begin{equation}
\hspace{-1cm}\mbox{où }\ R_{N,1}=1-\frac{1}{3}\left[\frac{N+1}{N+1+P}+\frac{N+1}{N+1+Q}\right]=1-\frac{1}{3}\left[\frac{N+1}{N_1+1}+\frac{N+1}{N_2+1}\right].\label{eq:maa2}
\end{equation}}
\end{teo}
\noindent 4°) Nous notons que, pour $\ \displaystyle e=h,\ $ le théorème~\ref{thm:teo11114} se réduit au théorème~\ref{thm:teo11113}.
\end{rmq}
\noindent Dans le paragraphe~\ref{ma-sectionb} qui suit, nous rappelons des
résultats sur la classe des processus gaussiens centrés étudiée par Deheuvels et Derzko\cite{MR2728435},
qui inclut les processus limites des théorèmes~\ref{thm:teo11112},~\ref{thm:teo11113} et~\ref{thm:teo11114}. Dans le paragraphe~\ref{ma-sectionc}, nous présentons les arguments principaux de nos démonstrations. Le paragraphe~\ref{ma-sectiond}
est consacré à l'approximation forte des processus empiriques étudiés. Enfin, au paragraphe~\ref{ma-sectione}, nous assemblons les résultats précédents pour établir les théorèmes~\ref{thm:teo11112},~\ref{thm:teo11113} et~\ref{thm:teo11114}.
\section{Processus gaussiens centrés sur une moyenne}\label{ma-sectionb}
\noindent Pour $\ \displaystyle d\geq 1,\ $ considérons un processus gaussien $\ \displaystyle \left\{\mathfrak{X}(t):t\in[0,1]^d\right\},\ $ avec des trajectoires continues et une fonction de covariance $\ \displaystyle K(s,t)=\mathbb{E}\left(\mathfrak{X}(s)\right.$ $\displaystyle\left.\times\mathfrak{X}(t)\right),\ $ vérifiant la condition
\begin{equation}
0<\sigma_{\mathfrak{X}}^2:=\int_{[0,1]^d\times[0,1]^d} K(s,t)\ \ud s\ud t<\infty.\label{eq:ma}
\end{equation}
\noindent Par le théorème 3, p.104 dans Talagrand (1987)\cite{MR906527}, la continuité de $\ \displaystyle \mathfrak{X}(.)\ $ sur $\ \displaystyle [0,1]^d\ $ implique la continuité de $\ \displaystyle K(.,.)\ $ sur $\ \displaystyle [0,1]^d\times[0,1]^d, $ et donc la convergence de l'intégrale~(\ref{eq:ma}). Nous supposons que $\ \displaystyle \sigma^{2}_{\mathcal{X}}>0,\ $  équivaut au fait que $\ \displaystyle \mathfrak{X}(.)\ $ n'est pas identiquement nul sur $\ \displaystyle [0,1]^d.\ $ Pour tout $\ \displaystyle t=(t_{1},\ldots,t_{d})\in \mathbb{R}^d\ $ ou ($\ \displaystyle s=(s_{1},\ldots,s_{d})\in \mathbb{R}^d\ $), nous désignons par $\ \displaystyle \ud t\ $ ou ($\ \displaystyle \ud s\ $) la mesure de Lebesgue dans $\ \displaystyle \mathbb{R}^d.\ $ De plus, en général, nous posons pour chaque fonction mesurable $\ \displaystyle \phi(.)\ $ sur $\ \displaystyle[0,1]^d,\ $
\begin{displaymath}
\int_{[0,1]^d\times\ldots\times[0,1]^d}\phi(s)\ \ud s=\int_{0}^{t}\phi(s)\ \ud s=\int_{0}^{t_{1}}\ldots\int_{0}^{t_{d}}\phi(s_{1},\ldots,s_{d})\ \ud s_{1}\ldots \ud s_{d}.
\end{displaymath}
\noindent Introduisons la fonction de $\ \displaystyle t\in[0,1]^d\ $ définie par
\begin{displaymath}
\Psi(t):=\int_{[0,1]^d} K(s,t)\ \ud s=\mathbb{E}\left(\mathfrak{X}(t)\int_{[0,1]^d}\mathfrak{X}(s)\ \ud s\right),\ \mbox{ pour }\ t\in[0,1]^d.
\end{displaymath}
\noindent En gardant à l'esprit, compte tenu de~(\ref{eq:ma}) que
\begin{equation}
\sigma_{\mathcal{X}}^{2}=\int_{[0,1]^d} \Psi(t)\ \ud t={\rm Var}\left(\int_{[0,1]^d} \mathfrak{X}(s)\ \ud s\right),\label{eq:mb}
\end{equation}
\noindent nous notons $\ \displaystyle \mathfrak{D}=\mathfrak{D}(\mathfrak{X})=\left\{\mathfrak{X}_{C}:C\in\mathbb{R}\right\},\ $ la classe  de tous les processus gaussiens de la forme
\begin{equation}
{\mathfrak{X}}_{C}(t):=\mathfrak{X}(t)-\frac{C\ \Psi(t)}{\sigma_{\mathfrak{X}}^{2}}\int_{[0,1]^d} \mathfrak{X}(s)\ \ud s,\ \mbox{ pour }\ t\in[0,1]^d,\label{eq:mc}
\end{equation}
\noindent où $\ \displaystyle C\in\mathbb{R}\ $ est une constante. Les processus gaussiens de cette forme sont appelés {\it processus gaussiens centrés sur une moyenne}.\vskip5pt
\noindent $\displaystyle \mathfrak{E}=\mathfrak{E}(\mathfrak{X})=\left\{\lambda{\mathfrak{X}}_{C}:\ \lambda\in\mathbb{R},\ C\in\mathbb{R}\right\}\ $ coïncide avec l'espace vectoriel engendré par $\ \displaystyle\mathfrak{D}.\ $ Pour chaque $\ \displaystyle C\in\mathbb{R},\ $ nous définissons une application linéaire $\ \displaystyle J_{C}:\mathfrak{Y}\in\mathfrak{E}\longrightarrow {\mathfrak{Y}}_{C}\in\mathfrak{E}\ $ de $\displaystyle \mathfrak{E}\ $ sur lui, même par $\ \displaystyle \mathfrak{Y}\in\mathfrak{E},$
\begin{equation}
\hspace{-1cm}\left\{J_{C}\ \mathfrak{Y}\right\}(t)={\mathfrak{Y}}_{C}(t)=\mathfrak{Y}(t)-\frac{C\ \Psi(t)}{\sigma_{\mathfrak{Y}}^{2}}\int_{[0,1]^d} \mathfrak{Y}(s)\ \ud s,\ \mbox{ pour }\ t\in[0,1]^d.\label{eq:md}
\end{equation}
\noindent Les propriétés de $\ \displaystyle \left\{J_{C}:C\in\mathbb{R}\right\}\ $ sont précisées dans le lemme suivant. Nous désignons par $\ \displaystyle f\circ g=f(g)\ $ la composition de $\ \displaystyle g\ $ par  $\ \displaystyle f,\ $ et par $\ \displaystyle f^{-1},\ $ l'élément inverse de $\ \displaystyle f\ $ (quand il existe).
\begin{lemme}
\label{lem:le900}
\noindent Nous avons les identités suivantes. Pour chaque $\ \displaystyle A,B,C\in \mathbb{R},$
\begin{flalign}
&J_{A}\circ J_{B}=J_{A}\circ J_{B}=J_{A+B-AB};\label{eq:me}\\
&J_{A}\circ (J_{B}\circ J_{C})=(J_{A}\circ J_{B})\circ J_{C}=J_{A+B+C-AB-AC-BC+ABC};\label{eq:mf}\\
&J_{0}\circ J_{A}=J_{A}\circ J_{0}=J_{A};\ J_{A}\circ J_{B}=J_{1}\iff A=1\ \mbox{ ou }\ B=1;\label{eq:mg}\\
&J_{1}\circ J_{A}=J_{A}\circ J_{1}=J_{1};\ A\neq 1\Longrightarrow J_{A}^{-1}=J_{A/(A-1)}\label{eq:mh}\\
&J_{A}+J_{B}=2J_{(A+B)/2}.\label{eq:mi}
\end{flalign}
\end{lemme}
\noindent{\bf Démonstration.} Simple donc omise.$\Box$\vskip5pt
\noindent Le lemme~\ref{lem:le900} montre que $\ \displaystyle \left\{J_{C}:C\neq 1\right\}\ $ est un groupe abélien d'endomorphismes de $\ \displaystyle\mathfrak{E}.\ $ Ce groupe est une représentation du groupe
commutatif défini sur $\ \displaystyle \mathbb{R}-\left\{1\right\}\ $ par le produit interne $\ \displaystyle x*y=x+y-xy.\ $\vskip5pt
\begin{lemme}
\label{lem:le901}
\noindent Pour tout $\ \displaystyle C\in\mathbb{R},\ $ nous avons les identités en loi
\begin{equation}
J_{C}({\mathfrak{X}}_{0})={\mathfrak{X}}_{C}\stackrel{d}{=}J_{C}({\mathfrak{X}}_{2})={\mathfrak{X}}_{2-C}.\label{eq:mj}
\end{equation}
\end{lemme}
\noindent{\bf Démonstration.} Observons que
\begin{equation}
K_{C}(s,t)=\mathbb{E}\left({\mathfrak{X}}_{C}(s){\mathfrak{X}}_{C}(t)\right)=K(s,t)+(C^2-2C)\left\{\frac{\Psi(s)\Psi(t)}
{\sigma^{2}_{\mathfrak{X}}}\right\}.\label{eq:mk}
\end{equation}

\noindent La conclusion suit du fait que $\ \displaystyle (2-C)^2-2(2-C)=C^2-2C,\ $ qui implique, via~(\ref{eq:mk}) que $\ \displaystyle{\mathfrak{X}}_{2-C}\ $ et $\ \displaystyle{\mathfrak{X}}_{C}\ $ ont des fonctions de covariance identiques.$\Box$\vskip5pt
\noindent Le lemme simple suivant jouera un rôle instrumental dans la suite.\vskip5pt
\begin{lemme}
\label{lem:le902}
\noindent Le processus gaussien $\ \displaystyle\left\{\mathfrak{X}_{1}(t):\ t\in[0,1]^d\right\}\ $ et la variable aléatoire
\begin{displaymath}
\int_{[0,1]^d}\mathfrak{X}(s)\ \ud s\stackrel{d}{=}N(0,\sigma^2_{\mathfrak{X}})
\end{displaymath}
\noindent sont indépendants.\vskip5pt
\end{lemme}
\noindent {\bf Démonstration.} Nous observons que, pour chaque $\ \displaystyle t\in[0,1]^d,\ $

\begin{displaymath}
\mathbb{E}\left({\mathfrak{X}}_{C}(t)\int_{[0,1]^d}\mathfrak{X}(s)\ \ud s\right)=\Psi(t)-C\Psi(t),
\end{displaymath}

\noindent qui s'annule pour $\ \displaystyle C=1.\ $ Ceci suffit pour nos besoins. $\Box$\vskip5pt
\noindent Le théorème suivant est plus ou moins direct et facile, étant donné le lemme~\ref{lem:le902}.\vskip5pt
\begin{teo}\label{thm:teo1mpgc} Soit $\ \displaystyle Y\stackrel{d}{=}N(0,\sigma^2)\ $ qui désigne une variable aléatoire gaussienne indépendante de $\ \displaystyle \left\{\mathfrak{X}(t):t\in[0,1]^d\right\}.$ Alors, pour chaque $\ \displaystyle C\in\mathbb{R},\ $ nous avons l'égalité distributionnelle
\begin{equation}
\left\{{\mathfrak{X}}_{C}(t)+\Psi(t)Y:t\in[0,1]^d\right\}\stackrel{d}{=}\left\{{\mathfrak{X}}_{1\pm\sqrt{(1-C)^2+\sigma^{2}
\sigma^2_{\mathfrak{X}}}}(t):t\in[0,1]^d\right\}.\label{eq:ml}
\end{equation}
\end{teo}
\noindent {\bf Démonstration.} Nous utilisons la décomposition

\begin{displaymath}
{\mathfrak{X}}_{C}(t)+\Psi(t)Y={\mathfrak{X}}_{1}(t)+\frac{\Psi(t)}{\sigma^2_{\mathfrak{X}}}
\left\{Y\sigma^2_{\mathfrak{X}}+(1-C)\int_{[0,1]^d}\mathfrak{X}(s)\ \ud s\right\}=:{\mathfrak{X}}_{1}(t)+\frac{\Psi(t)}{\sigma^2_{\mathfrak{X}}}\ Z,
\end{displaymath}

\noindent où compte tenu de~(\ref{eq:mb}) et du lemme~\ref{lem:le902}, $\displaystyle Z\ $ est une v.a. normale centrée, indépendante de $\ \displaystyle {\mathfrak{X}}_{1}(.),\ $ avec une distribution donnée par

\begin{displaymath}
Z\stackrel{d}{=}-Z\stackrel{d}{=}N\left(0,\sigma_{\mathfrak{X}}^{2}\left\{(1-C)^2+\sigma^{2}\sigma^{2}_{\mathfrak{X}}\right\}\right)
\stackrel{d}{=}\pm\sqrt{(1-C)^2+\sigma^{2}\sigma^{2}_{\mathfrak{X}}}\int_{[0,1]^d}\mathfrak{X}(s)\ \ud s.
\end{displaymath}

\noindent Ceci à son tour, implique les relations

\begin{flalign*}
{\mathfrak{X}}_{C}(t)+\Psi(t)Y&={\mathfrak{X}}_{1}(t)+\frac{\Psi(t)}{\sigma^2_{\mathfrak{X}}}\ Z\\
&\stackrel{d}{=}{\mathfrak{X}}_{1}(t)\pm\frac{\Psi(t)}{\sigma^2_{\mathfrak{X}}}\sqrt{(1-C)^2+\sigma^2\sigma^2_{\mathfrak{X}}}
\int_{[0,1]^d}\mathfrak{X}(s)\ \ud s\\
&=\mathfrak{X}(t)-\frac{\Psi(t)}{\sigma^2_{\mathfrak{X}}}\left\{1\pm\sqrt{(1-C)^2+\sigma^2\sigma^2_{\mathfrak{X}}}\right\}
\int_{[0,1]^d}\mathfrak{X}(s)\ \ud s,
\end{flalign*}

\noindent ce qui donne~(\ref{eq:ml}), comme voulu.$\Box$\vskip5pt
\noindent Le corollaire suivant du théorème~\ref{thm:teo1mpgc} est d'un intérêt en soi.\vskip5pt
\begin{cor}\label{coro:cor1mpgc} Nous avons l'identité distributionnelle
\begin{equation}
\left\{\mathfrak{X}(t):t\in[0,1]^d\right\}\stackrel{d}{=}\left\{\mathfrak{X}(t)-\frac{2\Psi(t)}{\sigma^{2}_{\mathfrak{X}}}
\int_{[0,1]^d}\mathfrak{X}(s)\ \ud s:t\in[0,1]^d\right\}.\label{eq:mm}
\end{equation}
\end{cor}
\noindent{\bf Démonstration.} En posant $\ \displaystyle C=0\ $ et $\ \displaystyle \sigma^2=0\ $ dans~(\ref{eq:ml}), nous obtenons que $\ \displaystyle \mathfrak{X}(.)={\mathfrak{X}}_{0}(.)\stackrel{d}{=}{\mathfrak{X}}_{2}(.)\ $  qui est~(\ref{eq:mm}). Cette égalité suit aussi de ~(\ref{eq:mj}), pris avec $\ \displaystyle C=0.$ $\Box$\vskip5pt
\begin{rmq}\label{rmque:remqu1} Dans le présent article, nous serons principalement concernés avec le cas, où, pour un entier $\ \displaystyle m\geq 1\ $ fixé, $\ \displaystyle\mathfrak{X}(.)=B(H_{m}(.))\ $ est un pont Brownien, $\ \displaystyle H_{m}\ $ est la f.r. de la loi $\beta_{m,m}$, définie, pour $\ \displaystyle 0\leq x\leq 1,\ $ par
\begin{eqnarray}\label{eq:mn}
H_{m}(x)&:=& \frac{\Gamma(2m)}{(\Gamma(m))^2}\int_{0}^{x}(1-t)^{m-1}t^{m-1}\ud t\\
&=&\sum_{j=m}^{2m-1}\frac{(2m-1)!}{j!(2m-1-j)!}x^j(1-x)^{2m-1-j}\nonumber\\
&=&1-\sum_{j=0}^{m-1}\frac{(2m-1)!}{j!(2m-1-j)!}x^j(1-x)^{2m-1-j},
\nonumber
\end{eqnarray}et $\ \displaystyle d=1,$ car si $\ \displaystyle X^{\prime}\ $ est une variable aléatoire qui suit
une loi $\ \beta_{a,b},\ $ avec $\ \displaystyle a,b>0,\ $ alors sa fonction de répartition,
la fonction béta régularisée $\ \displaystyle I_{t}(a,b)\ $ est définie, pour tout $\ \displaystyle 0\leq t\leq 1,\ $ par

\begin{flalign}
I_{t}(a,b)&=\frac{\Gamma(a+b)}{\Gamma(a)\Gamma(b)}\int_{0}^{t}x^{a-1}\ (1-x)^{b-1}\ \ud x=\frac{B(t;a,b)}{B(a,b)}\label{eq:ppaak1}\\
&=\sum_{j=a}^{a+b-1}\frac{(a+b-1)!}{j!\ (a+b-1-j)!}\ t^j\ (1-t)^{a+b-1-j}\nonumber\\
&=1-\sum_{j=0}^{a-1}\frac{(a+b-1)!}{j!\ (a+b-1-j)!}\ t^j\ (1-t)^{a+b-1-j},\nonumber
\end{flalign}

\noindent où $\ \displaystyle B(t;a,b)\ $ est la fonction béta incomplète définie, pour tout $\ \displaystyle 0\leq t\leq 1,\ $ par

\begin{displaymath}
B(t;a,b)=\int_{0}^{t}x^{a-1}\ (1-x)^{b-1}\ \ud x\ \mbox{ où }\ B(1;a,b)=B(a,b).
\end{displaymath}

\noindent On a, alors dans ce cas,
\begin{displaymath}
K(s,t)=\mathbb{E}\left(B(H_{m}(s))B(H_{m}(t))\right)=H_{m}(s)\land H_{m}(t)-H_{m}(s)H_{m}(t),
\end{displaymath}pour $\ \displaystyle\ 0\leq s,t\leq 1,\ $ afin que, pour tout $\ \displaystyle t\in]0,1[,\ $
\begin{flalign*}
\Psi(t)&=\int_{[0,1]}K(s,t)\ \ud s=\int_{[0,1]}\mathbb{E}\left(B(H_{m}(s))B(H_{m}(t))\right)\ \ud s\\
&=\int_{0}^{1}\left\{H_{m}(s)\land H_{m}(t)-H_{m}(s)H_{m}(t)\right\}\ud s\\
&=\int_{0}^{1}\left\{H_{m}(s\land t)-H_{m}(s)H_{m}(t)\right\}\ud s.
\end{flalign*}

\noindent Soit $\ \displaystyle U^{\prime}\stackrel{d}{=}\mathcal{U}(0,1)\ $ et posons
\begin{flalign*}
\hspace{-1cm}&W^{\prime}=\int_{0}^{1}\left(H_{m}(s)-{\rm 1\!I}_{\left\{\displaystyle U^{\prime}\leq H_{m}(s)\right\}}\right)\ \ud s\ \mbox{ et pout tout }\ t\in[0,1],\\
\hspace{-1cm}&V^{\prime}=H_{m}(t)-{\rm 1\!I}_{\left\{\displaystyle U^{\prime}\leq H_{m}(t)\right\}}.
\end{flalign*}

\noindent Notons que
\begin{flalign*}
\hspace{-1cm}&W^{\prime}=H_{m}^{-1}(U^{\prime})-\mathbb{E}H_{m}^{-1}(U^{\prime})\ \mbox{ et pour tout }\ t\in[0,1],\\
\hspace{-1cm}&V^{\prime}=H_{m}(t)-{\rm 1\!I}_{\left\{\displaystyle H_{m}^{-1}(U^{\prime})\leq t\right\}}.
\end{flalign*}
\noindent Nous voyons facilement que, pour tout $\ \displaystyle t\in(0,1),\ $
{\small
\begin{displaymath}
\Psi(t)={\rm Cov }\left(W^{\prime},V^{\prime}\right)=-\mathbb{E}H_{m}^{-1}(U^{\prime}){\rm 1\!I}_{\left\{\displaystyle H_{m}^{-1}(U^{\prime})\leq t\right\}}+\mathbb{E}H_{m}^{-1}(U^{\prime})\mathbb{E}{\rm 1\!I}_{\left\{\displaystyle H_{m}^{-1}(U^{\prime})\leq t\right\}}.
\end{displaymath}}Or $\ H_{m}^{-1}(U^{\prime})\stackrel{d}{=}\beta_{m,m},\ $ donc, on a,
{\small
\begin{flalign}
\hspace{0cm}\Psi(t)&={\rm Cov }\left(W^{\prime},V^{\prime}\right)\label{eq:mo}\\
\hspace{-0cm}&=-\mathbb{E}H_{m}^{-1}(U^{\prime}){\rm 1\!I}_{\left\{\displaystyle H_{m}^{-1}(U^{\prime})\leq t\right\}}+\mathbb{E}H_{m}^{-1}(U^{\prime})\mathbb{E}{\rm 1\!I}_{\left\{\displaystyle H_{m}^{-1}(U^{\prime})\leq t\right\}}\nonumber\\
\hspace{0cm}&=-\frac{1}{B(m,m)}\int_{0}^{t}x^{m+1-1}(1-x)^{m-1}\ \ud x\nonumber\\
&\quad+\frac{H_{m}(t)}{B(m,m)}\int_{0}^{1}x^{m+1-1}(1-x)^{m-1}\ \ud x\nonumber\\
\hspace{0cm}&=\frac{B(m+1,m)}{B(m,m)}\left(-\frac{1}{B(m+1,m)}\int_{0}^{t}x^{m+1-1}(1-x)^{m-1}\ \ud x
+H_{m}(t)\right)\nonumber\\
\hspace{0cm}&=\frac{1}{2}\left(-\sum_{j=m+1}^{2m}\frac{(2m)!}{j!\ (2m-j)!}\ t^j\ (1-t)^{2m-j}\right.\nonumber\\
\hspace{0cm}&\left.\quad+\sum_{j=m}^{2m-1}\frac{(2m-1)!}{j!\ (2m-1-j)!}\ t^j\ (1-t)^{2m-1-j}\right)\nonumber\\
\hspace{0cm}&=\frac{(2m-1)!}{2m((m-1)!)^2}\ (t(1-t))^m,\nonumber
\end{flalign}}et
{\small
\begin{flalign}
\sigma^{2}_{\mathfrak{X}}&=\sigma^{2}_{B\circ H_{m}}\label{eq:mn2}\\
&={\rm Var}\left(\int_{0}^{1} B(H_{m}(x))\ \ud x\right)=\mathbb{E}\left(\int_{0}^{1} B(H_{m}(x))\ \ud x\int_{0}^{1} B(H_{m}(y))\ \ud y \right)\nonumber\\
&=\int_{0}^{1}\int_{0}^{1}\mathbb{E}\left(B(H_{m}(x))B(H_{m}(y))\right)\ \ud x\ud y\nonumber\\
&=\int_{0}^{1}\int_{0}^{1}(H_{m}(x)\land H_{m}(y)-H_{m}(x)H_{m}(y))\ \ud x\ud y\nonumber\\
&=\int_{0}^{1}\int_{0}^{1}(H_{m}(x\land y)-H_{m}(x)H_{m}(y))\ \ud x\ud y={\rm Var}(Y^{\prime})=\frac{1}{4(2m+1)},\nonumber
\end{flalign}}où $\ \displaystyle Y^{\prime}\ $ est une v.a. de loi $\ \displaystyle \beta_{m,m},\ $ c'est-à-dire $\ \displaystyle Y^{\prime}\stackrel{d}{=}\beta_{m,m}.\ $\vskip5pt

\noindent Nous obtenons aussi comme un cas particulier de~(\ref{eq:mm}), l'identité distributionnelle inattendue (voir par exemple, (1.3), dans [12, p. 1190])
\begin{flalign}
&\left\{(B\circ H_{m})(t):0\leq t\leq 1\right\}\label{eq:mp}\\
&\stackrel{d}{=}\left\{(B\circ H_{m})(t)-4(2m+1)\frac{(2m-1)!}{m((m-1)!)^2}(t(1-t))^m\right.\nonumber\\
&\left.\quad\quad\quad\times\int_{0}^{1}B(H_{m}(s))\ud s: 0\leq t\leq 1\right\}.\nonumber
\end{flalign}
\noindent De plus, le lemme~\ref{lem:le902} implique l'indépendance de
\begin{flalign}
&\hspace{-1cm}\left\{(B\circ H_{m})(t)-4(2m+1)\frac{(2m-1)!}{2m((m-1)!)^2}(t(1-t))^m\right.\label{eq:mq}\\
&\hspace{-1cm}\left.\times\int_{0}^{1}B(H_{m}(s))\ud s: 0\leq t\leq 1\right\}\ \mbox{ et }\ \left\{(B\circ H_{m})(t):0\leq t\leq 1\right\}.\nonumber
\end{flalign}
\noindent Posons en accord avec~(\ref{eq:mc}),~(\ref{eq:mn}) et~(\ref{eq:mo}), pris avec $\ \displaystyle \mathfrak{X}=B\circ H_{m}\ $ et $\ d=1,\ $
{\small
\begin{flalign}
\hspace{-0,5cm}(B\circ H_{m})_{C}(t)&=B(H_{m}(t))\label{eq:mr}\\
\hspace{-0,5cm}&\quad-2(2m+1)C\ \frac{(2m-1)!}{m((m-1)!)^2}(t(1-t))^m\int_{0}^{1}B(H_{m}(s))\ud s,\nonumber
\end{flalign}
}pour $\displaystyle0\leq t\leq 1.$\vskip5pt
\noindent Compte tenu de~(\ref{eq:mq}), nous voyons que
\begin{flalign}
&\mathbb{E}\left((B\circ H_{m})_{C}(s)(B\circ H_{m})_{C}(t)\right)\label{eq:ms}\\
&=H_{m}(s)\land H_{m}(t)-H_{m}(s)H_{m}(t)\nonumber\\
&\quad+4(2m+1)(C^2-2C)\frac{((2m-1)!)^2}{(2m((m-1)!)^2)^2}(st(1-s)(1-t))^m,\nonumber
\end{flalign}
\noindent pour $\ \displaystyle 0\leq s,t\leq 1.$\vskip5pt
\end{rmq}
\section{Le Processus empirique de base}\label{ma-sectionc}
\subsection{Un théorème d'approximation forte}\label{ma-soussectionc}
\noindent Dans ce qui suit, nous construisons, sur le même
espace de probabilité $(\Omega,\mathcal{A},\mathbb{P})$, une suite
$\{(\zeta_{\ell,N}, \xi_{\ell,N}):\ell\geq1\}$, $N=1,2,\ldots$, de
répliques indépendantes et idem-distribuées [i.i.d.] d'un vecteur
aléatoire bivarié $\left(\zeta,\xi\right)$, où $\zeta$ et
$\xi$ sont des v.a. indépendantes de loi exponentielle de moyenne 1.
Cette loi est désignée par la suite par $\Gamma(1)$, et nous notons
ceci $\zeta\stackrel{d}{=}\xi\stackrel{d}{=}\Gamma(1)$. Nous
considérons ensuite une suite de vecteurs
$\left\{\left(Z_{k,N},Z_{k,N}^{\prime}\right): k\geq 0\right\}$,
$N=1,2\ldots$, telle que, pour tout $\displaystyle k\geq 0,\ $
\begin{equation}
Z_{k,N}:=\sum_{\ell=km+1}^{(k+1)m}\zeta_{\ell,N}\ \mbox{ et }\ Z_{k,N}^{\prime}:=\sum_{\ell=km+1}^{(k+1)m}\xi_{\ell,N}.\label{eq:iiaaaaaaaaaf}
\end{equation}
Cette définition implique que, pour tout $k\geq 0$,
$(Z_{k,N},Z_{k,N}^{\prime})\stackrel{d}{=}(Z,Z^{\prime})$, où
$Z:=Z_{0,N}$ et $Z^{\prime}:=Z_{0,N}^{\prime}$ sont indépendantes et
de même loi $ \Gamma(m,1)$. Cette dernière propriété est notée
$Z\stackrel{d}{=}Z^{\prime}\stackrel{d}{=}\Gamma(m,1)$. La fonction
de répartition [f.r.] jointe de $Z$ et $Z^{\prime}$ est donc donnée par
\begin{equation}
P\left(Z\leq x,Z^{\prime}\leq y\right)=G_{1}^{(m)}(x)G_{1}^{(m)}(y),\
\mbox{ pour }\ x,y\in\mathbb{R},
\end{equation}
où, pour $x\geq 0$, la f.r. de la loi
$\Gamma(m,1)$ est donnée par
\begin{equation}
G_{1}^{(m)}(x):=\frac{1}{\Gamma(m)}\int_{0}^{x}t^{m-1}e^{-t}\ \ud
t=1-e^{-x}\displaystyle\sum_{j=0}^{m-1}\frac{x^j}{j!}.\label{eq:iaaa}
\end{equation}
Désignons, respectivement, la f.r. de la loi $\beta_{m,m}$, pour
$0\leq x\leq 1$, par
\begin{eqnarray}\label{eq:iaab}
G_{3}^{(m)}(x)&:=& \frac{\Gamma(2m)}{(\Gamma(m))^2}\int_{0}^{x}(1-t)^{m-1}t^{m-1}\ud t\\
&=&\sum_{j=m}^{2m-1}\frac{(2m-1)!}{j!(2m-1-j)!}x^j(1-x)^{2m-1-j}\nonumber\\
&=&1-\sum_{j=0}^{m-1}\frac{(2m-1)!}{j!(2m-1-j)!}x^j(1-x)^{2m-1-j},
\nonumber
\end{eqnarray}
et la f.r. de la loi $\Gamma(2m,1)$, pour $x\geq 0$, par
\begin{eqnarray}\label{eq:iaac}
G_{2}^{(m)}(x)&:=&\frac{1}{\Gamma(2m)}\int_{0}^{x}t^{2m-1}e^{-t}\ud t=1-e^{-x}\sum_{j=0}^{2m-1}\frac{x^j}{j!}.
\end{eqnarray}
Désignons, respectivement, les fonctions de quantiles des lois $\beta_{m,m}$ et $\Gamma(2m,1)$, pour $0<u,v<1$, par
\begin{eqnarray*}
Q_{3}^{(m)}(u)&=&
\inf\{x\geq 0:G_{3}^{(m)}(x)\geq u\},\\
Q_{2}^{(m)}(v)&=&\inf\{x\geq 0:G_{2}^{(m)}(x)\geq v\}.
\end{eqnarray*}
Compte tenu de $(\ref{eq:iaab})$ et $(\ref{eq:iaac})$, on a, les relations réciproques
\begin{eqnarray}
G_{3}^{(m)}\left(Q_{3}^{(m)}(v)\right)&=&v,\ \mbox{ pour }\ 0< v<1,\label{eq:iaad}\\
G_{2}^{(m)}\left(Q_{2}^{(m)}(w)\right)&=&w,\ \mbox{ pour }\ 0< w<1,\label{eq:iaae}\\
\noalign{\hbox{et}}
Q_{3}^{(m)}\left(G_{3}^{(m)}(y)\right)&=&y,\ \mbox{ pour }\ 0<y<1,\label{eq:iaaf}\\
Q_{2}^{(m)}\left(G_{2}^{(m)}(y)\right)&=&y,\ \mbox{ pour }\ y\geq 0.\label{eq:iaag}
\end{eqnarray}
Désignons les densités des lois $\beta_{m,m}$ et $\Gamma(2m,1)$, par
{\small
\begin{eqnarray}
\hspace{0cm}g_{3}^{(m)}(y)&=&\frac{\ud}{\ud y}\ G_{3}^{(m)}(y)=\frac{1}{\beta(m,m)}(1-y)^{m-1}y^{m-1},\ \mbox{pour}\ 0<y< 1,\label{eq:iaaav}\\
\hspace{0cm}g_{2}^{(m)}(x)&=&\frac{\ud}{\ud x}\ G_{2}^{(m)}(x)=\frac{1}{\Gamma(2m)}x^{2m-1}e^{-x},\ \mbox{ pour }\ x>0.\label{eq:iaaaw}
\end{eqnarray}
}On déduit de $(\ref{eq:iaac})$, $(\ref{eq:iaae})$,
$(\ref{eq:iaad})$ et $(\ref{eq:iaag})$ que les densités de quantile $
q_{3}^{(m)}(v)=\frac{\ud}{\ud u}Q_{3}^{(m)}(v)$ et $q_{2}^{(m)}(w)=\frac{\ud}{\ud v}Q_{2}^{(m)}(w)$
sont continues en $v\in (0,1)$ et $w\in (0,1)$, et
telles que, pour $0< v,w<1$,
\begin{eqnarray}
q_{3}^{(m)}(v)&=&\frac{\ud}{\ud v}\ Q_{3}^{(m)}(v)=\frac{1}{g_{3}^{(m)}(Q_{3}^{(m)}(v))}\in(0,\infty),\label{eq:iaah}\\
q_{2}^{(m)}(w)&=&\frac{\ud}{\ud v}\
Q_{2}^{(m)}(w)=\frac{1}{g_{2}^{(m)}(Q_{2}^{(m)}(w))))}\in(0,\infty).\label{eq:iaai}
\end{eqnarray}
Notons aussi que,
\begin{eqnarray}
\beta(m,m)&=&\frac{\Gamma(m)^2}{\Gamma(2m)}=\frac{\{(m-1)!\}^2}{
(2m-1)!}=\frac{m!m!2m}{m^2(2m)!}
=\frac{2^{1-2m}m!}{m\left({\textstyle{\frac{1}{2}}}\right)_m}\,,\label{eq:iaaax}
\end{eqnarray}
relation dans laquelle nous avons fait usage du symbole de
Pochhammer
\[(a)_m=\frac{\Gamma(a+m)}{\Gamma(a)}=a\times (a+1)\times\ldots\times(a+m-1).\]
Comme la loi $\ \displaystyle\beta_{m,m}\ $ est symétrique par rapport à  $\ \displaystyle 1/2$,
on a,
\begin{equation}
g_{3}^{(m)}(1-y)=g_{3}^{(m)}(y),\ \hbox{pour}\quad 0<y<1,\label{eq:iaaay}
\end{equation}
et donc, pour $0<y<1$, $G_{3}^{(m)}(y)=1-G_{3}^{(m)}(1-y)$, et pour $0<v<1$, $Q_{3}^{(m)}(v)=1-Q_{3}^{(m)}(1-v)$,
de sorte que
\begin{equation}
q_{3}^{(m)}(v)=q_{3}^{(m)}(1-v),\ \hbox{ pour }\
0<v<1.\label{eq:iaaaz}
\end{equation}
La proposition suivante décrivant le comportement asymptotique de
$q_{3}^{(m)}(v)$ et $Q_{3}^{(m)}(v)$ lorsque $v\downarrow 0$ ou $v\uparrow 1$,
ainsi que celui de $q_{2}^{(m)}(w)$ et $Q_{2}^{(m)}(w)$ lorsque
$w\downarrow 0$ ou $w\uparrow 1$, sera utile par la suite.
\begin{prop}\label{propos:propo250}Nous avons
{\small
\begin{flalign}
&q_{3}^{(m)}(v)=\frac{1}{m}
\left\{\frac{\left(\frac{1}{2}\right)_{m}}{2^{1-2m}m!}\right\}^{-1/m}
v^{1/m-1}\left(1+o(1)\right),\; v\downarrow 0,\label{eq:iaaj0}\\
&q_{3}^{(m)}(v)=\frac{1}{m}\left\{\frac{\left(\frac{1}{2}\right)_{m}}
{2^{1-2m}m!}\right\}^{-1/m}
(1-v)^{1/m-1}\left(1+o(1)\right),\; v\uparrow 1,\qquad\qquad\label{eq:iaaj}\\
&q_{2}^{(m)}(w)=\frac{1+o(1)}{2m}
\left\{\Gamma(2m+1)\right\}^{1/(2m)}w^{1/(2m)-1},\; w\downarrow
0\label{eq:iaak0}\\
&q_{2}^{(m)}(w)=\frac{1+o(1)}{1-w},\; w\uparrow
1,\label{eq:iaak}\end{flalign}}ainsi que
{\small
\begin{flalign}
&Q_{3}^{(m)}(v)=\left\{\frac{\left(\frac{1}{2}\right)_{m}}
{2^{1-2m}m!}\right\}^{-1/m}(1+o(1))v^{1/m},\;
v\downarrow 0,\label{eq:iaal0}\qquad\qquad\\
&Q_{3}^{(m)}(v)=1-\left\{\frac{\left(\frac{1}{2}\right)_{m}}{2^{1-2m}m!}
\right\}^{-1/m}(1+o(1))\left(1-v\right)^{1/m},\; v\uparrow
1.\label{eq:iaal}\qquad\qquad\\
&Q_{2}^{(m)}(w)=\left\{\Gamma(2m+1)\right\}^{1/(2m)}\left(1+o(1)\right)w^{1/(2m)},
\; w\downarrow 0,\label{eq:iaam0}\\
&Q_{2}^{(m)}(w)=\int_{0}^{w}q_{2}^{(m)}(s)\ \ud
s=\left(1+o(1)\right)\left\{-\log(1-w)\right\},\; w\uparrow
1.\label{eq:iaam}
\end{flalign}}
\end{prop}
\noindent{\bf Démonstration.} Par~(\ref{eq:iaaav}), on obtient que, lorsque $y\downarrow
 0$,
\begin{eqnarray}
\hspace{-1cm}&&g_{3}^{(m)}(y)=\frac{1}{\beta(m,m)}\Big\{y^{m-1}-\binom{m-1}{1}y^{m}+\ldots\nonumber\\
\hspace{-1cm}&&\quad+(-1)^{m-1}\binom{m-1}{m-1}y^{2m-2}\Big\}=\left\{\frac{m\left({\textstyle{\frac{1}{2}}}\right)_m}{2^{1-2m}m!}
\right\}y^{m-1}(1+o(1)).\label{eq:iaaaaa}
\end{eqnarray}
Par conséquent, on a, lorsque $y\downarrow0$,
\[G_{3}^{(m)}(y)
=\left\{\frac{\left({\textstyle{\frac{1}{2}}}\right)_m}{2^{1-2m}m!}
\right\}y^{m}(1+o(1)).\] En posant $\ \displaystyle
y=Q_{3}^{(m)}(u)\ $ dans cette relation , nous voyons que, lorsque
$u\downarrow 0$, $Q_{3}^{(m)}(u)\downarrow 0$, et
\[G_{3}^{(m)}(Q_{3}^{(m)}(u))=u=\left\{\frac{
\left({\textstyle{\frac{1}{2}}}\right)_m}{2^{1-2m}m!}
\right\}\left(Q_{3}^{(m)}(u)\right)^{m}(1+o(1)).\] On en déduit
(\ref{eq:iaal0}). En appliquant à (\ref{eq:iaal0}), les
relations~(\ref{eq:iaah}) et~(\ref{eq:iaaaaa}), on obtient que,
lorsque $ u\downarrow 0$, \begin{eqnarray*}q_{3}^{(m)}(u)
&=&\frac{1}{g_{3}^{(m)}(Q_{3}^{(m)}(u))}\\
&=&\left\{\frac{2^{1-2m}m!}
{m\left({\textstyle{\frac{1}{2}}}\right)_m}
\right\}\left(\left\{\frac{\left({\textstyle{\frac{1}{2}}}\right)_m}{2^{1-2m}m!}
\right\}^{-1/m}u^{1/m}\right)^{1-m}(1+o(1)).\end{eqnarray*} On en
déduit (\ref{eq:iaaj0}). En posant $ u=1-v$ dans (\ref{eq:iaaj0}),
on obtient que, lorsque $v\to 1$,
\[q_{3}^{(m)}(1-v)=\frac{1}{m}\left\{
\frac{\left({\textstyle{\frac{1}{2}}}\right)_m}{2^{1-2m}m!}
\right\}^{-1/m}(1-v)^{1/m-1}(1+o(1)).\] Ceci, via~(\ref{eq:iaaaz}),
établit (\ref{eq:iaaj}).  Par ailleurs, en posant $u=1-v$  dans
(\ref{eq:iaal0}), on obtient que, lorsque $ v\to 1$,
\[Q_{3}^{(m)}(1-v)=\left\{\frac{\left({\textstyle{\frac{1}{2}}}\right)_m}{2^{1-2m}m!}
\right\}^{-1/m}(1-v)^{1/m}(1+o(1)).\] Ceci, via~(\ref{eq:iaaaz}),
établit(\ref{eq:iaal}). Par~(\ref{eq:iaae}) et compte tenu de
l'expression de la densité $g_{2}^{(m)}$, on obtient que, lorsque
$s\to\infty$,
\begin{displaymath}
1-G_{2}^{(m)}(s)=\left\{1+\sum_{j=1}^{2m-1}\frac{(2m-1)!}{(2m-1-j)!s^j}\right\}g_{2}^{(m)}(s)=\left(1+o(1)\right)g_{2}^{(m)}(s).
\end{displaymath}En posant $s=Q_{2}^{(m)}(t)$ dans cette relation, nous voyons que,
lorsque $t\uparrow 1$, $Q_{2}^{(m)}(t)\uparrow\infty$, et $
1-G_{2}^{(m)}(Q_{2}^{(m)}(t))=1-t=\left(1+o(1)\right)g_{2}^{(m)}
\left(Q_{2}^{(m)}(t)\right)$. Compte tenu de~(\ref{eq:iaai}), ceci
suffit pour justifier (\ref{eq:iaak}). Par (\ref{eq:iaae}) et
l'expression de la densité $g_{2}^{(m)}$, un développement de Taylor
de $G_{2}^{(m)}(s)$ dans un voisinage à droite de $0$ montre que,
lorsque $s\downarrow0$,
\begin{eqnarray*}
G_{2}^{(m)}(s)&=&1-e^{-s}\left(1+\sum_{j=1}^{2m-1}\frac{s^j}{j!}\right)\\
&=&1-\left(1+\sum_{k=1}^{2m}\frac{(-1)^{j}s^j}{j!}+o\left(s^{2m}\right)\right)
\left(1+\sum_{j=1}^{2m-1}\frac{s^j}{j!}\right)\\
&=&\left(1+o(1)\right)\frac{s^{2m}}{\Gamma(2m+1)},
\end{eqnarray*}
et
\begin{equation}
g_{2}^{(m)}(s)=\frac{1}{\Gamma(2m)}\left(1+o(1)\right)s^{2m-1}.\label{eq:iaaaab}
\end{equation}Ceci implique que, lorsque $t\downarrow 0$,
$Q_{2}^{(m)}(t)\downarrow 0$ et \[
G_{2}^{(m)}\left(Q_{2}^{(m)}(t)\right)=t=\left(1+o(1)\right)
\frac{\left(Q_{2}^{(m)}(t)\right)^{2m}}{\Gamma(2m+1)}.\] On en
déduit (\ref{eq:iaam0}). De plus, en appliquant à (\ref{eq:iaam0}),
les relations~(\ref{eq:iaai}) et~(\ref{eq:iaaaab}), on obtient que,
lorsque $u\downarrow 0$,
\[q_{2}^{(m)}(t)=\frac{1}{g_{2}^{(m)}(Q_{2}^{(m)}(t))}=(1+o(1))\Gamma(2m)
\Big\{\Big\{\Gamma(2m+1)\Big\}^{1/(2m)}t^{1/(2m)}\Big\}^{1-2m}.\] On
en déduit (\ref{eq:iaak0}). Finalement, nous voyons que, lorsque
$t\uparrow1$,

\begin{flalign*}
\hspace{-0,5cm}-\log(1-t)&=-\log\left\{1-G_{2}^{(m)}\left(Q_{2}^{(m)}(t)\right)\right\}\\
\hspace{-0,5cm}&=-\log\left\{\left(1+\sum_{j=1}^{2m-1}\frac{\left(Q_{2}^{(m)}(t)\right)^{j}}{j!}\right)e^{-Q_{2}^{(m)}(t)}\right\}\\
\hspace{-0,5cm}&=Q_{2}^{(m)}(t)-\log\left(1+\sum_{j=1}^{2m-1}\frac{\left(Q_{2}^{(m)}(t)\right)^j}{j!}\right)=\left(1+o(1)\right)Q_{2}^{(m)}(t),
\end{flalign*}

\noindent qui complète la démonstration
de~(\ref{eq:iaam}).$\Box$\vskip5pt
\noindent En se rappelant de la définition~(\ref{eq:iiaaaaaaaaaf}) de la suite i.i.d. de vecteurs
$\left\{\left(Z_{k,N},Z^{\prime}_{k,N}\right):\right.$ $\left.k\geq 0,\ N\geq 1\right\},\ $ ainsi que $\ Z\ $ et $\ Z^{\prime},\ $
nous introduisons la paire aléatoire $\ \displaystyle (R,T)\ $ et une suite i.i.d. $
\left\{\left(R_{k,N},T_{k,N}\right):k\geq 0,\ N\geq 1\right\},\ $ de répliques indépendantes de $\ \displaystyle (R,T),\ $
en posant
\begin{equation}
R=\frac{Z}{Z+Z^{\prime}},\ T=Z+Z^{\prime},\ \mbox{ et}\label{eq:iaan}
\end{equation}
\begin{equation}
R_{k,N}=\frac{Z_{k,N}}{Z_{k,N}+Z_{k,N}^{\prime}},\;
T_{k,N}=Z_{k,N}+Z_{k,N}^{\prime},\ \mbox{ pour }\ k\geq 0\ \mbox{ et
}\ N\geq1.\label{eq:iaao}
\end{equation}
\noindent Les variables $R,T$ dans~(\ref{eq:iaan}) (de même que $
R_{k,N},T_{k,N}$ dans~(\ref{eq:iaao})) sont indépendantes,
respectivement, de lois $\beta_{m,m}$ et $\Gamma(2m,1)$. Cette
relation est notée $R\stackrel{d}{=}\beta_{m,m}$ et
$T\stackrel{d}{=}\Gamma(2m,1)$. La loi jointe de $R,T$ est donc
donnée par
\begin{equation}
P\left(R\leq v,\ T\leq y\right)=G_{3}^{(m)}(v)\ G_{2}^{(m)}(y),\ \mbox{ pour }\ v,y\in\mathbb{R},\label{eq:iaap}
\end{equation}
où $G_{3}^{(m)}(.)$ et $G_{2}^{(m)}(.)$ sont définies
par~(\ref{eq:iaab}) et~(\ref{eq:iaac}). Compte tenu de la
transformation de quantiles (voir le théorème 1, pp. 3-4 de
Shorack et Wellner\cite{MR838963}), nous posons, pour $k\geq 0$
et $N\geq 1$,
\begin{eqnarray}
V&=&G_{3}^{(m)}(R)\ \mbox{ et }\
V_{k,N}=G_{3}^{(m)}\left(R_{k,N}\right),\label{eq:iaaq}\\
\noalign{\hbox{et}} W&=&G_{2}^{(m)}(T)\ \mbox{ et }\
W_{k,N}=G_{2}^{(m)}\left(T_{k,N}\right).\label{eq:iaar}
\end{eqnarray}
\noindent $\displaystyle
\left\{\left(V_{k,N},W_{k,N}\right):k\geq 0,\ N\geq1\right\}\ $
définit une suite i.i.d. de répliques indépendantes de $(V,W)$, où
$V\stackrel{d}{=}\mathcal{U}(0,1)$ et
$W\stackrel{d}{=}\mathcal{U}(0,1)$ sont de loi uniforme $\mathcal{U}(0,1)$.
~(\ref{eq:iaad}),~(\ref{eq:iaaf}),~(\ref{eq:iaae}),~(\ref{eq:iaag})
et~(\ref{eq:iaao}) impliquent que, pour tout $\
\displaystyle0\leq k\leq N$,
\begin{flalign}
&Z_{k,N}=R_{k,N}\ T_{k,N}=Q_{3}^{(m)}\left(V_{k,N}\right)Q_{2}^{(m)}\left(W_{k,N}\right)\ \mbox{ et }\label{eq:iaas}\\
&Z_{k,N}^{\prime}=\left(1-R_{k,N}\right)T_{k,N}=\left(1-Q_{3}^{(m)}\left(V_{k,N}\right)\right)Q_{2}^{(m)}\left(W_{k,N}\right).\nonumber
\end{flalign}
Pour $\ \displaystyle N\geq 1,\ $ la mesure empirique basée
sur $\ \displaystyle
\left\{\left(V_{k,N},W_{k,N}\right)\right.:\left.0\leq k\leq
N\right\}\ $ est notée
\begin{equation}
\lambda_{N}(.)=\frac{1}{N+1}\sum_{k=0}^
{N}\delta_{\left(V_{k,N},W_{k,N}\right)}(.),\label{eq:iaat}
\end{equation}
 où $\ \displaystyle \delta_{z}(.)\ $ désigne la mesure de Dirac en $\ \displaystyle z\in\mathbb{R}^2.\ $ En dénotant
par $\ \displaystyle \lambda(.)\ $ la mesure de Lebesgue sur $\
\displaystyle[0,1]^2,\ $ le processus empirique uniforme indexé par
des ensembles est défini par
\begin{equation}
\alpha_{N}(A):=(N+1)^{1/2}
\left(\lambda_{N}(A)-\lambda(A)\right),\label{eq:iaau}
\end{equation}
 pour chaque borélien $A\subseteq  [0,1]^2$. Lorsque
 $A$ est un produit
d'intervalles, nous définissons la version continue à droite de la
f.r. empirique bivariée basée sur
$\left\{\left(V_{k,N},W_{k,N}\right):0\leq k\leq N\right\}$, en
posant, pour $\ \displaystyle0\leq v,w\leq 1,\ $
\begin{flalign}
\mathbf{U}_{N}(v,w)&=\lambda_{N}\left([0,v]\times[0,w]\right)\label{eq:iaav}\\
&=\frac{1}
{N+1}\#\left\{V_{k,N}\leq v,W_{k,N}\leq w:0\leq k\leq N\right\},\nonumber
\end{flalign}
$\# E $ désignant, ici, le nombre d'éléments de $E$. Pour $0\leq
v,w\leq 1$, on note
\begin{equation}
\alpha_{N}(v,w)=\alpha_{N}\left([0,v]\times [0,w]\right)=(N+1)^{1/2}\left(\mathbf{U}_{N}
(v,w)-vw\right),\label{eq:iaaw}
\end{equation}
la version continue à droite du processus empirique uniforme
engendré par  $\left\{\left(V_{k,N},W_{k,N}\right):0\leq k\leq
N\right\}$. La version continue à gauche du processus empirique
engendré par $\left\{\left(1-V_{k,N},1-W_{k,N}\right):0\leq k\leq
N\right\}$ sera notée par analogie, $\ \displaystyle\left\{
\alpha_{N}^{\ast}(v,w):0\leq v,w\leq 1\right\}$. Comme $\
\displaystyle \alpha_{N}(1,1)=0,\ $ ce processus vérifie  les
relations, pour $\ \displaystyle 0\leq u,v\leq 1,\ $
\begin{eqnarray}
\alpha_{N}^{\ast}(1-v,1-w)&
=&\alpha_{N}^{\ast}\left((v,1]\times (w,1]\right)\label{eq:iaax}\\
&=&\alpha_{N}(v,w)-\alpha_{N}(v,1)
-\alpha_{N}(1,w).\nonumber
\end{eqnarray}
On a l'identité en loi
\begin{displaymath}
\left\{\alpha_{N}^{\ast}(s,t):0\leq s,t\leq 1\right\}\stackrel{d}{=}\left\{\lim_{v\uparrow s, w\uparrow t}\alpha_{N}(v,w):0\leq s,t\leq
1\right\}.
\end{displaymath}
On notera que les versions continues à droite des processus
empiriques marginaux
\begin{eqnarray}
\alpha_{N:1}(v)&:=&\alpha_{N}(v,1),\ \mbox{ pour}\ 0\leq v\leq
1,\label{eq:iaay}
\\
\alpha_{N:2}(w)&:=&\alpha_{N}(1,w),\ \mbox{ pour}\ 0\leq w\leq
1,\label{eq:iaaz}
\end{eqnarray}
engendrés, respectivement par $\ \displaystyle \left\{V_{k,N}:0\leq
k\leq N\right\}\ $ et $\ \displaystyle \left\{W_{k,N}:0\leq k\leq
N\right\}\ $, sont des processus empiriques uniformes indépendants,
sur $[0,1]$.\vskip5pt
\noindent Nous allons étudier le comportement limite de $\{
\alpha_{N:1}(v):0\leq v\leq 1\}$, processus empirique engendré, par
$\left\{V_{k,N}:0\leq k\leq N\right\}$, et des statistiques
\begin{eqnarray}
\Delta_{N}&:=&\overline{Z}_{N}
-\overline{Z}_{N}^{\prime}=\frac{1}{N+1} \sum_{k=0}^{N}\left\{Z_{k,N}
-Z_{k,N}^{\prime}\right\},\label{eq:iaaaa}\\
\noalign{\hbox{et}} \Theta_{N}&:=&\overline{Z}_{N}
+\overline{Z}_{N}^{\prime}-2m=\frac{1}{N+1}
\sum_{k=0}^{N}\left\{Z_{k,N}
+Z_{k,N}^{\prime}-2m\right\}\label{eq:iaaab}\\
&=&\frac{1}{N+1} \sum_{k=0}^{N}\left\{T_{k,N}
-2m\right\}.\nonumber\end{eqnarray} Nous établissons le
théorème suivant qui généralise à $\ \displaystyle m\geq 1\ $ le théorème 3.1 de de
Deheuvels et Derzko~\cite{MR2728435}.\vskip5pt
\begin{teo}
\label{thm:teo11111}
\noindent Sous une version appropriée de
$(\Omega,\mathcal{A},\mathbb{P})$, il existe une suite $\
\left\{\left(B_{N}(.),\phi_{N},\psi_{N}\right):N\geq 1\right\}\ $\
vérifiant les propriétés suivantes.\vskip5pt
(i) Pour chaque $\ N\geq 1,\ \left\{B_{N}(v):0\leq v\leq 1\right\}\ $ est un pont Brownien et $\phi_{N}\stackrel{d}{=}\psi_{n}\stackrel{d}{=}N(0,1)\ $ sont deux v.a. normales standard.\vskip5pt
(ii) Pour chaque $N\geq 1$, $\{B_{N}(v):0\leq v\leq 1\}$, et les
variables aléatoires $\phi_{N}$ et $\psi_{N}$ sont
indépendants.\vskip5pt
(iii) Nous avons, lorsque $\ N\to\infty,\ $
\begin{eqnarray}
&&\hspace{-0,5cm}\left\vert\left\vert\alpha_{N;1}-B_{N}\right\vert\right\vert=\sup_{0\leq
v\leq 1}\left\vert
\alpha_{N;1}(v)-B_{N}(v)\right\vert=O_{\mathbb{P}}\left(\frac{\left(\log
N\right)^2}{\sqrt{N}}\right),
\label{eq:iaaaaax}\\
&&\hspace{-0,5cm}\left\vert
(N+1)^{1/2}\Delta_{N}-\phi_{N}\sqrt{\frac{2m}{2m+1}}
+4m\int_{0}^{1}B_{N}(G_{3}^{(m)}(x))\ \ud x\right\vert\label{eq:iaaaaay}\\
&&\qquad\qquad=O_{\mathbb{P}}\left(\frac{\left(\log N\right)^3}{\sqrt{N}}\right),\nonumber\\
\noalign{\hbox{et}}&&\hspace{-0,5cm}\left\vert
(N+1)^{1/2}\Theta_{N}-\psi_{N}\sqrt{2m}\right\vert
=O_{\mathbb{P}}\left(\frac{\left(\log
N\right)^3}{\sqrt{N}}\right).\label{eq:iaaaaaz}
\end{eqnarray}
\end{teo}
\noindent{\bf Démonstration.} La démonstration du théorème~\ref{thm:teo11111} est reportée au paragraphe~\ref{ma-soussectionk}.
Les ingrédients de base nécessaires pour cette preuve sont donnés dans les paragraphes~\ref{ma-soussectiond},~\ref{ma-soussectione}
et dans le paragraphe~\ref{ma-sectiond} ci-dessous.
\subsection{Décompositions du processus empirique}\label{ma-soussectiond}
\noindent La première étape dans la description de
$\alpha_{N:1}(.)$, $\Delta_{N}$ et $\Theta_{N}$, définis par
(\ref{eq:iaay}),~(\ref{eq:iaaaa}) et~(\ref{eq:iaaab}), est fournie
par les propositions~\ref{propos:propo251} et~\ref{propos:propo252}
suivantes, qui généralisent au cas $\ \displaystyle m\geq 1$ les
propositions 3.1 et 3.2 de Deheuvels et Derzko~\cite{MR2728435}.
Nous commençons par un lemme qui jouera un rôle important
par la suite. Rappelons la définition~(\ref{eq:iaax}) de $\
\displaystyle \alpha_{N}^{\ast}(.,.).\ $
\begin{lemme}
\label{lem:le250}
 Nous avons l'égalité
\begin{equation}
\int_{0}^{1}\alpha_{N:1}(v)\ud v=\int_{0}^{1}\alpha_{N}(v,1)\ud
v=-\int_{0}^{1}\alpha_{N}^{\ast}(z,1)\ud z.\label{eq:iaaac}
\end{equation}
\end{lemme}
\noindent{\bf Démonstration.}
 Nous déduisons de~(\ref{eq:iaax}) que $\ \displaystyle \alpha_{N:1}(v)
 =\alpha_{N}(v,1)
=-\alpha_{N}^{\ast}(1-v,1),\ $ pour $\ \displaystyle 0\leq v\leq 1.\
$ La conclusion suit après le changement de variables $\
\displaystyle v=1-z\ $ dans~(\ref{eq:iaaac}).$\Box$
\begin{prop}\label{propos:propo251}Nous avons les identités
{\small
\begin{eqnarray}
\hspace{0cm}&&(N+1)^{1/2}\Delta_{N}=\int\int_{[0,1]^2}\left(2\
Q_{3}^{(m)}(s)-1\right)
\ Q_{2}^{(m)}(t)\ \alpha_{N}(\ud s,\ud t)\label{eq:iaaad}\\
\hspace{0cm}&&=2\int\int_{[0,1]^2}q_{3}^{(m)}(v)\
q_{2}^{(m)}(w)\left\{\alpha_{N}(v,w)
-\alpha_{N}(v,1)-\alpha_{N}(1,w)\right\}\ \ud v\ \ud w\nonumber\\
\hspace{0cm}&&\quad+\int_{0}^{1}q_{2}^{(m)}(w)\ \alpha_{N}(1,w)\ \ud w\nonumber\\
\hspace{0cm}&&=2\int\int_{[0,1]^2}q_{3}^{(m)}(1-v)\ q_{2}^{(m)}(1-w)\ \alpha_{N}^{\ast}(v,w)\ \ud v\ \ud w\nonumber\\
\hspace{0cm}&&\quad-\int_{0}^{1}q_{2}^{(m)}(1-w)\
\alpha_{N}^{\ast}(1,w)\ \ud w,\nonumber\\
\noalign{\hbox{et}}
\hspace{0cm}&&(N+1)^{1/2}\ \Theta_{N}=\int_{0}^{1}Q_{2}^{(m)}(t)\ \alpha_{N}(1,\ud t)\label{eq:iaaae}\\
&&=-\int_{0}^{1}q_{2}^{(m)}(w)\ \alpha_{N}(1,w)\ \ud
w=\int_{0}^{1}q_{2}^{(m)}(1-w)\ \alpha_{N}^{\ast}(1,w)\ \ud
w.\nonumber
\end{eqnarray}}
\end{prop}
\noindent{\bf Démonstration.} Par~(\ref{eq:iaas}),~(\ref{eq:iaav})
et~(\ref{eq:iaaaa}), on a,
\begin{equation}
\Delta_{N}=\int\int_{[0,1]^2}\left(2\ Q_{3}^{(m)}(s)-1\right)
\ Q_{2}^{(m)}(t)\ \mathbf{U}_{N}(\ud s,\ud t)=:\Delta_{N,1}+\Delta_{N,2},
\label{eq:iaaaf}
\end{equation}
où
\begin{equation}
\Delta_{N,1}:=2\int\int_{[0,1]^2}Q_{3}^{(m)}(s)\ Q_{2}^{(m)}(t)\ \mathbf{U}_{N}(\ud s,\ud t),\label{eq:iaaag}
\end{equation}
et, compte tenu de~(\ref{eq:iaas}) et~(\ref{eq:iaaab}),
\begin{equation}
\Delta_{N,2}:=-\int_{0}^{1}Q_{2}^{(m)}(t)\ \mathbf{U}_{N}(1,\ud t)=-\left\{\Theta_{N}+2m\right\}.\label{eq:iaaah}
\end{equation}
Compte tenu de~(\ref{eq:iaah}),~(\ref{eq:iaai}),~(\ref{eq:iaal0})
et~(\ref{eq:iaam0}), nous appliquons le théorème de Fubini
dans~(\ref{eq:iaaag}) pour obtenir les relations
{\footnotesize
\begin{flalign}
\hspace{0cm}\Delta_{N,1}&=2\int\int_{[0,1]^2}Q_{3}^{(m)}(s)\ Q_{2}^{(m)}(t)\ \mathbf{U}_{N}(\ud s,\ud t)\label{eq:iaaai}\\
\hspace{0cm}&=2\int\int_{[0,1]^2}\left\{\int\int_{[0,1]^2}{\rm
1\!I}_{[0,s)}(v)q_{3}^{(m)}(v){\rm
1\!I}_{[0,t)}(w)q_{2}^{(m)}(w)\ \ud v\ud w\right\}\mathbf{U}_{N}(\ud s,\ud t)\nonumber\\
\hspace{0cm}&=2\int\int_{[0,1]^2}q_{3}^{(m)}(v)\ q_{2}^{(m)}(w)\left\{\int\int_{[0,1]^2}{\rm
1\!I}_{[0,s)}(v){\rm
1\!I}_{[0,t)}(w)\ \mathbf{U}_{N}(\ud s,\ud t)\right\}\ud v\ \ud w\nonumber\\
\hspace{0cm}&=2\int\int_{[0,1]^2}q_{3}^{(m)}(v)\ q_{2}^{(m)}(w)\left\{{\mathbf{\lambda}}_{N}\left((v,1]\times(w,1]\right)\right\}\ud v\ \ud w\nonumber\\
\hspace{0cm}&=2\int\int_{[0,1]^2}q_{3}^{(m)}(v)\ q_{2}^{(m)}(w)\left\{\mathbf{U}_{N}(1,1)
-\mathbf{U}_{N}(v,1)\right.\nonumber\\
\hspace{0cm}&\left.\hspace{5,2cm}-\mathbf{U}_{N}(1,w)
+\mathbf{U}_{N}(v,w)\right\}\ud v\ \ud w.\nonumber
\end{flalign}}
Les hypothèses du théorème de Fubini exigent que la fonction $\
\displaystyle {\rm 1\!I}_{[0,s)}(v)q_{3}^{(m)}(v)\times$ ${\rm
1\!I}_{[0,t)}(w)q_{2}^{(m)}(w)\ $ de $\ \displaystyle
s,t,v,w\in[0,1]\ $ soit intégrable sur $\ [0,1]^4,\ $ par rapport à
$\ \displaystyle \ud v\ud w \ \mathbf{U}_{N}(\ud v,\ud w).\ $ Cette
propriété découle du fait que la mesure empirique $\ \displaystyle
{\mathbf{\lambda}}_{N}={\mathbf{U}}_{N}(\ud v,\ud w)\ $ s'annule sur l'intervalle $
(\max_{0\leq k\leq N}V_{k,N},1]\times(\max_{0\leq k\leq
N}W_{k,N},1]$, presque sûrement de mesure de Lebesgue positive. Les
mêmes arguments, compte tenu de~(\ref{eq:iaaah}), donnent les
relations
\begin{flalign}
\hspace{0cm}\Delta_{N,2}&=-\int_{0}^{1}Q_{2}^{(m)}(t)\ \mathbf{U}_{N}(1,\ud t)\label{eq:iaaaj}\\
\hspace{0cm}&=-\int_{0}^{1}\left\{\int_{0}^{1}{\rm
1\!I}_{[0,t)}(w)q_{2}^{(m)}(w)\ \ud v\right\}\mathbf{U}_{N}(1,\ud t)\nonumber\\
\hspace{0cm}&=-\int_{0}^{1}q_{2}^{(m)}(w)\left\{\int_{0}^{1}{\rm
1\!I}_{[0,t)}(w)\right\}\mathbf{U}_{N}(1,\ud t)\ud w\nonumber\\
\hspace{0cm}&=-\int_{0}^{1}q_{2}^{(m)}(w)\left\{\mathbf{U}_{N}(1,1)-\mathbf{U}_{N}
(1,w)\right\}\ud w=-\left\{\Theta_{N}+2m\right\}.\nonumber
\end{flalign}
Le théorème de Fubini s'applique, du fait que la mesure empirique $\
\displaystyle {\mathbf{U}}_{N}(1,\ud t)\ $ est nulle sur l'intervalle $
\left(\max_{0\leq k\leq N}W_{k,N},1\right]$, presque sûrement de
mesure de Lebesgue positive.\vskip5pt

\noindent En se rappelant de~(\ref{eq:iaaf}) $\left(\right.$resp.~(\ref{eq:iaag})$\left.\right)$,~(\ref{eq:iaaav})   $\left(\right.$resp.~(\ref{eq:iaaaw})$\left.\right)$ et de~(\ref{eq:iaaj})
$\left(\right.$resp.~(\ref{eq:iaak})$\left.\right)$, nous faisons usage du changement de variable $\ \displaystyle w=1-v\ $ $\left(\right.$resp.
$\ \displaystyle v=1-w $ $\left.\right)$ et $\ \displaystyle v=G_{3}^{(m)}(x)=1-\sum_{j=0}^{m-1}\frac{(2m-1)!}{j!(2m-1-j)!}x^j(1-x)^{2m-1-j}\ $
$\bigg($resp. $\ \displaystyle w=G_{2}^{(m)}(x)=1-e^{-x}.\displaystyle\sum_{j=0}^{2m-1}\frac{x^j}{j!} $ $\bigg)$, pour obtenir les égalités $\
\displaystyle q_{3}^{(m)}(v)=q_{3}^{(m)}(G_{3}^{(m)}(x))=1/g_{3}^{(m)}(x),\ \ud v=g_{3}^{(m)}(x)\ \ud x\ \mbox{ et }\ q_{3}^{(m)}(v)\ \ud v=\ud
x\ $ $\left(\right.$ resp. $\ \displaystyle q_{2}^{(m)}(w)=q_{2}^{(m)}(G_{2}^{(m)}(x))=1/g_{2}^{(m)}(x),\ \ud w=g_{2}^{(m)}(x)\ \ud x\ \mbox{ et
}\ q_{2}^{(m)}(w)\ \ud v=\ud x $ $\left.\right)$, pour $\ \displaystyle 0<v<1\ $ $\left(\right.$ resp. $\ \displaystyle 0<w<1$ $\left.\right)$.
Ceci, à son tour, donne, via~(\ref{eq:iaab}) $\left(\right.$resp.~(\ref{eq:iaac})$\left.\right)$,

{\small
\begin{flalign}
\int_{0}^{1}w\ q_{3}^{(m)}(1-w)\ \ud w&=\int_{0}^{1}(1-v)\ q_{3}^{(m)}(v)\ \ud v=\int_{0}^{1}\left(1-G_{3}^{(m)}(x)\right)\ \ud x\label{eq:iaaak}\\
&=\sum_{j=0}^{m-1}\frac{(2m-1)!}{j!(2m-1-j)!}\int_{0}^{1}x^{j+1-1}(1-x)^{2m-1-j}\ \ud x\nonumber\\
&=\sum_{j=0}^{m-1}\frac{(2m-1)!}{j!(2m-1-j)!}\frac{j!(2m-1-j)!}{(2m)!}=\frac{1}{2}.\nonumber
\end{flalign}
}

{\footnotesize
\begin{flalign}
\hspace{0cm}\Bigg(\mbox{resp.}\ \int_{0}^{1}v\ q_{2}^{(m)}(1-v)\ \ud v&=\int_{0}^{1}(1-w)\ q_{2}^{(m)}(w)\ \ud
v=\int_{0}^{1}\left(1-G_{2}^{(m)}(x)\right)\ \ud x\label{eq:iaaal}\\
&=\int_{0}^{\infty}e^{-x}.\displaystyle\sum_{j=0}^{2m-1}\frac{x^j}{j!}=\displaystyle\sum_{j=0}^{2m-1}\frac{1}{j!}\int_{0}^{\infty}e^{-x}x^j\ \ud
x\nonumber\\
&\left.=\displaystyle\sum_{j=0}^{2m-1}\frac{1}{j!}\ \Gamma(j+1)=2m\right).\nonumber
\end{flalign}}

\noindent En combinant~(\ref{eq:iaaai}) avec~(\ref{eq:iaaaj}),~(\ref{eq:iaaak}) et~(\ref{eq:iaaal}), nous obtenons les relations,

\begin{flalign}
\mathbb{E}\left(\Delta_{N,1}\right)&=2\int\int_{[0,1]^2}q_{3}^{(m)}(v)q_{2}^{(m)}(w)\left\{(1-v)(1-w)\right\}\ \ud v\ud w\label{eq:iaaam}\\
&=\int_{0}^{1}(1-w)\ q_{2}^{(m)}(w)\ \ud w=-\mathbb{E}\left(\Delta_{N,2}\right)=2m.\nonumber
\end{flalign}

\noindent En combinant~(\ref{eq:iaaai}),~(\ref{eq:iaaaj}) et~(\ref{eq:iaaam}), avec la définition~(\ref{eq:iaaw}) de $\ {\mathbf{\alpha}}_{N}(v,w),\ $ nous
obtenons les 2 premières égalités de~(\ref{eq:iaaad}).
Pour prouver la dernière égalité de~(\ref{eq:iaaad}), nous rappelons d'abord par~(\ref{eq:iaax}) que $\ \displaystyle{\mathbf{\alpha}}_{N}^{\ast}(1-v,1-w)
={\mathbf{\alpha}}_{N}(v,w)-{\mathbf{\alpha}}_{N}(v,1)
-{\mathbf{\alpha}}_{N}(1,w),\ $ qui pour $\ \displaystyle v=0,\ $ donne $\ \displaystyle
{\mathbf{\alpha}}_{N}^{\ast}(1,1-w)=-{\mathbf{\alpha}}_{N}(1,w).\ $
Ceci à son tour, montre que

\begin{flalign*}
(N+1)^{1/2}\Delta_{N}&=2\int\int_{[0,1]^2}q_{3}^{(m)}(v)\ q_{2}^{(m)}(w)\ {\mathbf{\alpha}}_{N}^{\ast}(1-v,1-w)\ \ud v\ud w\\
&\quad-\int_{0}^{1}q_{2}^{(m)}(w)\ {\mathbf{\alpha}}_{N}^{\ast}(1,1-w)\ \ud w.
\end{flalign*}

\noindent La preuve de~(\ref{eq:iaaad}) est complétée en faisant les changements de variables $\ (v,w)\to (1-v,1-w)\ $
dans cette dernière relation. En combinant~(\ref{eq:iaaaj}) avec~(\ref{eq:iaaam}) et le changement de variables $\ \displaystyle w\to 1-w,\ $
nous voyons que

\begin{displaymath}
(N+1)^{1/2}\ \Theta_{N}=-\int_{0}^{1}q_{2}^{(m)}(w)\ {\mathbf{\alpha}}_{N}(1,w)\ \ud w=\int_{0}^{1}q_{2}^{(m)}(1-w)\ {\mathbf{\alpha}}_{N}(1,1-w)\ \ud w.
\end{displaymath}

\noindent En se rappelant de~(\ref{eq:iaax}), que $\ \displaystyle {\mathbf{\alpha}}_{N}(1,1-w)=-{\mathbf{\alpha}}_{N}^{\ast}(1,w),\ $ pour $\ 0\leq w\leq 1,\ $ nous
déduisons facilement~(\ref{eq:iaaae}) de cette relation. $\Box$\vskip5pt
\begin{prop}\label{propos:propo252}Nous avons
{\small
\begin{flalign}
(N+1)^{1/2}\Delta_{N}&=2\int\int_{[0,1]^2}q_{3}^{(m)}(1-v)\ q_{2}^{(m)}(1-w)\ \left\{\alpha_{N}^{\ast}(v,w)
-v\ \alpha_{N}^{\ast}(1,w)\right.\label{eq:iaaan}\\
&\left.\hspace{6cm}-w\ \alpha_{N}^{\ast}(v,1)\right\}\ \ud v\ud w\nonumber\\
&\quad+4m\int_{0}^{1}q_{3}^{(m)}(1-v)\ \alpha_{N}^{\ast}(v,1)\ \ud v.\nonumber
\end{flalign}}
\end{prop}
\noindent{\bf Démonstration.} Nous réécrivons~(\ref{eq:iaaad}) en

{\footnotesize
\begin{flalign}
\hspace{0cm}(N+1)^{1/2}\Delta_{N}&=2\int\int_{[0,1]^2}q_{3}^{(m)}(1-v)\ q_{2}^{(m)}(1-w)\ \left\{{\mathbf{\alpha}}_{N}^{\ast}(v,w)
-v\ {\mathbf{\alpha}}_{N}^{\ast}(1,w)\right.\label{eq:iaaap}\\
\hspace{0cm}&\left.\hspace{6cm}-w\ {\mathbf{\alpha}}_{N}^{\ast}(v,1)\right\}\ \ud v\ud w\nonumber\\
\hspace{0cm}&\quad+2\left\{\int_{0}^{1}v\ q_{3}^{(m)}(1-v)\ud v\right\}\int_{0}^{1}q_{2}^{(m)}(1-w)\ {\mathbf{\alpha}}_{N}^{\ast}(1,w)\ \ud w\nonumber\\
\hspace{0cm}&\quad+2\left\{\int_{0}^{1}w\ q_{2}^{(m)}(1-w)\ud w\right\}\int_{0}^{1}q_{3}^{(m)}(1-v)\ {\mathbf{\alpha}}_{N}^{\ast}(v,1)\ \ud v\nonumber\\
\hspace{0cm}&\quad-\int_{0}^{1}q_{2}^{(m)}(1-w)\ {\mathbf{\alpha}}_{N}^{\ast}(1,w)\ \ud w,\nonumber
\end{flalign}}

\noindent où nous avons utilisé~(\ref{eq:iaaak}) et~(\ref{eq:iaaal}). Ceci donne~(\ref{eq:iaaan}), comme de\-man\-dé. $\Box$\vskip5pt
\noindent Après avoir obtenu, dans les propositions~\ref{propos:propo251} et~\ref{propos:propo252}, les représentations appropriées de $\ \displaystyle
(N+1)^{1/2}\Delta_{N}\ $ et $\ \displaystyle (N+1)^{1/2}\ \Theta_{N}$ en termes de $\ \displaystyle \alpha_{N}\ $ et $\ \alpha_{N}^{\ast},\ $ nous allons maintenant approcher dans le paragraphe~\ref{ma-soussectionh} ces statistiques par des homologues gaussiens. Nous
devons d'abord obtenir dans le paragraphe~\ref{ma-soussectione} des bornes supérieures pour les versions pondérées de ces
processus empiriques. Ces préliminaires sont rendus nécessaires du fait que la fonction $\ \displaystyle q_{2}^{(m)}(.)\ $
dans~(\ref{eq:iaaad}),~(\ref{eq:iaaae}) et~(\ref{eq:iaaan}) n'est pas bornée sur $\ \displaystyle (0,1).\ $
\subsection{Processus empiriques pondérés} \label{ma-soussectione}
\noindent Dans ce sous-paragraphe, nous allons prouver certains résultats techniques, indiqués dans les lemmes~\ref{lem:le251}
et~\ref{lem:le252} ci-dessous. Ces lemmes joueront un rôle crucial dans la preuve de la prochaine proposition~\ref{propos:propo256},
donnée plus tard, dans le sous-paragraphe~\ref{ma-soussectionk}. La notation suivante sera nécessaire. Posons, pour chaque $\ \displaystyle N\geq
1,\ $
\begin{equation}
I_{N}^{\prime}:=\int_{0}^{1}\int_{0}^{1/N}q_{3}^{(m)}(1-v)\ q_{2}^{(m)}(1-w)\ \alpha_{N}^{\ast}(v,w)\ \ud v\ud w,\label{eq:iaaaq}
\end{equation}
\noindent et
\begin{equation}
K_{N}^{\prime}:=\int_{0}^{1/N}q_{2}^{(m)}(1-w)\ \alpha_{N}^{\ast}(1,w)\ \ud w.\label{eq:iaaar}
\end{equation}
\noindent Nous obtiendrons des limites de bornes supérieures pour $\ \displaystyle\left\vert I_{N}^{\prime}\right\vert\ $ et $\
\displaystyle\left\vert K_{N}^{\prime}\right\vert.\ $ Vers ce but, et compte tenu de~(\ref{eq:iaav}), nous introduisons la version continue à
droite de la f.r. empirique basée sur la suite $\ \displaystyle \left\{\left(V_{k,1},W_{k,1}\right):0\leq k\leq N\right\},\ $ de vecteurs
aléatoires i.i.d. avec des lois uniformes sur $\ \displaystyle [0,1]^2.\ $ Cette f.r. empirique est définie, pour $\ \displaystyle\ 0\leq
v,w\leq1,\ $ par
\begin{equation}
U_{N}^{[s]}(v,w)=\frac{1}
{N+1}\#\left\{V_{k,1}\leq v,W_{k,1}\leq w:0\leq k\leq N\right\}.\label{eq:iaaat}
\end{equation}
\noindent Le processus empirique correspondant est défini, pour $\ \displaystyle0\leq v,w\leq1,\ $ par
\begin{equation}
\alpha_{N}^{[s]}(v,w)=\left(N+1\right)^{1/2}
\left(U_{N}^{[s]}
(v,w)-vw\right).\label{eq:iaaau}
\end{equation}
\noindent Nos arguments dépendront du fait~\ref{fai:fa200} suivant. Nous notons que, dans ce fait, la relation~(\ref{eq:iaaaac}) est une
conséquence du corollaire 2 de Einmahl et Mason\cite{MR807337}, pris avec $\ \displaystyle\nu=\frac{1}{2}\ $ et $\ \displaystyle d=2,\ $ alors
que la relation~(\ref{eq:iaaaad}) est due à Cs\'aki\cite{MR0391233,MR719704}.
\begin{fait}
\label{fai:fa200}
\noindent Nous avons,
\begin{equation}
\limsup_{N\to \infty}\frac{1}{\log\log N}\log
\left\{\sup_{0<v,w\leq 1}\left\vert\frac{\alpha^{[
s]}_{N}(v,w)}{\sqrt{vw(1-vw)}}\right\vert\right\}=1\quad\mbox{ p.s.}
\label{eq:iaaaac}
\end{equation}
\noindent et
\begin{equation}
\limsup_{N\to \infty}\frac{1}{\log\log N}\log \left\{\sup_{0<w\leq
1}\left\vert\frac{\alpha^{[
s]}_{N}(1,w)}{\sqrt{w(1-w)}}\right\vert\right\}=\frac{1}{2}\quad\mbox{ p.s.}\label{eq:iaaaad}
\end{equation}
\end{fait}
\begin{lemme}
\label{lem:le251}
\noindent Pour chaque $\varepsilon>0$, nous avons, lorsque
$\ N\to \infty,\ $
\begin{equation}
\vert I^{\prime}_{N}\vert=O_{\mathbb{P}}\left(\frac{(\log
N)^{1+\varepsilon}}{\sqrt{N}}\right).\label{eq:iaaaae}
\end{equation}
\end{lemme}
\noindent{\bf Démonstration.} Par le fait~\ref{fai:fa200}, pour chaque $\ \displaystyle C_{0}>0\ \mbox{ et }\ \varepsilon>0,\ $ il existe un $\
\displaystyle N_{0}<\infty\ \ \mbox{ p.s. },\ $ tel que, pour tout $\ \displaystyle N\geq N_{0}\ $ et $\ \displaystyle 0<v,w\leq 1,\ $

\begin{equation}
\left\vert{\mathbf{\alpha}}_{N}^{[s]}(v,w)\right\vert\leq C_{0}
\ (\log N)^{1+\varepsilon}\ \sqrt{vw}.\label{eq:iaaaaf}
\end{equation}

\noindent En se rappelant de la définition~(\ref{eq:iaaaq}) de $\ \displaystyle I_{N}^{\prime},\ $ nous observons que

\begin{displaymath}
I_{N}^{\prime}\stackrel{d}{=}I_{N}^{\prime [s]}:=\int_{0}^{1}\int_{0}^{1/N}q_{3}^{(m)}(1-v)\ q_{2}^{(m)}(1-w)\
{\mathbf{\alpha}}_{N}^{[s]}(v,w)\ \ud v\ \ud w.
\end{displaymath}

\noindent Comme suit de~(\ref{eq:iaaaaf}), nous avons, pour tout $\ \displaystyle N\geq N_{0},\ $

\begin{equation}
\left\vert I_{N}^{\prime [s]}\right\vert\leq C_{0}\ (\log
n)^{1+\varepsilon}
\int_{0}^{1}\int_{0}^{1/N}q_{3}^{(m)}(1-v)\ q_{2}^{(m)}(1-w)\ \sqrt{vw}\ \ud v\ \ud w.\label{eq:iaaaag}
\end{equation}

\noindent En se rappelant de~(\ref{eq:iaak}), que $\ \displaystyle q_{2}^{(m)}(1-w)=\frac{1+o(1)}{w}\ $ lorsque $\ w\downarrow 0,\ $ nous
déduisons de~(\ref{eq:iaaaag}) que, pour chaque $\ \displaystyle1<C_{0}<C_{1}<\infty,\ $ il existe un $\ \displaystyle N_{1}\geq N_{0},\ $ tel
que, pour tout $\ \displaystyle N\geq N_{1},\ $

\begin{displaymath}
\left\vert I_{N}^{\prime [s]}\right\vert\leq C_{1}\ (\log
N)^{1+\varepsilon}\left\{\int_{0}^{1}\sqrt{v}\ q_{3}^{(m)}(1-v)\ \ud v\right\} \int_{0}^{1/N}\frac{\ud w}{\sqrt{w}}.
\end{displaymath}

\noindent Or $\ \displaystyle v\mapsto q_{3}^{(m)}(1-v)\ $ est bornée sur $\ \displaystyle(0,1).\ $ De plus,

\begin{equation}
\forall\ v\in(0,1),\ 0<q_{3}^{(m)}(1-v)<1,\label{eq:iaaaaaax}
\end{equation}

\noindent d'où,

\begin{displaymath}
\left\vert I_{N}^{\prime [s]}\right\vert\leq C_{1}\ (\log
n)^{1+\varepsilon}\sup_{0<v<1}q_{3}^{(m)}(1-v)\left\{\int_{0}^{1}\sqrt{v}\ \ud v\right\} \int_{0}^{1/N}\frac{\ud w}{\sqrt{w}}=\frac{4}{3}\ \
C_{1}\ \frac{(\log N)^{1+\varepsilon}}{\sqrt{N}}.
\end{displaymath}

\noindent Ceci, avec le fait que la convergence presque sûre implique la convergence en probabilité, implique que lorsque $\ \displaystyle
N\to\infty\ $ $\ \displaystyle \left\vert I_{N}^{\prime [s]}\right\vert=O_{\mathbb{P}}\left(N^{-1/2}(\log N)^{1+\varepsilon}\right).\ $ Puisque
$\ \displaystyle I_{N}^{\prime}\stackrel{d}{=}I_{N}^{\prime [s]},\ $ ceci implique~(\ref{eq:iaaaae}). $\Box$\vskip5pt
\begin{lemme}
\label{lem:le252}
\noindent Pour chaque $\varepsilon>0$, nous avons, lorsque
$\ N\to\infty,\ $
\begin{equation}
\vert K^{\prime}_{N}\vert=O_{\mathbb{P}}\bigg(\frac{(\log
N)^{\displaystyle 1/2+\varepsilon}}{\sqrt{N}}\bigg).\label{eq:iaaaah}
\end{equation}
\end{lemme}
\noindent{\bf Démonstration.} En se rappelant de la définition~(\ref{eq:iaaar}) de $\ \displaystyle K^{\prime}_{N},\ $ nous observons que

\begin{equation}
K_{N}^{\prime}\stackrel{d}{=}K_{N}^{\prime [s]}:=\int_{0}^{1/N}q_{2}^{(m)}(1-w)
\ {\mathbf{\alpha}}_{N}^{[s]}(1,w)\ \ud w.\label{eq:iaaaai}
\end{equation}

\noindent Par~(\ref{eq:iaaaad}) dans le fait~\ref{fai:fa200}, pour chaque $\ \displaystyle C_{2}>0\ \mbox{ et }\ \varepsilon>0,\ $ il existe un
$\ \displaystyle N_{2}<\infty\ \ \mbox{ p.s.,}$ tel que, pour tout $\ \displaystyle N\geq N_{2}\ $ et $\ \displaystyle 0\leq w\leq 1,\ $

\begin{equation}
\left\vert{\mathbf{\alpha}}_{N}^{[s]}(1,w)\right\vert\leq C_{2}\
(\log N)^{\displaystyle\frac{1}{2}+\varepsilon}\sqrt{w},\label{eq:iaaaaj}
\end{equation}

\noindent d'où, par~(\ref{eq:iaaaai}),

\begin{equation}
\left\vert K_{N}^{\prime [s]}\right\vert\leq C_{2}\ (\log
N)^{\displaystyle\frac{1}{2}+\varepsilon}\int_{0}^{1/N}q_{2}^{(m)}(1-w)\ \sqrt{w}\ \ud w.\label{eq:iaaaak}
\end{equation}

\noindent En se rappelant de~(\ref{eq:iaak}), que $\ \displaystyle q_{2}^{(m)}(1-w)=\frac{1+o(1)}{w}\ $ lorsque $\ w\downarrow 0,\ $
nous déduisons de~(\ref{eq:iaaaak}) que, pour chaque $\ \displaystyle1<C_{2}<C_{3}<\infty,\ $ il existe un $\ \displaystyle N_{3}\geq N_{2}\ $
tel que, pour tout $\ \displaystyle N\geq N_{3},\ $

\begin{displaymath}
\left\vert K_{N}^{\prime [s]}\right\vert\leq C_{3}\ (\log
N)^{\displaystyle\frac{1}{2}+\varepsilon}\int_{0}^{1/N}\frac{\ud w}{\sqrt{w}}\leq 2\ C_{3}\ \frac{(\log
N)^{\displaystyle\frac{1}{2}+\varepsilon}}{\sqrt{N}}.
\end{displaymath}

\noindent Puisque la convergence presque sûre implique la convergence en probabilité, cela implique que $\ \displaystyle \left\vert
K_{N}^{\prime [s]}\right\vert=O_{\mathbb{P}}\left(N^{-1/2}(\log N)^{1/2+\varepsilon}\right)\ $ lorsque $\ \displaystyle N\to\infty.\ $ De plus,
on a,  $\ \displaystyle K_{N}^{\prime}\stackrel{d}{=}K_{N}^{\prime [s]},\ $ donc en déduit~(\ref{eq:iaaaah}).$\Box$\vskip5pt
\section{Approximations gaussiennes}\label{ma-sectiond}
\subsection{Ponts browniens et processus de Wiener}\label{ma-soussectionf}
\noindent La notation et les faits suivants concernant les processus de Wiener et les ponts browniens seront nécessaires. Nous rappelons les
définitions usuelles d'un processus de Wiener (univarié)$\ \displaystyle \left\{W(u):u\geq 0\right\} $ et d'un pont Brownien $\
\displaystyle\left\{B(u):0\leq u\leq 1\right\},\ $ tous deux définis comme des processus gaussiens centrés avec un échantillon
de trajectoires continues et des covariances vérifiant respectivement,

\begin{equation}
\mathbb{E}\left(W(u)\ W(v)\right)=u\land v,\ \mbox{ pour }\ u,v\geq 0,\label{eq:iaaaal}
\end{equation}

\noindent et

\begin{equation}
\mathbb{E}\left(B(u)\ B(v)\right)=u\land v-uv,\ \mbox{ pour }\ 0\leq u,v\leq 1.\label{eq:iaaaam}
\end{equation}

\noindent Nous pouvons définir ces processus sur le même espace de probabilité, en posant

\begin{equation}
B(u)=W(u)-u\ W(1),\ \mbox{ pour }\ 0\leq u\leq 1.\label{eq:iaaaan}
\end{equation}

\noindent Nous notons que pour un usage ultérieur, lorsque $\ \displaystyle B(.) \ \mbox{ et }\ W(.)\ $ sont tels que dans~(\ref{eq:iaaaan}),
$\ \displaystyle \left\{B(u):0\leq u\leq 1\right\}\ $ et $\ \displaystyle W(1)\stackrel{d}{=}N(0,1)\ $ sont indépendants.\vskip5pt

\noindent Ces définitions s'étendent aux processus bivariés comme suit. Par un {\it processus bivarié}, on sous-entend une fonction aléatoire de
$\ \displaystyle (u,v)\in A,\ $ où $\ \displaystyle A\ $ désigne un sous-ensemble de $\ \mathbb{R}^2.\ $ Nous définissons un {\it processus de
Wiener bivarié} $\ \displaystyle \left\{\mathbf{W}(u,v):u,v\geq 0\right\}\ $ comme un processus gaussien centré sur $\ \mathbb{R}^2_{+},\ $ avec
un échantillon de trajectoires continues, et une fonction de covariance vérifiant, pour tout $\ \displaystyle\ u^{\prime},
u^{\prime\prime},v^{\prime},v^{\prime\prime}\geq 0,\ $

\begin{equation}
\mathbb{E}\left(\mathbf{W}(u^{\prime},v^{\prime})\ \mathbf{W}(u^{\prime\prime},v^{\prime\prime})\right)=\left(u^{\prime}\land
u^{\prime\prime}\right)\left(v^{\prime}\land v^{\prime\prime}\right).\label{eq:iaaaao}
\end{equation}

\noindent Un pont Brownien bivarié $\ \displaystyle\left\{\mathbf{B}(u,v):0\leq u,v\leq 1\right\}\ $ est défini à son tour comme un processus
gaussien centré avec un échantillon de trajectoires continues sur $\ \displaystyle[0,1]^2\ $ et une fonction de covariance telle que

\begin{equation}
\mathbb{E}\left(\mathbf{B}(u^{\prime},v^{\prime})\ \mathbf{B}(u^{\prime\prime},v^{\prime\prime})\right)=\left(u^{\prime}\land
u^{\prime\prime}\right)\left(v^{\prime}\land v^{\prime\prime}\right)-u^{\prime}u^{\prime\prime}v^{\prime}v^{\prime\prime},\label{eq:iaaaap}
\end{equation}

\noindent pour $\ \displaystyle 0\leq u^{\prime}, u^{\prime\prime},v^{\prime},v^{\prime\prime}\leq 1.\ $ Nous pouvons définir $\ \displaystyle
\left\{\mathbf{B}(u,v):0\leq u,v\leq 1\right\}\ $ et $\ \displaystyle \left\{\mathbf{W}(u,v):u,v\geq 0\right\}\ $ sur le même espace de
probabilité en posant

\begin{equation}
\mathbf{B}(u,v)=\mathbf{W}(u,v)-uv\ \mathbf{W}(1,1),\ \mbox{ pour }\ 0\leq u,v\leq 1.\label{eq:iaaaaq}
\end{equation}

\noindent Quand $\ \displaystyle\mathbf{B}(.,.)\ $ et $\ \displaystyle \mathbf{W}(.,.)\ $ sont comme dans~(\ref{eq:iaaaaq}),
$\ \displaystyle\left\{\mathbf{B}(u,v):0\leq u,v\leq 1\right\}\ $ et $\ \displaystyle \mathbf{W}(1,1)\stackrel{d}{=}N(0,1)\ $ sont
indépendants.\vskip5pt

\noindent Il sera commode d'écrire $\ \displaystyle\mathbf{W}(s,t)\ $ en fonction de la {\it mesure de Wiener} $\ \displaystyle \mathbf{W}(\ud
u,\ud v)\ $ (voir par exemple Lifshits\cite{MR1472736}, p. 107), soit par l'une des intégrales

\begin{displaymath}
\mathbf{W}(s,t)=\int\int_{[0,s]\times[0,t]}\mathbf{W}(\ud u,\ud v)=\int\int_{[0,s)\times[0,t)}\mathbf{W}(\ud u,\ud v),
\end{displaymath}

\noindent pour $\ s,t\leq 0.\ $ En définissant la mesure du pont Brownien par $\ \displaystyle\mathbf{B}(\ud u,\ud v)=\mathbf{W}(\ud u,\ud
v)-\mathbf{W}(1,1)\ \ud u\ud v,\ $ nous obtenons la représentation intégrale, pour $\ 0\leq s,t\leq 1,\ $

\begin{displaymath}
\mathbf{B}(s,t)=\int\int_{[0,s]\times[0,t]}\mathbf{B}(\ud u,\ud v)=\int\int_{[0,s)\times[0,t)}\mathbf{B}(\ud u,\ud v).
\end{displaymath}

\noindent En se rappelant les égalités $\ \displaystyle\mathbf{B}(0,0)=\mathbf{B}(1,1)=0,\ $ le processus défini, pour $\ \displaystyle 0\leq
u,v\leq 1,\ $ par
\begin{equation}
\mathbf{B}^{\ast}(1-u,1-v)=\int\int_{(u,1]\times(v,1]}\mathbf{B}(\ud u,\ud v)=\mathbf{B}(u,v)-\mathbf{B}(u,1)-\mathbf{B}(1,v)\label{eq:iaaaar}
\end{equation}
\noindent vérifie l'identité distributionnelle
\begin{equation}
\left\{\mathbf{B}(u,v):0\leq u,v\leq 1\right\}\stackrel{d}{=}\left\{\mathbf{B}^{\ast}(u,v):0\leq u,v\leq 1\right\}.\label{eq:iaaaas}
\end{equation}
\noindent Ceci vaut en tant que conséquence de la stationnarité de la mesure de Wiener $\ \mathbf{W}(\ud u,\ud v)\ $ dans $\
\displaystyle\mathbb{R}^2_{+},\ $ avec l'invariance de la mesure de Lebesgue par translations sur $\ \mathbb{R}^2.\ $

\noindent Nous définissons un {\it pont Brownien réduit} $\ \displaystyle \left\{\mathbf{B}_{[0]}(u,v):0\leq u,v\leq 1\right\}\ $ comme un
processus gaussien centré avec un échantillon de trajectoires continues sur $\ \displaystyle [0,1]^2,\ $ et une fonction de covariance

\begin{equation}
\mathbb{E}\left(\mathbf{B}_{[0]}(u^{\prime},v^{\prime})\mathbf{B}_{[0]}(u^{\prime\prime},v^{\prime\prime})\right)=\left(u^{\prime}\land
u^{\prime\prime}-u^{\prime}u^{\prime\prime}\right)\left(v^{\prime}\land v^{\prime\prime}-v^{\prime}v^{\prime\prime}\right),\label{eq:iaaaat}
\end{equation}

\noindent pour tout $\ \displaystyle 0\leq u^{\prime}, u^{\prime\prime},v^{\prime},v^{\prime\prime}\leq 1.\ $\vskip5pt

\noindent Nous pouvons définir un pont Brownien réduit en fonction du pont Brownien bivarié (usuel) $\ \displaystyle\left\{\mathbf{B}(u,v):0\leq
s,t\leq 1\right\},\ $ par la relation

\begin{equation}
\mathbf{B}_{[0]}(u,v)=\mathbf{B}(u,v)-u\mathbf{B}(1,v)-v\mathbf{B}(u,1),\ \mbox{ pour }\ 0\leq u,v\leq 1.\label{eq:iaaaau}
\end{equation}

\noindent Les processus marginaux de $\ \displaystyle\mathbf{B}(u,v),\ $ sont obtenus en posant

\begin{equation}
\mathrm{B}_{[1]}(u)=\mathbf{B}(u,1),\ \mbox{ pour }\ 0\leq u\leq 1,\label{eq:iaaaav}
\end{equation}

\noindent et

\begin{equation}
\mathrm{B}_{[2]}(v)=\mathbf{B}(1,v),\ \mbox{ pour }\ 0\leq v\leq 1,\label{eq:iaaaaw}
\end{equation}

\noindent qui sont des ponts Browniens, avec des fonctions de covariance comme dans\\*
~(\ref{eq:iaaaam}). Les processus définis
par~(\ref{eq:iaaaau}),~(\ref{eq:iaaaav}) et~(\ref{eq:iaaaaw}) sont liés comme suit.

\begin{lemme}\label{lem:le253}
Les processus $\ \displaystyle\mathbf{B}_{[0]}(.,.),\ \mathrm{B}_{[1]}(.)\ \mbox{ et }\ \mathrm{B}_{[2]}(.)\ $ sont indépendants, et
constituent, respectivement un pont Brownien bivarié réduit et deux ponts Browniens univariés vérifiant, pour $\ 0\leq u,v\leq 1,\ $
\begin{equation}
\mathbf{B}(u,v)=\mathbf{B}_{[0]}(u,v)+v\ \mathrm{B}_{[1]}(u)+u\ \mathrm{B}_{[2]}(v).\label{eq:iaaaax1}
\end{equation}
\end{lemme}
\noindent{\bf Démonstration.} Voir le théorème 1, p. 105 dans Deheuvels\cite{MR612295} et par exemple, le fait 2.1. p. 499 dans Deheuvels,
Peccati et Yor\cite{MR2199561}.$\Box$\vskip5pt
\subsection{Identités distributionnelles pour les intégrales de ponts Browniens}\label{ma-soussectiong}
\noindent Dans ce paragraphe, nous établissons des propriétés distributionnelles in\-té\-res\-san\-tes d'intégrales de fonctionnelles de ponts
browniens qui sont nécessaires dans la preuve de nos théorèmes. La proposition~\ref{propos:propo253} et le corollaire~\ref{coro:cor201} sont respectivement les mêmes que la proposition 4.1 et le corollaire 4.1 de Deheuvels et Derzko~\cite{MR2728435}. Pour les démonstrations des propositions~\ref{propos:propo253},~\ref{propos:propo254} et du corollaire~\ref{coro:cor201},  nous nous référons à Deheuvels et Derzko~\cite{MR2728435}. Soit $\ \displaystyle X\geq 0,\ $ une v.a. non négative avec une f.r. $\ \displaystyle F (x) = P (X\leq x),\ $ et $\ \displaystyle \left\{B(t): 0 \leq t\leq 1\right\}\ $ qui désigne un pont Brownien (univarié).
\begin{prop}\label{propos:propo253}Supposons que $\ \displaystyle \sigma_{F}^2:={\rm Var}(X)<\infty.\ $ Alors,
\begin{equation}
\int_{0}^{\infty}\mathrm{B}(F(t))\ \ud t\stackrel{d}{=}N(0,\sigma_{F}^2).\label{eq:iaaaax}
\end{equation}
\end{prop}

\noindent Les lemmes~\ref{lem:le254} et~\ref{lem:le255} suivants sont nécessaires pour la preuve de la proposition~\ref{propos:propo253}.

\begin{lemme}\label{lem:le254}Chaque fois que $\ \displaystyle\mathbb{E}(X)<\infty,\ $ nous avons

\begin{equation}
\lim_{x\to\infty}x(1-F(x))=0\ \mbox{ et }\ \mathbb{E}(X)=\int_{0}^{\infty}(1-F(t))\ \ud t.\label{eq:iaaaay}
\end{equation}

\noindent Si, en plus $\ \displaystyle\mathbb{E}(X^2)<\infty,\ $ alors nous avons

\begin{equation}
\lim_{x\to\infty}x^2(1-F(x))=0,\ \lim_{x\to\infty}\int_{x}^{\infty}\left(1-F(y)\right)\ud y=0,\label{eq:iaaaaz}
\end{equation}

\noindent et

\begin{equation}
\mathbb{E}(X^2)=2\int_{0}^{\infty}y(1-F(y))\ \ud y.\label{eq:iaaaaaa}
\end{equation}

\end{lemme}

\noindent{\bf Démonstration.} Pour $\ \displaystyle k=1,2,\ldots,\ $ l'hypothèse que $\ \displaystyle\mathbb{E}(X^k)<\infty\ $ implique que
$\ \displaystyle\mathbb{E}(X^k{\rm 1\!I}_{\left\{X\geq x\right\}})\to 0\ $ lorsque $\ \displaystyle x\to\infty.\ $ L'inégalité de Markov
$\ \displaystyle x^k(1-F(x))\leq\mathbb{E}(X^k{\rm 1\!I}_{\left\{X\geq x\right\}}),\ $ implique aussi que $\ x^k(1-F(x))\to 0\ $ lorsque $\ x\to
\infty.\ $ La preuve de~(\ref{eq:iaaaay}) est obtenue en posant $\ \displaystyle k=1\ $ dans cette relation, puis, en intégrant par parties pour
obtenir

\begin{flalign*}
\mathbb{E}(X)&=\lim_{x\to\infty}\int_{0}^{x}y\ \ud F(y)=\lim_{x\to\infty}\left\{\left[-y(1-F(y))\right]_{0}^{x}\right.\\
&\left.+\int_{0}^{x}(1-F(y))\ \ud y\right\}=\int_{0}^{\infty}(1-F(y))\ \ud y.
\end{flalign*}

\noindent Pou établir~(\ref{eq:iaaaaz}), nous appliquons l'argument ci-dessus avec $\ \displaystyle k=2,\ $ pour obtenir que, chaque fois $\
\displaystyle  \mathbb{E}(X^2)<\infty,\ $ nous avons $\ x^2(1-F(x))\to 0\ $ lorsque $\ x\to\infty.\ $ Ceci, combiné à une intégration par
parties montre que

\begin{flalign*}
\mathbb{E}(X^2)&=\lim_{x\to\infty}\int_{0}^{x}y^2\ \ud F(y)=\lim_{x\to\infty}\left\{\left[-y^2(1-F(y))\right]_{y=0}^{y=x}\right.\\
&\left.+2\int_{0}^{x}y(1-F(y))\ \ud y\right\}=2\int_{0}^{\infty}y(1-F(y))\ \ud y,
\end{flalign*}

\noindent ce qui donne~(\ref{eq:iaaaaaa}) comme demandé.$\Box$\\*
\\*

\begin{lemme}\label{lem:le255} Sous l'hypothèse que $\ \displaystyle \sigma^{2}_{F}={\rm Var}(X)<\infty,\ $ nous avons

\begin{equation}
\sigma^{2}_{F}=\int_{0}^{\infty}\int_{0}^{\infty}\left\{F(x)\land F(y)-F(x)F(y)\right\}\ud x\ud y.\label{eq:iaaaaab}
\end{equation}

\end{lemme}

\noindent{\bf Démonstration.} Supposons que $\ \displaystyle \sigma^{2}_{F}={\rm Var}(X)<\infty\ $ Nous avons la chaîne d'égalités

\begin{flalign}
I&:=\int_{0}^{\infty}\int_{0}^{\infty}\left\{F(x)\land F(y)-F(x)F(y)\right\}\ud x\ud y\label{eq:iaaaaac}\\
&=\int_{0}^{\infty}\left\{\left(1-F(x)\right)\int_{0}^{x}F(y)\ \ud y+F(x)\int_{x}^{\infty}(1-F(y))\ud y\right\}\ud x\nonumber\\
&=\int_{0}^{\infty}\left\{-\left(1-F(x)\right)\int_{0}^{x}\left(1-F(y)\right)\ \ud y+x\left(1-F(x)\right)\right.\nonumber\\
&\left.\quad\quad-\left(1-F(x)\right)\int_{x}^{\infty}\left(1-F(y)\right)\ \ud y+\int_{x}^{\infty}\left(1-F(t)\right)\ \ud t\right\}\ud
x,\nonumber\\
&=\int_{0}^{\infty}x\left(1-F(x)\right)\ \ud x-\left\{\int_{0}^{\infty}\left(1-F(x)\right)\ \ud x\right\}^2\nonumber\\
&\quad\quad+\int_{0}^{\infty}\left\{\int_{x}^{\infty}\left(1-F(y)\right)\ \ud y\right\}\ud x\nonumber
\end{flalign}

\noindent qui, compte tenu de~(\ref{eq:iaaaay}) et~(\ref{eq:iaaaaz}), est égal à

\begin{equation}
I=\frac{1}{2}\mathbb{E}(X^2)-\left(\mathbb{E}(X)\right)^2+\int_{0}^{\infty}\left\{\int_{x}^{\infty}\left(1-F(y)\right)\ \ud y\right\}\ud
x.\label{eq:iaaaaad}
\end{equation}

\noindent Par intégration par parties, nous déduisons de~(\ref{eq:iaaaaz}) et~(\ref{eq:iaaaaaa}) que, lorsque $\ z\to\infty,\ $

\begin{flalign*}
\hspace{0cm}\int_{0}^{z}\left\{\int_{x}^{\infty}\left(1-F(y)\right)\ \ud y\right\}\ud x&=\left[x\int_{x}^{\infty}\left(1-F(y)\right)\ \ud
y\right]_{x=0}^{x=z}+\int_{0}^{z}x\left(1-F(x)\right)\ \ud x\\
&\longrightarrow\int_{0}^{\infty}y\left(1-F(y)\right)\ \ud y=\frac{1}{2}\mathbb{E}(X^2).
\end{flalign*}

\noindent Ceci, combiné à ~(\ref{eq:iaaaaac}) et~(\ref{eq:iaaaaad}), achève la démonstration de~(\ref{eq:iaaaaab}) .$\Box$\\*
\\*

\noindent{\bf Démonstration de la proposition~\ref{propos:propo253}.} Pour chaque $\ \displaystyle 0<T<\infty,\ $ la variable aléatoire

\begin{displaymath}
I_{T}:=\int_{0}^{T}B(F(x))\ \ud x\stackrel{d}{=}N\left(0,\sigma^{2}_{F;T}\right),
\end{displaymath}

\noindent est normale centrée avec une variance donnée par

{\small
\begin{flalign}
\hspace{0cm}\sigma^{2}_{F;T}&:=\mathbb{E}\left(\left\{\int_{0}^{T}B(F(x))\ \ud
x\right\}^2\right)=\mathbb{E}\left(\int_{0}^{T}\int_{0}^{T}B(F(x))B(F(y))\ \ud x\ud y\right)\label{eq:iaaaaae}\\
\hspace{0cm}&=\int_{0}^{T}\int_{0}^{T}\left\{F(x)\land F(y)-F(x)F(y)\right\}\ \ud x\ud y,\nonumber
\end{flalign}}

\noindent où nous avons utilisé le théorème de Fubini. En faisant $\ \displaystyle T\to\infty\ $ dans~(\ref{eq:iaaaaae}), nous voyons que $\
\displaystyle \sigma^{2}_{F;T}\to\sigma^{2}_{F}\ $ lorsque $\ \displaystyle T\to\infty,\ $ où $\ \displaystyle\sigma^{2}_{F}\ $ est comme
dans~(\ref{eq:iaaaaab}). Nous faisons alors l'utilisation du lemme~\ref{lem:le255}, pour montrer que $\ \displaystyle\sigma^{2}_{F}={\rm Var
X},\ $ qui à son tour, donne~(\ref{eq:iaaaax}). $\Box$\vskip5pt
\noindent Soit $\ \displaystyle \left\{B(u):0\leq u\leq 1\right\}\ $ un pont Brownien et soit $\ \displaystyle \left\{B_{[0]}(u):0\leq u\leq 1\right\}\ $
un pont Brownien bivarié réduit. Nous rappelons que les fonctions de covariance de ces processus sont données, respectivement,
par~(\ref{eq:iaaaam}) et~(\ref{eq:iaaaat}).
\begin{prop}\label{propos:propo254} Chaque fois que $\ \displaystyle \sigma^{2}_{F}:=\rm{Var}(X)<\infty,\ $ nous avons l'identité
distributionnelle,
\begin{equation}
\left\{\int_{0}^{\infty}\mathbf{B}_{[0]}(F(x),w)\ \ud x:0\leq w\leq 1\right\}\stackrel{d}{=}\left\{\sigma_{F}B(w):0\leq w\leq
1\right\}.\label{eq:iaaaaaae}
\end{equation}
\end{prop}
\noindent{\bf Démonstration.} Le membre de gauche de~(\ref{eq:iaaaaaae}) définit un processus gaussien centré avec une fonction de covariance

\begin{flalign*}
&\mathbb{E}\left(\left\{\int_{0}^{\infty}\mathbf{B}_{[0]}(F(x),w)\ \ud x\right\}\left\{\int_{0}^{\infty}\mathbf{B}_{[0]}(F(y),v)\ \ud
y\right\}\right)\\
&=\int_{0}^{\infty}\int_{0}^{\infty}\mathbb{E}\left(\mathbf{B}_{[0]}(F(x),w)\mathbf{B}_{[0]}(F(y),v)\right)\ \ud x\ud y\\
&=\left(w\land v-wv\right)\int_{0}^{\infty}\int_{0}^{\infty}\left\{F(x)\land F(y)-F(x)F(y)\right\}\ud x\ud y\\
&=\left(w\land v-wv\right)\sigma_{F}^{2},\ \mbox{ pour }\ 0\leq w,v\leq 1,
\end{flalign*}

\noindent où nous avons utilisé~(\ref{eq:iaaaaab}). En se rappelant de la forme~(\ref{eq:iaaaam}) de la fonction de covariance d'un pont
Brownien, nous voyons que ce dernier résultat est équivalent à~(\ref{eq:iaaaaaae}).$\Box$\\*
\\*
\begin{cor}\label{coro:cor201} Soit $\ \displaystyle X\geq 0\ \mbox{ et }\ Y\geq 0\ $  des variables aléatoires non négatives, avec des f.r. $\
\displaystyle F(x)=\mathbb{P}\left(X\leq x\right)\ $ et $\ \displaystyle G(y)=\mathbb{P}\left(Y\leq y\right).\ $ Supposons que
$\ \displaystyle\sigma_{F}^{2}:={\rm Var}(X)<\infty,\ $ et $\ \displaystyle\sigma_{G}^{2}:={\rm Var}(Y)<\infty,\ $ Alors, nous avons
\begin{equation}
\int_{0}^{\infty}\int_{0}^{\infty}\mathbf{B}_{[0]}(F(x),G(y))\ \ud x\ud
y\stackrel{d}{=}N\left(0,\sigma_{F}^{2}\sigma_{G}^{2}\right).\label{eq:iaaaaag}
\end{equation}
\end{cor}
\noindent{\bf Démonstration.} Par~(\ref{eq:iaaaaaae}), nous avons l'identité distributionnelle

{\small
\begin{equation}
\left\{\int_{0}^{\infty}\mathbf{B}_{[0]}(F(x),G(y))\ \ud x:0\leq y<\infty\right\}\stackrel{d}{=}\left\{\sigma_{F}B(G(y)):0\leq
y<\infty\right\},\label{eq:iaaaaah}
\end{equation}
}

\noindent où $\ \displaystyle \left\{B(u):0\leq u\leq 1\right\}\ $ désigne un pont Brownien. En se rappelant de~(\ref{eq:iaaaax}) que

\begin{displaymath}
\int_{0}^{\infty}B(G(y))\ \ud y\stackrel{d}{=}N\left(0,\sigma_{G}^{2}\right),
\end{displaymath}

\noindent nous déduisons facilement~(\ref{eq:iaaaaag}) de~(\ref{eq:iaaaaah}).$\Box$\vskip5pt
\subsection{Décompositions de ponts browniens.}\label{ma-soussectionh}
\noindent Rappelons les définitions~(\ref{eq:iaad}),~(\ref{eq:iaaf}),~(\ref{eq:iaae}),~(\ref{eq:iaag}),~(\ref{eq:iaah}) et~(\ref{eq:iaai}) de $\
\displaystyle Q_{3}^{(m)},\ Q_{2}^{(m)},$ $\displaystyle\ q_{3}^{(m)}\ \mbox{ et }\ q_{2}^{(m)}.\ $ Soit $\
\displaystyle\left\{\mathbf{B}(v,w):0\leq v,w\leq 1\right\}\ $ qui désigne un pont brownien bivarié, et, compte tenu de~(\ref{eq:iaaaas}), soit
$\ \displaystyle \mathbf{B}^{\ast}(.,.)\ $ un pont Brownien bivarié défini, comme dans~(\ref{eq:iaaaar}) par la relation
\begin{equation}
\mathbf{B}^{\ast}(1-v,1-w)=\mathbf{B}(v,w)-\mathbf{B}(v,1)-\mathbf{B}(1,w),\label{eq:iaaaaai}
\end{equation}
\noindent pour $\ \displaystyle 0\leq v,w\leq 1.\ $ Compte tenu de~(\ref{eq:iaaaax}), nous posons $\
\displaystyle \mathbf{B}_{[0]}^{\ast}(.,.),B_{[1]}^{\ast}(.)\ \mbox{ et }\ B_{[2]}^{\ast}(.), $ un pont Brownien réduit et des ponts Browniens
univariés, définis en fonction de $\ \displaystyle \mathbf{B}^{\ast},\ $ par les relations, pour $\ \displaystyle 0\leq v,w\leq 1,\ $
\begin{equation}
\mathbf{B}_{[0]}^{\ast}(v,w)=\mathbf{B}^{\ast}(v,w)-v\mathbf{B}^{\ast}(1,w)-w\mathbf{B}^{\ast}(v,1),\label{eq:iaaaaaj}
\end{equation}
\noindent et
\begin{equation}
B_{[1]}^{\ast}(v)=\mathbf{B}^{\ast}(v,1) \mbox{ et }\ B_{[2]}^{\ast}(w)=\mathbf{B}^{\ast}(1,w).\label{eq:iaaaaak}
\end{equation}
\noindent La proposition suivante peut être interprétée comme une version limite des propositions~\ref{propos:propo251}
et~\ref{propos:propo252}, lorsque $\ N\to\infty.\ $
\begin{prop}\label{propos:propo255} Sous la notation et les hypothèses ci-dessus, nous avons les relations,
{\scriptsize
\begin{flalign}
&\int\int_{[0,1]^2}\left(2\ Q_{3}^{(m)}(v)-1\right)\ Q_{2}^{(m)}(w)\ \mathbf{B}(\ud v,\ud w)\label{eq:iaaaaal}\\
&=2\int\int_{[0,1]^2} q_{3}^{(m)}(1-v)\ q_{2}^{(m)}(1-w)\ \mathbf{B}^{\ast}(v,w)\ \ud v\ \ud w-\int_{0}^{1}q_{2}^{(m)}(1-w)\
\mathbf{B}^{\ast}(1,w)\ \ud w\nonumber\\
&=2\int\int_{[0,1]^2} q_{3}^{(m)}(1-v)\ q_{2}^{(m)}(1-w)\ \mathbf{B}^{\ast}_{[0]}(v,w)\ \ud v\ \ud w\nonumber\\
&\quad\quad+4m\int_{0}^{1}q_{3}^{(m)}(1-v)\ B^{\ast}_{[1]}(v)\ \ud v\nonumber\\
&=2\int\int_{[0,1]\times[0,\infty)}\mathbf{B}^{\ast}_{[0]}(1-G_{3}^{(m)}(x),1-G_{2}^{(m)}(y))\ \ud x\ \ud y\nonumber\\
&\quad\quad+4m\int_{0}^{1}B^{\ast}_{[1]}(1-G_{3}^{(m)}(x))\ \ud u\stackrel{d}{=}N(0,2m),\nonumber
\end{flalign}}
\noindent et
\begin{flalign}
&\hspace{-1cm}\int\int_{[0,1]^2} Q_{2}^{(m)}(w)\ \mathbf{B}(\ud v,\ud w)=-\int_{0}^{1}q_{2}^{(m)}(1-w)\ \mathbf{B}(1,w)\ \ud w\label{eq:iaaaaam}\\
&\hspace{-1cm}=\int_{0}^{1}q_{2}^{(m)}(1-w)\ {\mathbf{B}}^{\ast}(1,w)\ \ud w=\int_{0}^{1}q_{2}^{(m)}(1-w)\ B^{\ast}_{[2]}(w)\ \ud w,\nonumber
\end{flalign}
\noindent où les composantes aléatoires
{\small
\begin{flalign}
&2\int\int_{[0,1]^2}\mathbf{B}^{\ast}_{[0]}(v,w)\ q_{3}^{(m)}(1-v)\ q_{2}^{(m)}(1-w)\ \ud v\ \ud w\label{eq:iaaaaan}\\
&=2\int\int_{[0,1]\times[0,\infty)}\mathbf{B}^{\ast}_{[0]}(1-G_{3}^{(m)}(x),1-G_{2}^{(m)}(y))\ \ud x\ \ud
y\stackrel{d}{=}N\left(0,\frac{2m}{2m+1}\right),\nonumber
\end{flalign}}
\begin{flalign}
&\hspace{-1cm}4m\int_{0}^{1}q_{3}^{(m)}(1-v)\ B^{\ast}_{[1]}(v)\ \ud v=4m\int_{0}^{1}B^{\ast}_{[1]}(1-G_{3}^{(m)}(x))\ \ud
x\label{eq:iaaaaao}\\
&\hspace{-1cm}=-4m\int_{0}^{1}\mathbf{B}(G_{3}^{(m)}(x),1)\ \ud x\stackrel{d}{=}N\left(0,\frac{4m^2}{2m+1}\right),\nonumber
\end{flalign}
\noindent et
\begin{equation}
\int_{0}^{1}q_{2}^{(m)}(1-w)\ B^{\ast}_{[2]}(w)\ \ud w\stackrel{d}{=}N(0,2m),\label{eq:iaaaaap}
\end{equation}
\noindent sont indépendantes.
\end{prop}
\noindent{\bf Démonstration.} Soit $\ \displaystyle R\stackrel{d}{=}\beta_{m,m}\ $ $\left(\right.$resp. $\ \displaystyle
T\stackrel{d}{=}\Gamma(2m,1)\left.\right).\ $ Compte tenu de~(\ref{eq:iaaq}) $\left(\right.$resp.~(\ref{eq:iaar})$\left.\right)$, nous faisons
usage de l'identité distributionnelle $\ \displaystyle R\stackrel{d}{=}Q_{3}^{(m)}(V)\ $ $\left(\right.$resp. $\ \displaystyle
T\stackrel{d}{=}Q_{2}^{(m)}(W)$ $\left.\right)$, où $\ \displaystyle V\stackrel{d}{=}\mathcal{U}(0,1)\ $ $\left(\right.$resp. $\ \displaystyle
W\stackrel{d}{=}\mathcal{U}(0,1)\left.\right).\ $
\begin{equation}
\mathbb{E}(R)=\int_{0}^{1}Q_{3}^{(m)}(v)\ \ud v=\frac{1}{2}\label{eq:iaaaaaq}
\end{equation}
\begin{equation}
\hspace{-0,28cm}\left(\mbox{resp.}\ \mathbb{E}(T)=\int_{0}^{1}Q_{2}^{(m)}(w)\ \ud w=2m\right).\label{eq:iaaaaar}
\end{equation}
\noindent Nous déduisons de~(\ref{eq:iaaaaaq}) $\left(\right.$resp.~(\ref{eq:iaaaaar})$\left.\right)$, une preuve alternative de
~(\ref{eq:iaaak}) $\left(\right.$resp.~(\ref{eq:iaaal})$\left.\right)$. En se rappelant de~(\ref{eq:iaal0}) et~(\ref{eq:iaal})
$\left(\right.$resp.~(\ref{eq:iaam0}) et\\*
~(\ref{eq:iaam}) $\left.\right)$ que $\ \displaystyle Q_{3}^{(m)}(v)\longrightarrow 0\ $ $\left(\right.$ resp. $\ \displaystyle
Q_{2}^{(m)}(w)\longrightarrow 0\ $ $\left.\right)$ lorsque $\ \displaystyle v\downarrow 0\ $ $\left(\right.$ resp. $\ \displaystyle w\downarrow 0\ $ $\left.\right),\ $ et $\ \displaystyle (1-v)\ Q_{3}^{(m)}(v)\longrightarrow0\ $ $\left(\right.$ resp. $\ \displaystyle (1-w)\
Q_{2}^{(m)}(w)\longrightarrow0\ $ $\left.\right)$ lorsque $\ \displaystyle v\uparrow 1\ $ $\left(\right.$ resp. $\ \displaystyle w\uparrow 1\ $
$\left.\right),\ $ par~(\ref{eq:iaaaaaq}) $\left(\right.$resp.~(\ref{eq:iaaaaar})$\left.\right)$ et une intégration par parties, nous obtenons
\begin{equation}
\int_{0}^{1}(1-v)\ q_{3}^{(m)}(1-v)\ \ud v=\left[(1-v)\ Q_{3}^{(m)}(v)\right]_{0}^{1}+\int_{0}^{1}Q_{3}^{(m)}(v)\ \ud v\label{eq:iaaaaas}
\end{equation}
{\footnotesize
\begin{equation}
\hspace{-0,9cm}\left(\mbox{resp.}\ \int_{0}^{1}(1-w)\ q_{2}^{(m)}(1-w)\ \ud w=\left[(1-w)\
Q_{2}^{(m)}(w)\right]_{0}^{1}+\int_{0}^{1}Q_{2}^{(m)}(w)\ \ud w\right).\label{eq:iaaaaat}
\end{equation}
}Nous faisons usage ensuite d'un argument de type Fubini, dans l'esprit de Donati-Martin et Yor\cite{MR1118933,MR1457611}(voir par
exemple le paragraphe 3.3 de Deheuvels, Pecccati et Yor\cite{MR2199561}). Nous écrivons, par~(\ref{eq:iaaaaai}) et les changements de variables
$\ \displaystyle (v,w)\longrightarrow(1-v,1-w),\ $

{\footnotesize
\begin{flalign}
\hspace{0cm}J_{1}&:=2\int\int_{[0,1]^2}Q_{3}^{(m)}(v)\ Q_{2}^{(m)}(w)\ \mathbf{B}(\ud v,\ud w)\label{eq:iaaaaau}\\
\hspace{0cm}&=2\int\int_{[0,1]^2}\left\{\int\int_{[0,1]^2}q_{3}^{(m)}(s)\ {\rm
1\!I}_{[0,v)}(s)\ q_{2}^{(m)}(t)\ {\rm
1\!I}_{[0,w)}(t)\ \ud s\ud t\right\}\mathbf{B}(\ud v,\ud w)\nonumber\\
\hspace{0cm}&=2\int\int_{[0,1]^2}q_{3}^{(m)}(s)\ q_{2}^{(m)}(t)\left\{\int\int_{[0,1]^2}{\rm
1\!I}_{[0,v)}(s)\ {\rm
1\!I}_{[0,w)}(t)\ \mathbf{B}(\ud v,\ud w)\right\}\ud s\ \ud t\nonumber\\
\hspace{0cm}&=2\int\int_{[0,1]^2}q_{3}^{(m)}(s)\
q_{2}^{(m)}(t)\left\{\mathbf{B}(s,t)-\mathbf{B}(s,1)-\mathbf{B}(1,t)\right\}\ud s\ \ud t\nonumber\\
\hspace{0cm}&=2\int\int_{[0,1]^2}q_{3}^{(m)}(s)\ q_{2}^{(m)}(t)\ {\mathbf{B}}^{\ast}(1-s,1-t)\ \ud s\ \ud t\nonumber\\
\hspace{0cm}&=2\int\int_{[0,1]^2}q_{3}^{(m)}(1-v)\ q_{2}^{(m)}(1-w)\ {\mathbf{B}}^{\ast}(v,w)\ \ud v\ \ud w.\nonumber
\end{flalign}}

\noindent En faisant usage de~(\ref{eq:iaaaaaj}),~(\ref{eq:iaaaaas}) et~(\ref{eq:iaaaaat}), nous obtenons aussi que

\begin{flalign}
J_{1}&=2\int\int_{[0,1]^2}q_{3}^{(m)}(1-v)\ q_{2}^{(m)}(1-w)\ \mathbf{B}_{[0]}^{\ast}(v,w)\ \ud v\ \ud w\label{eq:iaaaaav}\\
&\quad\quad+2\int\int_{[0,1]^2}v\ q_{3}^{(m)}(1-v)\ q_{2}^{(m)}(1-w)\ {\mathbf{B}}^{\ast}(1,w)\ \ud v\ \ud w\nonumber\\
&\quad\quad+2\int\int_{[0,1]^2}q_{3}^{(m)}(1-v)\ w\ q_{2}^{(m)}(1-w)\ {\mathbf{B}}^{\ast}(v,1)\ \ud v\ \ud w\nonumber\\
&=2\int\int_{[0,1]^2}q_{3}^{(m)}(1-v)\ q_{2}^{(m)}(1-w)\ \mathbf{B}_{[0]}^{\ast}(v,w)\ \ud v\ \ud w\nonumber\\
&\quad\quad+\int_{0}^{1}q_{2}^{(m)}(1-w)\ B_{[2]}^{\ast}(w)\ \ud w\nonumber\\
&\quad\quad+4m\int_{0}^{1}q_{3}^{(m)}(1-v)\ B_{[1]}^{\ast}(v)\ \ud v.\nonumber
\end{flalign}

\noindent En se rappelant de~(\ref{eq:iaaaaai}) que $\ \displaystyle \mathbf{B}^{\ast}(1,w)=-\mathbf{B}(1,1-w),\ $ un argument similaire nous
permet d'écrire par le changement de variables $\ \displaystyle t\longrightarrow 1-w,\ $ la relation

\begin{flalign}
\hspace{0cm}J_{2}&:=\int\int_{[0,1]^2}Q_{2}^{(m)}(w)\ \mathbf{B}(\ud v,\ud w)\label{eq:iaaaaaw}\\
\hspace{0cm}&=\int\int_{[0,1]^2}\left\{\int_{0}^{1}{\rm
1\!I}_{[0,w)}(t)\ q_{2}^{(m)}(t)\ \ud t\right\}\mathbf{B}(\ud v,\ud w)\nonumber\\
\hspace{0cm}&=\int_{0}^{1}\left\{\int\int_{[0,1]^2}{\rm
1\!I}_{[0,w)}(t)\ \mathbf{B}(\ud v,\ud w)\right\}q_{2}^{(m)}(t)\ \ud t\nonumber\\
\hspace{0cm}&=-\int_{0}^{1}q_{2}^{(m)}(t)\ \mathbf{B}(1,t)\ \ud t=-\int_{0}^{1}q_{2}^{(m)}(t)\ B_{[2]}(t)\ \ud t\nonumber\\
\hspace{0cm}&=\int_{0}^{1}q_{2}^{(m)}(1-w)\ {\mathbf{B}}^{\ast}(1,w)\ \ud w=\int_{0}^{1}q_{2}^{(m)}(1-w)\ B_{[2]}^{\ast}(w)\ \ud
w.\nonumber
\end{flalign}

\noindent En combinant~(\ref{eq:iaaaaau}) avec~(\ref{eq:iaaaaav}) et~(\ref{eq:iaaaaaw}), nous obtenons facilement que

{\small
\begin{flalign*}
\hspace{0cm}&\int\int_{[0,1]^2}\left(2\ Q_{3}^{(m)}(v)-1\right)Q_{2}^{(m)}(w)\ \mathbf{B}(\ud v,\ud w)=J_{1}-J_{2}\\
\hspace{0cm}&=2\int\int_{[0,1]^2} q_{3}^{(m)}(1-v)\ q_{2}^{(m)}(1-w)\ \mathbf{B}^{\ast}(v,w)\ \ud v\ \ud w-\int_{0}^{1}q_{2}^{(m)}(1-w)\
\mathbf{B}^{\ast}(1,w)\ \ud w\\
\hspace{0cm}&=2\int\int_{[0,1]^2} q_{3}^{(m)}(1-v)\ q_{2}^{(m)}(1-w)\ \mathbf{B}^{\ast}_{[0]}(v,w)\ \ud v\ \ud
w+4m\int_{0}^{1}q_{3}^{(m)}(1-v)\ B^{\ast}_{[1]}(v)\ \ud v,
\end{flalign*}
}

\noindent ce qui donne les deux premières égalités de~(\ref{eq:iaaaaal}). Le fait que les composantes
aléatoires~(\ref{eq:iaaaaan}),~(\ref{eq:iaaaaao}) et~(\ref{eq:iaaaaap}) sont indépendantes provient du lemme~\ref{lem:le253}. Pour conclure la
preuve de~(\ref{eq:iaaaaal}), nous utilisons les arguments suivants. En premier lieu, nous faisons les changements de variables
$\ \displaystyle(v,w)\longrightarrow(1-v,$ $\displaystyle 1-w),\ $ et puis nous posons $\ \displaystyle v=G_{3}^{(m)}(x)\ $ avec $\
\displaystyle q_{3}^{(m)}(v)\ \ud v=\ud x\ $ et $\ \displaystyle w=G_{2}^{(m)}(y)\ $ avec $\ \displaystyle q_{2}^{(m)}(w)\ \ud w=\ud y,\ $ pour
obtenir l'égalité,

\begin{flalign*}
&\int\int_{[0,1]^2} q_{3}^{(m)}(1-v)\ q_{2}^{(m)}(1-w)\ \mathbf{B}^{\ast}_{[0]}(v,w)\ \ud v\ \ud w\\
&=\int\int_{[0,1]^2}\mathbf{B}^{\ast}_{[0]}(1-v,1-w)\ \ud v\ q_{3}^{(m)}(v)\ \ud w\ q_{2}^{(m)}(w)\\
&=\int\int_{[0,1]\times[0,\infty)}\mathbf{B}^{\ast}_{[0]}(1-G_{3}^{(m)}(x),1-G_{2}^{(m)}(y))\ \ud x\ \ud y.
\end{flalign*}

\noindent Alors, nous observons que $\ \displaystyle\left\{\mathbf{B}_{[0]}^{\ast}(v,w):0\leq v,w\leq
1\right\}\stackrel{d}{=}\left\{\mathbf{B}_{[0]}^{\ast}(1-v,1-w)\right.$ $ \displaystyle:0\leq v,w\leq 1\big\}.\ $ Puisque, la variance de $\
\displaystyle R\stackrel{d}{=}\beta_{m,m}\ $ est égale à $\ \displaystyle \sigma_{R}^{2}=\frac{1}{4(2m+1)}\ $ et la variance de $\
T\stackrel{d}{=}\Gamma(2m,1)\ $ est égale à $\ \sigma^{2}_{G_{2}^{(m)}}=2m,\ $ nous pouvons appliquer~(\ref{eq:iaaaaag}) pour obtenir que

{\small
\begin{flalign*}
\hspace{0cm}&\int\int_{[0,1]\times[0,\infty)}\mathbf{B}^{\ast}_{[0]}(1-G_{3}^{(m)}(x),1-G_{2}^{(m)}(y))\ \ud x\ \ud y\\
\hspace{0cm}&\stackrel{d}{=}\int_{0}^{1}\int_{0}^{\infty}\mathbf{B}^{\ast}_{[0]}(G_{3}^{(m)}(x),G_{2}^{(m)}(y))\ \ud x\ud y
\stackrel{d}{=}N\left(0,\sigma_{R}^{2}\ \sigma_{G_{2}^{(m)}}^{2}\right)\stackrel{d}{=}N\left(0,\frac{m}{2\left(2m+1\right)}\right),
\end{flalign*}
}

\noindent ce qui donne~(\ref{eq:iaaaaan}). De même, en faisant le changement de variable
$\ \displaystyle v\longrightarrow 1-v\ $ et puis en posant $\ \displaystyle v=G_{3}^{(m)}(x)\ $ et $\ \displaystyle q_{3}^{(m)}(v)\ \ud v=\ud
x,\ $ on obtient l'égalité

\begin{displaymath}
\int_{0}^{1}q_{3}^{(m)}(1-v)\ B^{\ast}_{[1]}(v)\ \ud v=\int_{0}^{1}B^{\ast}_{[1]}(1-v)\ q_{3}^{(m)}(v)\ \ud
v=\int_{0}^{1}B^{\ast}_{[1]}(1-G_{3}^{(m)}(x))\ \ud x.
\end{displaymath}

\noindent Or, d'après~(\ref{eq:iaaaaai}), on a, pour tout $\ \displaystyle 0\leq x\leq 1,\ $ $ \displaystyle
B^{\ast}_{[1]}(1-G_{3}^{(m)}(x))=-B(G_{3}^{(m)}(x),1)\ $  et comme la variance de $\ \displaystyle R\stackrel{d}{=}\beta_{m,m}\ $ est égale à $\
\displaystyle \sigma_{R}^{2}=\frac{1}{4(2m+1)},\ $ nous pouvons appliquer~(\ref{eq:iaaaax}) pour obtenir que,

{\small
\begin{displaymath}
\int_{0}^{1}B^{\ast}_{[1]}(1-G_{3}^{(m)}(x))\ \ud x=-\int_{0}^{1}\mathbf{B}(G_{3}^{(m)}(x),1)\ \ud
x\stackrel{d}{=}N\left(0,\sigma_{R}^{2}\right)\stackrel{d}{=}N\left(0,\frac{1}{4(2m+1)}\right),
\end{displaymath}
}

\noindent ce qui donne~(\ref{eq:iaaaaao}). La preuve est complétée par l'observation que si\\*
$\ \displaystyle Y\stackrel{d}{=}N\left(0,\frac{m}{2\left(2m+1\right)}\right)\ $ et $\ \displaystyle
Z\stackrel{d}{=}N\left(0,\frac{1}{4(2m+1)}\right)\ $ sont indépendantes, alors $\ \displaystyle 2\ Y+4m\
Z\stackrel{d}{=}N\left(0,\frac{4m}{2\left(2m+1\right)}+\frac{16 m^2}{4(2m+1)}\right)\stackrel{d}{=}N\left(0,2m\right).\ \Box$\vskip5pt
\subsection{Ponts Browniens pondérés}\label{ma-soussectioni}
\noindent Dans ce paragraphe, nous donnons des bornes pour des ponts browniens pondérés parallèlement à celles obtenues précédemment pour les
processus empiriques pondérés dans le paragraphe~\ref{ma-soussectione}. Nous supposons que $\ \displaystyle \left\{\mathbf{B}(v,w)\right.$
$\left.\displaystyle:0\leq v,w\leq 1\right\}\ $ désigne un pont Brownien bivarié. A la suite de~(\ref{eq:iaaaax1}), nous définissons un pont
Brownien réduit $\ \displaystyle \ \left\{{\mathbf{B}}_{[0]}(v,w):0\leq v,w\leq 1\right\}\ $ en posant
\begin{equation}
\mathbf{B}(v,w)={\mathbf{B}}_{[0]}(v,w)+w\mathbf{B}(v,1)+v\mathbf{B}(1,w).\label{eq:iaaaaaaa}
\end{equation}
\noindent Dans l'esprit de~(\ref{eq:iaaaq}) et~(\ref{eq:iaaar}), nous posons, pour chaque $\ \displaystyle N\geq 1,\ $
\begin{equation}
J_{N}^{\prime}:=\int_{0}^{1}\int_{0}^{1/N}q_{3}^{(m)}(1-v)\ q_{2}^{(m)}(1-w)\ \mathbf{B}(v,w)\ \ud v\ \ud w,\label{eq:iaaaaaab}
\end{equation}
\noindent et
\begin{equation}
L_{N}^{\prime}:=\int_{0}^{1/N}q_{2}^{(m)}(1-w)\ \mathbf{B}(1,w)\ \ud w\ \ud w.\label{eq:iaaaaaac}
\end{equation}
\begin{lemme}\label{lem:le256} Nous avons, lorsque $\ N\to\infty,\ $
\begin{equation}
\left\vert J_{N}^{\prime}\right\vert=O_{\mathbb{P}}\left(\frac{1}{\sqrt{N}}\right)\ \mbox{ et }\ \left\vert
L_{N}^{\prime}\right\vert=O_{\mathbb{P}}\left(\frac{1}{\sqrt{N}}\right).\label{eq:iaaaaaad}
\end{equation}
\end{lemme}
\noindent{\bf Démonstration.} Nous observons tout d'abord, par l'intermédiaire de~(\ref{eq:iaaaaaaa}), que, pour $\ \displaystyle 0\leq w\leq
1,\ $

{\small
\begin{flalign*}
&\hspace{0cm}\int_{0}^{1}q_{3}^{(m)}(1-v)\ q_{2}^{(m)}(1-w)\ \mathbf{B}(v,w)\ \ud v=q_{2}^{(m)}(1-w)\int_{0}^{1}q_{3}^{(m)}(1-v)\ \mathbf{B}(v,w)\ \ud v\\
&\hspace{0cm}=q_{2}^{(m)}(1-w)\left\{\int_{0}^{1}q_{3}^{(m)}(1-v)\ {\mathbf{B}}_{[0]}(v,w)\ \ud v+\mathbf{B}(1,w)\int_{0}^{1}v\ q_{3}^{(m)}(1-v)\ \ud
v\right.\\
&\hspace{2,5cm}\left.+w\int_{0}^{1}q_{3}^{(m)}(1-v)\ \mathbf{B}(v,1)\ \ud v\right\}.
\end{flalign*}
}

\noindent En faisant le changement de variable $\ \displaystyle v\longrightarrow 1-v\ $ et en posant $\ v=G_{3}^{(m)}(x)\ $ et $\ \displaystyle
q_{3}^{(m)}(v)\ \ud v=\ud x,\ $ on obtient, pour $\ \displaystyle 0\leq w\leq 1,\ $ l'égalité

{\small
\begin{displaymath}
\int_{0}^{1}q_{3}^{(m)}(1-v)\ {\mathbf{B}}_{[0]}(v,w)\ \ud v=\int_{0}^{1}q_{3}^{(m)}(v)\ {\mathbf{B}}_{[0]}(1-v,w)\ \ud
v=\int_{0}^{1}{\mathbf{B}}_{[0]}(1-G_{3}^{(m)}(x),w)\ \ud x.
\end{displaymath}
}

\noindent Observons que

\begin{displaymath}
\left\{{\mathbf{B}}_{[0]}(v,w):0\leq v,w\leq 1\right\}\stackrel{d}{=}\left\{{\mathbf{B}}_{[0]}(1-v,w):0\leq v,w\leq 1\right\}.
\end{displaymath}

\noindent De plus, en faisant usage de~(\ref{eq:iaaaaaae}), avec $\ \displaystyle F=G_{3}^{(m)}\ $ qui est la f.r. d'une loi $\ \beta_{m,m}\ $
de variance $\ \displaystyle \sigma_{G_{3}^{(m)}}^{2}=\frac{1}{4(2m+1)},\ $ nous voyons que le processus $\ \left\{B_{0}(w):0\leq w\leq
1\right\}\ $ défini en posant, pour $\ \displaystyle 0\leq w\leq 1,\ $

{\footnotesize
\begin{equation}
\hspace{-0,5cm}\frac{1}{2\sqrt{2m+1}}\ B_{0}(w)=\int_{0}^{1}\ q_{3}^{(m)}(1-v)\ {\mathbf{B}}_{[0]}(v,w)\ \ud
v=\int_{0}^{1}{\mathbf{B}}_{[0]}(G_{3}^{(m)}(x),w)\ \ud x\label{eq:iaaaaaaf}
\end{equation}
}

\noindent est un pont Brownien. De même, en faisant un changement de variable\\*
$\displaystyle v\longrightarrow 1-v\ $ et en posant $\ v=G_{3}^{(m)}(x)\ $ et $\ \displaystyle q_{3}^{(m)}(v)\ \ud v=\ud x,\ $ on obtient
l'égalité

\begin{displaymath}
\int_{0}^{1}q_{3}^{(m)}(1-v)\ \mathbf{B}(v,1)\ \ud v=\int_{0}^{1}q_{3}^{(m)}(v)\ \mathbf{B}(1-v,1)\ \ud
v=\int_{0}^{1}\mathbf{B}(1-G_{3}^{(m)}(x),1)\ \ud x.
\end{displaymath}

\noindent Observons que

\begin{displaymath}
\left\{\mathbf{B}(v,1):0\leq v\leq 1\right\}\stackrel{d}{=}\left\{\mathbf{B}(1-v,1):0\leq v\leq 1\right\}.
\end{displaymath}

\noindent De plus, en faisant usage de~(\ref{eq:iaaaax}), avec $\ \displaystyle F=G_{3}^{(m)}\ $ qui est la f.r. d'une loi $\ \beta_{m,m}\ $ de
variance $\ \displaystyle \sigma_{G_{3}^{(m)}}^{2}=\frac{1}{4(2m+1)},\ $ on obtient,

\begin{displaymath}
\int_{0}^{1}q_{3}^{(m)}(1-v)\ \mathbf{B}(v,1)\ \ud v\stackrel{d}{=}\int_{0}^{1}\mathbf{B}(G_{3}^{(m)}(x),1)\ \ud
x\stackrel{d}{=}N\left(0,\frac{1}{4(2m+1)}\right).
\end{displaymath}

\noindent Une application du lemme~\ref{lem:le253}, montre que les ponts Browniens $\ \displaystyle \left\{B_{0}(w):\right. $
$\displaystyle\left.0\leq w\leq 1\right\}$ et $\ \displaystyle \left\{B_{1}(w):=\mathbf{B}(1,w):0\leq w\leq 1\right\}\ $
dans~(\ref{eq:iaaaaaaa}) et~(\ref{eq:iaaaaaaf}) ainsi que la variable $\ \displaystyle Z\ $ définie par

\begin{displaymath}
\frac{1}{2\sqrt{2m+1}}\ Z:=\int_{0}^{1}q_{3}^{(m)}(1-v)\ \mathbf{B}(v,1)\ \ud v\stackrel{d}{=}N\left(0,\frac{1}{4(2m+1)}\right)
\end{displaymath}

\noindent sont mutuellement indépendantes. Ceci avec~(\ref{eq:iaaak}) aussi montre que le processus $\ \displaystyle \left\{B_{3}(w):0\leq w\leq
1\right\},\ $ défini par l'identité, pour $\ \displaystyle 0\leq w\leq 1,\ $

\begin{displaymath}
\sqrt{\frac{m+1}{2(2m+1)}}\ B_{3}(w)=\frac{1}{2\sqrt{2m+1}}\ B_{0}(w)+\frac{1}{2}\ B_{1}(w),
\end{displaymath}

\noindent est un pont brownien indépendant de $\ \displaystyle Z\stackrel{d}{=}N(0,1).\ $ Par tout cela, nous avons l'égalité

{\small
\begin{flalign*}
&\int_{0}^{1}q_{3}^{(m)}(1-v)\ q_{2}^{(m)}(1-w)\ \mathbf{B}(v,w)\ \ud v=q_{2}^{(m)}(1-w)\left\{\sqrt{\frac{m+1}{2(2m+1)}}\ B_{3}(w)\right.\\
&\hspace{8,5cm}\left.+\frac{1}{2\sqrt{2m+1}}\ w\ Z\right\}.
\end{flalign*}
}

\noindent En faisant usage du théorème de Fubini, nous déduisons de ce dernier résultat que

\begin{flalign}
J_{N}^{\prime}&=\int_{0}^{1/N}\left\{\int_{0}^{1}q_{3}^{(m)}(1-v)\ q_{2}^{(m)}(1-w)\ \mathbf{B}(v,w)\ \ud v\right\}\ud w\label{eq:iaaaaaag}\\
&=\sqrt{\frac{m+1}{2(2m+1)}}\int_{0}^{1/N}q_{2}^{(m)}(1-w)\ B_{3}(w)\ \ud w\nonumber\\
&\quad+\frac{Z}{2\sqrt{2m+1}}\int_{0}^{1/N}w\ q_{2}^{(m)}(1-w)\ \ud w=:J_{N;1}^{\prime}+J_{N;2}^{\prime}.\nonumber
\end{flalign}

\noindent Par~(\ref{eq:iaak}), nous voyons que $\ \displaystyle w\ q_{2}^{(m)}(1-w)\to 1,\ $ lorsque $\ \displaystyle w\downarrow 0.\ $ Donc,
lorsque $\ \displaystyle w\downarrow 0,\ $

{\small
\begin{flalign}
J_{N;2}^{\prime}=\frac{Z}{2\sqrt{2m+1}}\int_{0}^{1/N}w\ q_{2}^{(m)}(1-w)\ \ud w&=(1+o(1))\ \frac{Z}{2N\sqrt{2m+1}}\label{eq:iaaaaaah}\\
&=O_{\mathbb{P}}\left(\frac{1}{N}\right).\nonumber
\end{flalign}}
\noindent En se rappelant de~(\ref{eq:iaaaaaac}) et~(\ref{eq:iaaaaaag}), nous voyons également que

{\footnotesize
\begin{equation}
J_{N;1}^{\prime}=\sqrt{\frac{m+1}{2(2m+1)}}\int_{0}^{1/N}q_{2}^{(m)}(1-w)\ B_{3}(w)\ \ud
w\stackrel{d}{=}\sqrt{\frac{m+1}{2(2m+1)}}L_{N}^{\prime}=N(0,\sigma_{N}^{2}),\label{eq:iaaaaaai}
\end{equation}
}
\noindent où nous utilisons le changement de variables $\ \displaystyle (w,v)=(s/N,t/N)\ $ pour écrire

{\small
\begin{flalign*}
\hspace{0cm}N\ \sigma^{2}_{N}&=N\ \mathbb{E}\left\{\left(\sqrt{\frac{m+1}{2(2m+1)}}\int_{0}^{1/N}q_{2}^{(m)}(1-w)\ B_{3}(w)\ \ud
w\right)^2\right\}\\
\hspace{0cm}&=\frac{N(m+1)}{2(2m+1)}\int_{0}^{1/N}\int_{0}^{1/N}(w\land v-wv)\ q_{2}^{(m)}(1-w)\ q_{2}^{(m)}(1-v)\ \ud w\ \ud v\\
\hspace{0cm}&=\frac{N(m+1)}{2(2m+1)}\int_{0}^{1}\int_{0}^{1}\left(\frac{s\land t}{st}-\frac{1}{N}\right)\left\{\frac{s}{N}\
q_{2}^{(m)}\left(1-\frac{s}{N}\right)\frac{t}{N}\ q_{2}^{(m)}\left(1-\frac{t}{N}\right)\right\}\ud s\ \ud t.
\end{flalign*}
}

\noindent En utilisant à nouveau le fait que, par~(\ref{eq:iaak}), $\ \displaystyle w\ q_{2}^{(m)}(1-w)\to 1,\ $ lorsque $\ \displaystyle
w\downarrow 0,\ $ nous déduisons de cette dernière égalité que, lorsque $\ \displaystyle N\to\infty,\ $

\begin{displaymath}
N\ \sigma^{2}_{N}=\frac{m+1}{2(2m+1)}\int_{0}^{1}\int_{0}^{1}\frac{s\land t}{st}\ \ud s\ \ud t-\frac{1+o(1)}{N}\ \frac{m+1}{2(2m+1)}.
\end{displaymath}

\noindent Un calcul facile montre que

\begin{flalign*}
\int_{0}^{1}\int_{0}^{1}\frac{s\land t}{st}\ \ud s\ \ud t&=\int_{0}^{1}\left\{\frac{1}{t}\int_{0}^{t}\ud s+\int_{t}^{1}\frac{\ud
s}{s}\right\}\ud t\\
&=\int_{0}^{1}\left\{1-\log t\right\}\ud t=1-\left[t\log t-t\right]_{0}^{1}=2.
\end{flalign*}

\noindent Ceci, aussi, montre que $\ \displaystyle N\ \sigma_{N}^{2}\to \frac{m+1}{2m+1},\ $ lorsque $\ N\to\infty.$
En se rappelant de~(\ref{eq:iaaaaaag}),~(\ref{eq:iaaaaaah}) et~(\ref{eq:iaaaaaai}), nous déduisons facilement~(\ref{eq:iaaaaaad}) de cette
propriété.$\Box$\vskip5pt

\subsection{Approximations fortes}\label{ma-soussectionj}
\noindent Le fait suivant est une version du théorème 2.3 de Castelle et Laurent-Bonvalot\cite{MR1900762} (voir par exemple
Tusn\'ady\cite{MR0443045}, et le théorème 1.1 de Castelle\cite{MR1632841}). Posons $\ \displaystyle\log_{+} v=\log(v\vee e),\ $ pour $\
\displaystyle v\in\mathbb{R},\ $ et posons $\ \displaystyle\left\vert\left\vert f\right\vert\right\vert=\sup_{z\in A}\left\vert
f(z)\right\vert,\ $ pour la norme-sup d'une fonction $\ \displaystyle f\ $ bornée, définie sur $\ \displaystyle A\ .$
\begin{fait}\label{fai:fa201}Sur un espace de probabilité convenable $\ \displaystyle (\Omega,\mathcal{A},\mathbb{P}),\ $ il est possible de
construire une suite de vecteurs aléatoires $\ \displaystyle \left\{\left(V_{k,N},W_{k,N}\right):0\leq k\leq N,\ N\geq1\right\}\ $ uniformément distribués sur $\ \displaystyle[0,1]^2\ $ et une suite de ponts Browniens bivariés $\ \left\{\mathbf{B}_{N}(v,w):0\leq
v,w\leq 1\right\}\ $ telle que la propriété suivante soit vérifiée. Pour des constantes appropriées $\ \displaystyle a>0,\ b>0\ \mbox{ et }\
c>0,\ $ nous avons, pour tout $\ \displaystyle x\geq 0\ \mbox{ et }\ N\geq 1,$
\begin{equation}
\mathbb{P}\left(\left\vert\left\vert\alpha_{N}-\mathbf{B}_{N}\right\vert\right\vert\geq \frac{\log_{+} N\left(a\log_{+}
N+x\right)}{\sqrt{N}}\right)\leq b\ e^{-cx}.\label{eq:iaaaaaaj}
\end{equation}
\end{fait}
\noindent Comme une conséquence évidente de~(\ref{eq:iaaaaaaj}), nous avons, lorsque $\ \displaystyle N\to\infty,\ $
\begin{equation}
\left\vert\left\vert\alpha_{N}-\mathbf{B}_{N}\right\vert\right\vert=O_{\mathbb{P}}
\left(\frac{(\log N)^2}{\sqrt{N}}\right).\label{eq:iaaaaaak}
\end{equation}
\noindent Compte tenu de~(\ref{eq:iaaaaai}),~(\ref{eq:iaaaaaj}) et~(\ref{eq:iaaaaak}), nous allons examiner ci-dessous, les processus définis
par
\begin{equation}
B_{N}(v)=\mathbf{B}_{N}(v,1)=-\mathbf{B}_{N}^{\ast}(1-v,1)=-B_{[1];N}^{\ast}(1-v),\label{eq:iaaaaaal}
\end{equation}
\noindent pour $\ \displaystyle 0\leq v\leq 1,\ $ où
\begin{flalign}
&\mathbf{B}_{N}^{\ast}(1-v,1-w)=\mathbf{B}_{N}(v,w)-\mathbf{B}_{N}(v,1)
-\mathbf{B}_{N}(1,w),\label{eq:iaaaaaam}\\
&\mathbf{B}_{[0];N}^{\ast}(v,w)=\mathbf{B}_{N}^{\ast}(v,w)-v\mathbf{B}_{N}^{\ast}(1,w)
-w\mathbf{B}_{N}^{\ast}(v,1),\label{eq:iaaaaaan}\\
&B_{[1];N}^{\ast}(v)
=\mathbf{B}_{N}^{\ast}(v,1)\ \mbox{ et }B_{[2];N}^{\ast}(w)=\mathbf{B}_{N}^{\ast}(1,w),\label{eq:iaaaaaao}
\end{flalign}
\noindent pour $\ \displaystyle 0\leq v,w\leq 1.$\vskip5pt
\noindent Le fait suivant est un cas particulier des principes d'invariance forts de Koml\'os, Major et Tusn\'ady\cite{MR0375412,MR0402883} pour
les sommes partielles. Nous nous référons au lemme A1 de Berkes et Philipp\cite{MR515811} pour des résultats qui nous permettent de combiner ces
constructions sur le même espace de probabilité.
\begin{fait}\label{fai:fa202} Sur un espace de probabilité $\ \displaystyle (\Omega,\mathcal{A},\mathbb{P}),\ $ soit $\ \displaystyle
\left\{\gamma_{i,n}:1\leq i\leq n\right\}\ $ qui désigne une suite de v.a. de loi $\ \displaystyle \Gamma(m,1),\ $ $\ \displaystyle
\gamma_{i,n}\stackrel{d}{=}\Gamma(m,1)\ $ avec $\ \displaystyle \mathbb{E}\left(\gamma_{i,n}\right)=m\ $ et $\ \displaystyle{\rm
Var}\left(\gamma_{i,n}\right)=m\ $ telle que, pour chaque $\ \displaystyle n\geq 1,\ $ $\ \displaystyle \gamma_{1,n},\ldots,\gamma_{n,n}\ $ sont
indépendantes. Alors sur une version convenablement élargie de $\ \displaystyle (\Omega,\mathcal{A},\mathbb{P}),\ $ il est possible de définir
une suite $\ \left\{\displaystyle W_{n}(t):t\geq 0,\ n\geq 1\right\} $ de processus de Wiener, telle que, pour chaque $\ n\geq 1\ \mbox{ et }\
x\geq 0,\ $
\begin{equation}
\mathbb{P}\left(\max_{1\leq j\leq n}\left\vert\sum_{i=1}^{j}\gamma_{i,n}-jm-\sqrt{m}\ W_{n}(j)\right\vert\geq x+A\ \log n\right)\leq B\
e^{-cx},\label{eq:iaaaaaap}
\end{equation}
\noindent où $\ \displaystyle A>0,\ B>0\ \mbox{ et }\ C>0\ $ sont des constantes universelles.
\end{fait}
\noindent Une conséquence facile du fait~\ref{fai:fa202} est que nous pouvons poser
{\small
\begin{equation}
\frac{1}{\sqrt{n}}\sum_{i=1}^{n}\left(\gamma_{i,n}-m\right)-n^{-1/2}\sqrt{m}\ W_{n}(n)=O_{\mathbb{P}}\left(\frac{\log n}{\sqrt{n}}\right),\
\mbox{ lorsque }\ n\to\infty.\label{eq:iaaaaaaq}
\end{equation}}

\subsection{Preuve du théorème 3.1}\label{ma-soussectionk}
\noindent Nous avons maintenant en main tous les ingrédients nécessaires à la preuve du théorème~\ref{thm:teo11111}. Nous posons $\ \displaystyle (\Omega,\mathcal{A},\mathbb{P})\ $ et $\ \displaystyle \mathbf{B}_{N}(.,.)\ $ comme dans le fait~\ref{fai:fa201}. Nous posons de
plus $\ \displaystyle B_{N}(.),\ \mathbf{B}_{N}^{\ast}(.,.),\ \mathbf{B}_{[0];N}^{\ast}(.,.),\ B_{[1];N}^{\ast}(.)$ $\displaystyle\mbox{et }\
B_{[2];N}^{\ast}(.)\ $ définis comme dans~(\ref{eq:iaaaaaal})$-$(\ref{eq:iaaaaaao}). Compte tenu de~(\ref{eq:iaaaaai}), ~(\ref{eq:iaaaaaj})
et~(\ref{eq:iaaaaak}), nous posons de plus,
{\small
\begin{flalign}
&\phi_{N}:=\sqrt{\frac{2(2m+1)}{m}}\int\int_{[0,1]^2} q_{3}^{(m)}(1-v)\ q_{2}^{(m)}(1-w)\ \mathbf{B}_{[0];N}^{\ast}(v,w)\ \ud v\
\ud w\label{eq:iaaaaaar}\\
&\mbox{et}\hspace{3cm}\psi_{N}=\frac{1}{\sqrt{2m}}\int_{0}^{1}q_{2}^{(m)}(1-w)\ B_{[2];N}^{\ast}(w)\ \ud w.\label{eq:iaaaaaat}
\end{flalign}}
\noindent Comme suit facilement de~(\ref{eq:iaaaaaak}), nous avons, lorsque $\ N\to\infty,\ $
\begin{flalign}
\left\vert\left\vert\alpha_{N;1}-B_{N}\right\vert\right\vert&=\sup_{0\leq v\leq
1}\left\vert\alpha_{N;1}(v,1)-\mathbf{B}_{N}(v,1)\right\vert\label{eq:iaaaaaau}\\
&\leq\left\vert\left\vert\alpha_{N}-\mathbf{B}_{N}\right\vert\right\vert=O_{\mathbb{P}}\left(\frac{\left(\log
N\right)^2}{\sqrt{N}}\right),\nonumber
\end{flalign}
\noindent qui est~(\ref{eq:iaaaaax}). Les arguments ci-dessous, capturés dans la proposition suivante complètent la preuve du
théorème~\ref{thm:teo11111}.
\begin{prop}\label{propos:propo256} Sous la notation ci-dessus, nous avons, lorsque $\ \displaystyle N\to\infty,\ $
{\footnotesize
\begin{equation}
\left\vert (N+1)^{1/2}\Delta_{N}-\phi_{N}\sqrt{\frac{2m}{2m+1}}
+4m\int_{0}^{1}B_{N}(G_{3}^{(m)}(x))\ \ud x\right\vert=O_{\mathbb{P}}\left(\frac{\left(\log N\right)^3}{N^{1/2}}\right),\label{eq:iaaaaaav}
\end{equation}}
\begin{equation}
\hspace{-1cm}et\hspace{1cm}\left\vert (N+1)^{1/2}\ \Theta_{N}-\psi_{N}\sqrt{2m}\right\vert
=O_{\mathbb{P}}\left(\frac{\left(\log N\right)^3}{N^{1/2}}\right),\label{eq:iaaaaaaw}
\end{equation}
\noindent où, pour chaque $\ N\geq 1,\ $ le pont Brownien $\ \displaystyle \left\{B_{N}(v):0\leq v\leq 1\right\} $ et les v.a. $\ \displaystyle \phi_{N}\stackrel{d}{=}N(0,1)\ $ et $\
\displaystyle\psi_{N}\stackrel{d}{=}N(0,1)\ $ sont indépendants.
\end{prop}
\noindent{\bf Démonstration.} Nous déduisons de~(\ref{eq:iaaad}) et~(\ref{eq:iaaaaal}) que
{\small
\begin{flalign*}
A_{N}:&=\left\vert (N+1)^{1/2}\Delta_{N}
-\int\int_{[0,1]^2}\left(2\ Q_{3}^{(m)}(v)-1\right)\ Q_{2}^{(m)}(w)\ {\mathbf{B}}_{N}(\ud v,\ud w)\right\vert\\
&=\left\vert (N+1)^{1/2}\Delta_{N}-2\int\int_{[0,1]^2} q_{3}^{(m)}(1-v)\ q_{2}^{(m)}(1-w)\ \mathbf{B}^{\ast}_{N}(v,w)\ \ud v\ \ud w\right.\\
&\hspace{4cm}\left.+\int_{0}^{1}q_{2}^{(m)}(1-w)\ \mathbf{B}_{N}^{\ast}(1,w)\ \ud w\right\vert\\
&=\left\vert2\int\int_{[0,1]^2} q_{3}^{(m)}(1-v)\ q_{2}^{(m)}(1-w)\ \left\{\alpha_{N}^{\ast}(v,w)-\mathbf{B}^{\ast}_{N}(v,w)\right\}\ \ud v\ \ud
w\right.\\
&\quad\quad\left.-\int_{0}^{1}q_{2}^{(m)}(1-w)\ \left\{\alpha_{N}^{\ast}(1,w)-\mathbf{B}^{\ast}_{N}(1,w)\right\}\ \ud w\right\vert.
\end{flalign*}}
\noindent Compte tenu de~(\ref{eq:iaah}),~(\ref{eq:iaai}), que $\ \displaystyle v\mapsto q_{3}^{(m)}(1-v)\ $ est bornée sur $\
\displaystyle(0,1)\ $ et que plus $\ \displaystyle\forall\ v\in(0,1),\ 0<q_{3}^{(m)}(1-v)<1\ $ et en faisant usage de~(\ref{eq:iaaaaaak})
et~(\ref{eq:iaaaaaau}), nous déduisons de cette expression et de l'inégalité triangulaire que
\begin{flalign*}
A_{N}\leq &3\ \left\{\int_{1/N}^{1}q_{2}^{(m)}(1-w)\ \ud w\right\}\left\vert\left\vert\alpha_{N}^{\ast}
-\mathbf{B}_{N}^{\ast}\right\vert\right\vert\\
&+2\left\vert\int_{0}^{1}q_{3}^{(m)}(1-v)\left\{\int_{0}^{1/N}q_{2}^{(m)}(1-w)\ \alpha_{N}^{\ast}(v,w)\ \ud w\right\}\ud v\right\vert\\
&+2\left\vert\int_{0}^{1}q_{3}^{(m)}(1-v)\left\{\int_{0}^{1/N}q_{2}^{(m)}(1-w)\ \mathbf{B}_{N}^{\ast}(v,w)\ \ud w\right\}\ud v\right\vert\\
&+\left\vert\int_{0}^{1/N}q_{2}^{(m)}(1-w)\ \alpha_{N}^{\ast}(1,w)\ \ud w\right\vert\\
&+\left\vert\int_{0}^{1/N}q_{2}^{(m)}(1-w)\ \mathbf{B}_{N}^{\ast}(1,w)\ \ud w\right\vert.
\end{flalign*}
\noindent En combinant~(\ref{eq:iaam}) et~(\ref{eq:iaaaaaau}), nous voyons que, lorsque $\ \displaystyle N\to\infty,\ $
\begin{flalign}
&3\ \left\{\int_{1/N}^{1}q_{2}^{(m)}(1-w)\ \ud w\right\}\left\vert\left\vert\alpha_{N}^{\ast}
-\mathbf{B}_{N}^{\ast}\right\vert\right\vert\label{eq:iaaaaaay}\\
&=3\ Q_{2}^{(m)}\left(1-\frac{1}{N}\right)\ \left\vert\left\vert\alpha_{N}^{\ast}
-\mathbf{B}_{N}^{\ast}\right\vert\right\vert=O_{\mathbb{P}}\left(\frac{\left(\log N\right)^3}{N^{1/2}}\right).\nonumber
\end{flalign}Compte tenu de~(\ref{eq:iaaaq}),~(\ref{eq:iaaaae}),~(\ref{eq:iaaar}),~(\ref{eq:iaaaah}),~(\ref{eq:iaaaaaab}),~(\ref{eq:iaaaaaac})
et~(\ref{eq:iaaaaaad}),
nous voyons que, pour chaque $\ \displaystyle\varepsilon>0,\ $ lorsque $\ N\to\infty,\ $
{\footnotesize
\begin{flalign}
&2\left\vert\int_{0}^{1}q_{3}^{(m)}(1-v)\left\{\int_{0}^{1/N}q_{2}^{(m)}(1-w)\ \alpha_{N}^{\ast}(v,w)\ \ud w\right\}\ud
v\right\vert=O_{\mathbb{P}}\left(\frac{\left(\log N\right)^{1+\varepsilon}}{N^{1/2}}\right),\label{eq:iaaaaaaz}\\
&\left\vert\int_{0}^{1/N}q_{2}^{(m)}(1-w)\ \alpha_{N}^{\ast}(1,w)\ \ud w\right\vert=O_{\mathbb{P}}\left(\frac{\left(\log
N\right)^{1+\varepsilon}}{N^{1/2}}\right),\nonumber\\
&2\left\vert\int_{0}^{1}q_{3}^{(m)}(1-v)\left\{\int_{0}^{1/N}q_{2}^{(m)}(1-w)\ \mathbf{B}_{N}^{\ast}(v,w)\ \ud w\right\}\ud
v\right\vert=O_{\mathbb{P}}\left(\frac{1}{N^{1/2}}\right)\label{eq:iaaaaaaaa}\\
&\mbox{et}\nonumber\\
&\left\vert\int_{0}^{1/N}q_{2}^{(m)}(1-w)\ \mathbf{B}_{N}^{\ast}(1,w)\ \ud
w\right\vert=O_{\mathbb{P}}\left(\frac{1}{N^{1/2}}\right)\nonumber.
\end{flalign}}En posant $\ \displaystyle\varepsilon=1\ $ dans les inégalités ci-dessus, nous concluons que, lorsque $\ N\to\infty,\ $
\begin{displaymath}
A_{N}=O_{\mathbb{P}}\left(\frac{\left(\log N\right)^3}{N^{1/2}}\right),
\end{displaymath}qui est~(\ref{eq:iaaaaaav}). Ensuite, nous faisons usage de la proposition~\ref{propos:propo255} pour réécrire la déclaration
ci-dessus en
{\small
\begin{flalign*}
A_{N}&=\left\vert (N+1)^{1/2}\Delta_{N}
-2\int\int_{[0,1]^2} q_{3}^{(m)}(1-v)\ q_{2}^{(m)}(1-w)\ \mathbf{B}^{\ast}_{[0];N}(v,w)\ \ud v\ \ud w\right.\\
&\quad\quad\left.-4m\int_{0}^{1}q_{3}^{(m)}(1-v)\ B^{\ast}_{[1];N}(v)\ \ud v\right\vert\\
&=\left\vert (N+1)^{1/2}\Delta_{N}
-\phi_{N}\ \sqrt{\frac{2m}{2m+1}}+4m\int_{0}^{1}B_{N}(G_{3}^{(m)}(x))\ \ud x\right\vert=O_{\mathbb{P}}\left(\frac{\left(\log
N\right)^3}{N^{1/2}}\right),
\end{flalign*}}lorsque $\ \displaystyle N\to\infty.\ $ Par des arguments similaires, nous déduisons de~(\ref{eq:iaaae}) et~(\ref{eq:iaaaaam}) que,
lorsque $\ \displaystyle N\to\infty,\ $
{\footnotesize
\begin{flalign*}
C_{N}:&=\left\vert (N+1)^{1/2}\ \Theta_{N}
-\int\int_{[0,1]^2} Q_{2}^{(m)}(w)\ \mathbf{B}_{N}(\ud v,\ud w)\right\vert\\
&=\left\vert\int_{0}^{1}q_{2}^{(m)}(1-w)\ \left\{\alpha_{N}^{\ast}(1,w)-{\mathbf{B}}^{\ast}_{N}(1,w)\right\}\ \ud w\right\vert\\
&\leq \left\{\int_{1/N}^{1}q_{2}^{(m)}(1-w)\ \ud w\right\}\left\vert\left\vert\alpha_{N}^{\ast}
-\mathbf{B}_{N}^{\ast}\right\vert\right\vert+\left\vert\int_{0}^{1/N}q_{2}^{(m)}(1-w)\ \alpha_{N}^{\ast}(1,w)\ \ud w\right\vert\\
&\quad+\left\vert\int_{0}^{1/N}q_{2}^{(m)}(1-w)\ {\mathbf{B}}^{\ast}_{N}(1,w)\ \ud w\right\vert.
\end{flalign*}}Compte tenu de~(\ref{eq:iaaaaaay}),~(\ref{eq:iaaaaaaz}) et~(\ref{eq:iaaaaaaaa}), nous déduisons de cette dernière relation que,
lorsque $\ \displaystyle N\to\infty,\ $
\begin{displaymath}
C_{N}=O_{\mathbb{P}}\left(\frac{\left(\log N\right)^3}{N^{1/2}}\right).
\end{displaymath}Ceci combiné à~(\ref{eq:iaaaaaat}) et~(\ref{eq:iaaaaam}) dans la proposition~\ref{propos:propo255}, implique que, lorsque $\
\displaystyle N\to\infty,\ $
\begin{flalign*}
C_{N}&=\left\vert (N+1)^{1/2}\ \Theta_{N}
-\int_{0}^{1}q_{2}^{(m)}(1-w)\ B^{\ast}_{[2];N}(w)\ \ud w\right\vert\\
&=\left\vert (N+1)^{1/2}\ \Theta_{N}
-\psi_{N}\sqrt{2m}\right\vert=O_{\mathbb{P}}\left(\frac{\left(\log N\right)^3}{N^{1/2}}\right),
\end{flalign*}qui est~(\ref{eq:iaaaaaaw}). Une application de la proposition~\ref{propos:propo255} montre que $\ \displaystyle \phi_{N}\stackrel{d}{=}N(0,1),\ $ $\displaystyle\ \psi_{N}\stackrel{d}{=}N(0,1)\ $ et $\ \displaystyle \left\{B_{N}(v):0\leq v\leq 1\right\} $ sont indépendants. Ceci complète la preuve de la proposition~\ref{propos:propo256}.$\Box$\vskip5pt

\section{Rapports d'espacements}\label{ma-sectione}

\subsection{Faits de base}\label{ma-soussectionl}

\noindent Ci-dessous, nous supposons que $\ \displaystyle F_{X}(t)=F_{Y}(t)=t,\ $ pour $\ \displaystyle 0\leq t\leq 1,\ $ et soit, pour $\
\displaystyle0\leq k\leq N,$ $\displaystyle S_{k,N;X}^{(m)}\ $ et $\ \displaystyle S_{k,N;Y}^{(m)}\ $ connu comme dans~(\ref{eq:iaaaaaaar}).
Dans toute la suite, nous allons considérer que dans~(\ref{eq:iaaaaaaar}), pour tout $\ \displaystyle 0\leq k\leq N,\ $ $\ \displaystyle i_{k,n_1}=i_{0,n_1}+km,\ $ et $\ \displaystyle j_{k,n_2}=j_{0,n_2}+km,\ $ où $\ i_{0,n_{1}}=j_{0,n_{2}}=0,\ $  et allons donc travailler avec les paires de $m$-espacements disjoints, $\displaystyle S_{k,N;X}^{(m)}\ $ et $\ \displaystyle S_{k,N;Y}^{(m)},\ $ définies de la façon suivante.
\begin{flalign*}
&S_{k,N;X}^{(m)}=D_{i_{k,n_{1}},n_1;X}^{(m)}=X_{\left(k+1\right)m,n_1}-X_{km,n_1}\ \mbox{ et }\\
&S_{k,N;Y}^{(m)}=D_{j_{k,n_{2}},n_2;Y}^{(m)}=Y_{\left(k+1\right)m,n_2}-Y_{km,n_2},\ \mbox{ pour }\ k=0,\ldots,N.
\end{flalign*}Nous ferons usage du fait suivant qui est une conséquence du théorème de représentation des $1$-espacements uniformes (voir par exemple,
Pyke\cite{MR0216622}). Rappelons les définitions~(\ref{eq:iaaaa}) et~(\ref{eq:iaaab}) de $\ \displaystyle \Delta_{N}\ $ et $\ \displaystyle
\Theta_{N}.\ $ Posons, pour tout entier $\ m\ $ fixé tel que $\ \displaystyle 1\leq m<n_{1}\land n_{2},\ N_{1}=\left\lfloor(n_{1}+1)/m\right\rfloor-1\ $ et $\
\displaystyle N_{2}=\left\lfloor(n_{2}+1)/m\right\rfloor-1\ $ comme dans~(\ref{eq:iaaaaaaaaf}) et posons aussi $\ \displaystyle P=P(N)\ $ et $\ \displaystyle Q=Q(N)\ $ comme dans~(\ref{eq:iaaaaaaaag}).
\begin{fait}\label{fai:fa203}  Pour chaque $\ \displaystyle 1\leq N<N_{1}\land N_{2},\ $ il existe deux suites indépendantes $\ \displaystyle
\left\{\zeta_{\ell,N}:\ell\geq 1\right\}\ $ et $\ \displaystyle \left\{\xi_{\ell,n}:\ell\geq 1\right\}\ $ de v.a. exponentiellement distribuées avec une moyenne 1,  telles que les relations suivantes soient vérifiées. Posons

{\small
\begin{flalign}
&\hspace{-1cm}\mathcal{T}_{N;X}=\sum_{k=0}^{N}\sum_{\ell=km+1}^{(k+1)m}\zeta_{\ell,N}=\sum_{k=0}^{N}Z_{k,N}=\frac{N+1}{2}
\left\{\Theta_{N}+2m+\Delta_{N}\right\},
\label{eq:iaaaaaaaah}\\
&\hspace{-1cm}\mathcal{R}_{N;X}=\sum_{k=N+1}^{N+P}\sum_{\ell=km+1}^{(k+1)m}\zeta_{\ell,N}=\sum_{k=N+1}^{N+P}Z_{k,N},\nonumber\\
&\hspace{-1cm}\mathcal{T}_{N;Y}=\sum_{k=0}^{N}\sum_{\ell=km+1}^{(k+1)m}\xi_{\ell,N}=\sum_{k=0}^{N}Z_{k,N}^{\prime}
=\frac{N+1}{2}\left\{\Theta_{N}+2m-\Delta_{N}\right\},
\label{eq:iaaaaaaaai}\\
&\hspace{-1cm}\mathcal{R}_{N;Y}=\sum_{k=N+1}^{N+Q}\sum_{\ell=km+1}^{(k+1)m}\xi_{\ell,N}=\sum_{k=N+1}^{N+Q}Z_{k,N}^{\prime},\nonumber
\end{flalign}}où, pour tout $\displaystyle k\geq 0,\ Z_{k,N}=\sum_{\ell=km+1}^{(k+1)m}\zeta_{\ell,N}\ $ et $\ \displaystyle
Z_{k,N}^{\prime}=\sum_{\ell=km+1}^{(k+1)m}\xi_{\ell,N}\ $ sont comme dans~(\ref{eq:iiaaaaaaaaaf}).\vskip5pt
\noindent Nous avons, pour tout $\ \displaystyle k=0,\ldots,N,\ $
\begin{flalign}
&\hspace{-1cm}\displaystyle
S_{k,N;X}^{(m)}\stackrel{d}{=}\frac{\displaystyle\sum_{\ell=km+1}^{(k+1)m}\zeta_{\ell,n}}{\mathcal{T}_{N;X}+\mathcal{R}_{N;X}}
\stackrel{d}{=}\frac{Z_{k,N}}{\mathcal{T}_{N;X}+\mathcal{R}_{N;X}}\ \mbox{ et }\label{eq:iaaaaaaab}\\
&\hspace{-1cm}S_{k,N;Y}^{(m)}\stackrel{d}{=}\frac{\displaystyle\sum_{\ell=km+1}^{(k+1)m}\xi_{\ell,n}}{\mathcal{T}_{N;Y}+\mathcal{R}_{N;Y}}
\stackrel{d}{=}\frac{Z_{k,N}^{\prime}}{\mathcal{T}_{N;Y}+\mathcal{R}_{N;Y}}.\nonumber
\end{flalign}
\end{fait}\noindent Compte tenu de~(\ref{eq:iaaaaaaat}) et en faisant usage du fait~\ref{fai:fa203}, nous observons que les égalités d'événements suivantes sont vérifiées. Posons, pour tout $\ \displaystyle t\in[0,1],\ $
\begin{displaymath}
\tau_{N}(t)=\left\{1+\left(\frac{1}{t}-1\right)\frac{\left\{\mathcal{T}_{N;Y}+\mathcal{R}_{N;Y}\right\}/(N+Q+1)}
{\left\{\mathcal{T}_{N;X}+\mathcal{R}_{N;X}\right\}/(N+P+1)}\right\}^{-1}.
\end{displaymath}En se rappelant des définitions~(\ref{eq:iaao}) et~(\ref{eq:iaaq}), de respectivement $\ \displaystyle R_{k,N}\ $ et $\ \displaystyle
V_{k,N},\ $ nous voyons que, pour tout $\ 0\leq k\leq N,\ $ et pour tout $\ \displaystyle t\in[0,1],\ $
{\scriptsize
\begin{flalign}
\left\{R_{k;n_1,n_2}\leq t\right\}&=\left\{\frac{(N+P+1)\ S_{k,N;X}^{(m)}}{(N+P+1)\ S_{k,N;X}^{(m)}+(N+Q+1)\ S_{k,N;Y}^{(m)}}\leq
t\right\}\label{eq:iaaaaaaac}\\
&=\left\{\frac{S_{k,N;Y}^{(m)}}{S_{k,N;X}^{(m)}}\geq\left(\frac{1}{t}-1\right)\ \frac{N+P+1}{N+Q+1}\right\}\nonumber\\
&=\left\{\frac{Z_{k,N}^{\prime}}{Z_{k,N}}\geq\left(\frac{1}{t}-1\right)\frac{\left\{\mathcal{T}_{N;Y}+\mathcal{R}_{N;Y}\right\}/(N+Q+1)}
{\left\{\mathcal{T}_{N;X}+\mathcal{R}_{N;X}\right\}/(N+P+1)}\right\}\nonumber\\
&=\left\{\frac{Z_{k,N}^{\prime}}
{Z_{k,N}}\geq\frac{1}{\tau_{N}(t)}-1\right\}=\left\{\frac{Z_{k,N}}{Z_{k,N}+Z_{k,N}^{\prime}}\leq\tau_{N}(t)\right\}=\left\{R_{k,N}\leq
\tau_{N}(t)\right\}\nonumber\\
&=\left\{G_{3}^{(m)}(R_{k,N})\leq G_{3}^{(m)}(\tau_{N}(t))\right\}=\left\{V_{k,N}\leq
G_{3}^{(m)}(\tau_{N}(t))\right\}.\nonumber
\end{flalign}}Rappelons les définitions~(\ref{eq:iaaaaaaau}) et~(\ref{eq:iaaaaaaav}) de $\ \displaystyle H_{N;n_1,n_2}(.)\ $ et $\
\gamma_{N;n_1,n_2}(.).\ $ Compte tenu de~(\ref{eq:iaav}),~(\ref{eq:iaaw}) et~(\ref{eq:iaay}), nous déduisons de~(\ref{eq:iaaaaaaac}) que, pour tout $\ t\in[0,1],\ $
{\small
\begin{flalign}
&\gamma_{N;n_1,n_2}(t)=\left(N+1\right)^{1/2}\left(H_{N;n_1,n_2}(t)-H_{m}(t)\right)\label{eq:iaaaaaaad}\\
&=\left(N+1\right)^{1/2}\left(U_{N}(G_{3}^{(m)}(\tau_{N}(t)),1)-H_{m}(t)\right)\nonumber\\
&=\alpha_{N;1}(G_{3}^{(m)}(\tau_{N}(t)))+\left(N+1\right)^{1/2}\left(G_{3}^{(m)}(\tau_{N}(t))-H_{m}(t)\right)\nonumber\\
&=\alpha_{N;1}(G_{3}^{(m)}(\tau_{N}(t)))+\left(N+1\right)^{1/2}\left(G_{3}^{(m)}(\tau_{N}(t))-G_{3}^{(m)}(t)\right),\nonumber
\end{flalign}}puisque, pour tout $\ t\in\mathbb{R},\ $
\begin{displaymath}
H_{m}(t)=G_{3}^{(m)}(t).
\end{displaymath}
\begin{lemme}
\label{lem:le257}Nous avons, uniformément en $\ \displaystyle t\in[0,1],\ $ et lorsque $\ \displaystyle N\to\infty,\ $
{\footnotesize
\begin{flalign}
&\gamma_{N;n_1,n_2}(t)=\alpha_{N;1}(G_{3}^{(m)}(\tau_{N}(t)))+\frac{(2m-1)!}{\left((m-1)!\right)^2}\
\left(t(1-t)\right)^{m-1}\left(N+1\right)^{1/2}\left(\tau_{N}(t)-t\right)\label{eq:iaaaaaaah}\\
&\quad\quad\quad\quad\quad\quad+O_{\mathbb{P}}\left(\frac{1}{\sqrt{N}}\right).\nonumber
\end{flalign}}
\end{lemme}\noindent Posons, pour tout $\ \displaystyle N\geq 1,\ $
\begin{displaymath}
\mathcal{Q}_{N}=\frac{\left\{\mathcal{T}_{N;Y}+\mathcal{R}_{N;Y}\right\}/(N+Q+1)}{\left\{\mathcal{T}_{N;X}+\mathcal{R}_{N;X}\right\}/(N+P+1)}-1.
\end{displaymath}Observons que
\begin{equation}
\tau_{N}(t)-t=-t(1-t)\ \mathcal{Q}_{N}\left\{1+(1-t)\mathcal{Q}_{N}\right\}^{-1}.\label{eq:iaaaaaaai}
\end{equation}Posons, de plus
{\scriptsize
\begin{equation}
\mathcal{D}_{N}=\mathcal{Q}_{N}+\frac{1}{2m}\left\{\frac{N+1}{N+P+1}+\frac{N+1}{N+Q+1}\right\}\Delta_{N}
+\frac{1}{2m}\left\{\frac{N+1}{N+P+1}-\frac{N+1}{N+Q+1}\right\}\Theta_{N}.\label{eq:iaaaaaaaaj}
\end{equation}
}\begin{lemme}
\label{lem:le260} Supposons que $\ \displaystyle 0\leq c\leq d\leq \infty.\ $ Alors, lorsque $\ \displaystyle N\to\infty,\ $
{\footnotesize
\begin{equation}
\left(\sqrt{N+1}\ \Delta_{N},\ \sqrt{N+1}\ \Theta_{N},\ \sqrt{N+1}\ \mathcal{D}_{N}\right)\stackrel{d}{\longrightarrow}\left(\sqrt{2m}\
\mathcal{X},\ \sqrt{2m}\ \mathcal{Y},\ \sigma_{1}(c,d)\ \mathcal{Z}\right),\label{eq:iaaaaaaaak}
\end{equation}}où $\ \displaystyle\mathcal{X}\stackrel{d}{=}N(0,1),\ \mathcal{Y}\stackrel{d}{=}N(0,1),\ \mathcal{Z}\stackrel{d}{=}N(0,1)\ $ sont
indépendantes et
\begin{equation}
\sigma_{1}^{2}(c,d)= \left\{ \begin{array}{lll}
\displaystyle\frac{1}{m}\left\{\frac{c}{(1+c)^2}+\frac{d}{(1+d)^2}\right\} & \textrm{si\quad $ 0\leq c\leq d<\infty$}\\
\displaystyle\frac{1}{m}\ \frac{c}{(1+c)^2}& \textrm{si\quad $ 0\leq c<d=\infty$}\\
\displaystyle0 & \textrm{si\quad $ c=d=\infty$}
\end{array} \right..\label{eq:iaaaaaaaal}
\end{equation}De plus, nous avons, lorsque $\ \displaystyle N\to\infty,\ $
\begin{equation}
(N+1)^{1/2}\mathcal{Q}_{N}\stackrel{d}{\longrightarrow}N\left(0,\sigma_{2}^{2}(c,d)\right),\label{eq:iaaaaaaaam}
\end{equation}où
\begin{equation}
\sigma_{2}^{2}(c,d)= \left\{ \begin{array}{lll}
\displaystyle\frac{1}{m}\left\{\frac{1}{1+c}+\frac{1}{1+d}\right\} & \textrm{si\quad $\ 0\leq c\leq d<\infty$}\\
\displaystyle\frac{1}{m}\ \frac{1}{1+c}& \textrm{si\quad $\ 0\leq c<d=\infty$}\\
\displaystyle0 & \textrm{si\quad $c=d=\infty$}
\end{array} \right..\label{eq:iaaaaaaaan}
\end{equation}
\end{lemme}
\noindent{\bf Démonstration.} Nous supposerons ci-dessous que~(\ref{eq:iaaaaaaas}) est vérifié avec\\*
$\displaystyle 0<c,d<\infty.\ $ La preuve du lemme dans les autres cas étant très similaire, sera omise. Posons
\begin{flalign}
&\mathcal{T}_{N;X}=m(N+1)+\eta_{N}^{\prime}\sqrt{N+1},\ \mathcal{R}_{N;X}=m\ P+\eta_{N}^{\prime\prime}\sqrt{P}\label{eq:iaaaaaaaap}\\
&\mathcal{T}_{N;Y}=m(N+1)+v_{N}^{\prime}\sqrt{N+1},\ \mathcal{R}_{N;Y}=m\ Q+v_{N}^{\prime\prime}\sqrt{Q}.\label{eq:iaaaaaaaaq}
\end{flalign}Compte tenu de~(\ref{eq:iaaaaaaaah})$-$(\ref{eq:iaaaaaaaai}) et~(\ref{eq:iaaaaaaaap})$-$(\ref{eq:iaaaaaaaaq}), nous notons que pour un usage ultérieur,
\begin{equation}
(N+1)^{1/2}\ \Delta_{N}=\eta_{N}^{\prime}-v_{N}^{\prime}\ \mbox{ et }\ (N+1)^{1/2}\ \Theta_{N}=\eta_{N}^{\prime}+v_{N}^{\prime}.\label{eq:iaaaaaaaar}
\end{equation}Or, de~(\ref{eq:iaaaaaaaah}) et~(\ref{eq:iaaaaaaaap}), on a,
\begin{displaymath}
m(N+1)+\eta_{N}^{\prime}\sqrt{N+1}=\sum_{k=0}^{N}Z_{k,N}\ \mbox{ et }\ m\ P+\eta_{N}^{\prime\prime}\sqrt{P}=\sum_{k=N+1}^{N+P}Z_{k,N}
\end{displaymath}i.e,
\begin{displaymath}
\eta_{N}^{\prime}=\sqrt{N+1}\left(\frac{1}{N+1}\sum_{k=0}^{N}Z_{k,N}-m\right)\ \mbox{ et }\
\eta_{N}^{\prime\prime}=\sqrt{P}\left(\frac{1}{P}\sum_{k=N+1}^{N+P}Z_{k,N}-m\right).
\end{displaymath}De même, de~(\ref{eq:iaaaaaaaai}) et~(\ref{eq:iaaaaaaaaq}), on a,
\begin{displaymath}
m(N+1)+v_{N}^{\prime}\sqrt{N+1}=\sum_{k=0}^{N}Z_{k,N}^{\prime}\ \mbox{ et }\ m\ Q+v_{N}^{\prime\prime}\sqrt{Q}
=\sum_{k=N+1}^{N+Q}Z_{k,N}^{\prime}
\end{displaymath}i.e,
\begin{displaymath}
v_{N}^{\prime}=\sqrt{N+1}\left(\frac{1}{N+1}\sum_{k=0}^{N}Z_{k,N}^{\prime}-m\right)\ \mbox{ et }\
v_{N}^{\prime\prime}=\sqrt{Q}\left(\frac{1}{Q}\sum_{k=N+1}^{N+Q}Z_{k,N}^{\prime}-m\right).
\end{displaymath}Donc, on a, lorsque $\ \displaystyle N\to\infty,\ $ la normalité asymptotique jointe suivante,
\begin{equation}
\left( \begin{array}{c}
\eta_{N}^{\prime}\\
\eta_{N}^{\prime\prime}\\
v_{N}^{\prime}\\
v_{N}^{\prime\prime}
\end{array} \right)\stackrel{d}{\longrightarrow} N\left(\left( \begin{array}{c}
0\\
0\\
0\\
0
\end{array} \right),\Sigma\right),\label{eq:iaaaaaaaas}
\end{equation}où
\begin{equation}
\Sigma=\left[ \begin{array}{cccc}
A & E & F & G\\
E & B & H &I\\
F & H & C & J\\
G & I & J & D\end{array}\right],\label{eq:iaaaaaaaat}
\end{equation}avec
{\small
\begin{flalign*}
\hspace{0cm}&A=\lim_{N\to\infty}{\rm Var}\left\{\frac{1}{\sqrt{N+1}}\left(\sum_{k=0}^{N}Z_{k,N}-m(N+1)\right)\right\},\\
\hspace{0cm}&B=\lim_{N\to\infty}{\rm Var}\left\{\frac{1}{\sqrt{P}}\left(\sum_{k=N+1}^{N+P}Z_{k,N}-m\ P\right)\right\},\\
\hspace{0cm}&C=\lim_{N\to\infty}{\rm Var}\left\{\frac{1}{\sqrt{N+1}}\left(\sum_{k=0}^{N}Z_{k,N}^{\prime}-m(N+1)\right)\right\},\\
\hspace{0cm}&D=\lim_{N\to\infty}{\rm Var}\left\{\frac{1}{\sqrt{Q}}\left(\sum_{k=N+1}^{N+Q}Z_{k,N}^{\prime}-m\ Q\right)\right\},\\
\hspace{0cm}&E=\lim_{N\to\infty}{\rm Cov}\left\{\frac{1}{\sqrt{N+1}}\left(\sum_{k=0}^{N}Z_{k,N}-m(N+1)\right),\frac{1}{\sqrt{P}}\left(\sum_{k=N+1}^{N+P}Z_{k,N}-m\ P\right)\right\}\\
\hspace{0cm}&F=\lim_{N\to\infty}{\rm Cov}\left\{\frac{1}{\sqrt{N+1}}\left(\sum_{k=0}^{N}Z_{k,N}-m(N+1)\right),\frac{1}{\sqrt{N+1}}\left(\sum_{k=0}^{N}Z_{k,N}^{\prime}-m\ (N+1)\right)\right\},\\
\hspace{-1cm}&G=\lim_{N\to\infty}{\rm Cov}\left\{\frac{1}{\sqrt{N+1}}\left(\sum_{k=0}^{N}Z_{k,N}-m(N+1)\right),\frac{1}{\sqrt{Q}}\left(\sum_{k=N+1}^{N+Q}Z_{k,N}^{\prime}-m\ Q\right)\right\},\\
\hspace{-1cm}&H=\lim_{N\to\infty}{\rm Cov}\left\{\frac{1}{\sqrt{P}}\left(\sum_{k=N+1}^{N+P}Z_{k,N}-m\
P\right),\frac{1}{\sqrt{N+1}}\left(\sum_{k=0}^{N}Z_{k,N}^{\prime}-m(N+1)\right)\right\},\\
\hspace{-1cm}&I=\lim_{N\to\infty}{\rm Cov}\left\{\frac{1}{\sqrt{P}}\left(\sum_{k=N+1}^{N+P}Z_{k,N}-m\
P\right),\frac{1}{\sqrt{Q}}\left(\sum_{k=N+1}^{N+Q}Z_{k,N}^{\prime}-m\ Q\right)\right\}
\ \mbox{ et }\\
\hspace{-1cm}&J=\lim_{N\to\infty}{\rm Cov}\left\{\frac{1}{\sqrt{N+1}}\left(\sum_{k=0}^{N}Z_{k,N}^{\prime}-m(N+1)\right),\frac{1}{\sqrt{Q}}\left(\sum_{k=N+1}^{N+Q}Z_{k,N}^{\prime}-m\ Q\right)\right\}.
\end{flalign*}}Or
\begin{equation}
A=\lim_{N\to\infty}\left\{\frac{1}{N+1}\sum_{k=0}^{N}{\rm Var}(Z_{0,N})\right\}={\rm Var}(Z_{0,N})=m,\label{eq:iaaaaaaaau}
\end{equation}
\begin{equation}
B=\lim_{N\to\infty}\left\{\frac{1}{P}\sum_{k=N+1}^{N+P}{\rm Var}(Z_{0,N})\right\}={\rm Var}(Z_{0,N})=m,\label{eq:iaaaaaaaau}
\end{equation}de même, on a,
\begin{equation}
C=D=m.\label{eq:iaaaaaaaav}
\end{equation}De plus, comme les suites de v.a. $\ \displaystyle \left\{Z_{k,N}:k\geq 0\right\}\ \mbox{ et }\ \left\{Z_{k,N}^{\prime}:k\geq
0\right\}\ $ sont indépendantes, on a,
\begin{equation}
E=F=G=H=I=J=0.\label{eq:iaaaaaaaaw}
\end{equation}Et donc, de~(\ref{eq:iaaaaaaaau})$-$(\ref{eq:iaaaaaaaaw}), on en déduit que,
\begin{flalign}
&\Sigma=m\left[ \begin{array}{cccc}
\displaystyle 1 & 0 & 0 & 0 \\
\displaystyle 0 & 1 & 0 & 0  \\
\displaystyle 0 & 0 & 1 & 0\\
\displaystyle 0 & 0 & 0 & 1\end{array}\right].\label{eq:iaaaaaaaax}
\end{flalign}D'où, de~(\ref{eq:iaaaaaaaas}) et~(\ref{eq:iaaaaaaaax}), on en déduit que, lorsque $\ N\to\infty,\ $
\begin{flalign}
&\left( \begin{array}{c}
\eta_{N}^{\prime}\\
\eta_{N}^{\prime\prime}\\
v_{N}^{\prime}\\
v_{N}^{\prime\prime}
\end{array} \right)\stackrel{d}{\longrightarrow}N\left(\left( \begin{array}{c}
0\\
0\\
0\\
0
\end{array} \right),
m\left[ \begin{array}{cccc}
\displaystyle 1 & 0 & 0 & 0 \\
\displaystyle 0 & 1 & 0 & 0  \\
\displaystyle 0 & 0 & 1 & 0\\
\displaystyle 0 & 0 & 0 & 1
\end{array}\right]\right).\label{eq:iaaaaaaaay}
\end{flalign}Par suite, on en déduit que, lorsque $\ N\to\infty,\ $
\begin{equation}
\left(\eta_{N}^{\prime},\ \eta_{N}^{\prime\prime},\ v_{N}^{\prime},\ v_{N}^{\prime\prime}\right)\stackrel{d}{\longrightarrow}
\left(\eta^{\prime},\ \eta^{\prime\prime},\ v^{\prime},\ v^{\prime\prime}\right),\label{eq:iaaaaaaaay1}
\end{equation}où $\ \displaystyle \eta^{\prime},\ \eta^{\prime\prime},\ v^{\prime},\ v^{\prime\prime}\ $ désignent des v.a.
indépendantes de loi $\ \displaystyle N(0,m).\ $ Nous obtenons donc que, lorsque $\ \displaystyle N\to\infty,\ $
\begin{flalign*}
\hspace{0cm}\mathcal{Q}_{N}&=\left\{m+\frac{v_{N}^{\prime}\sqrt{N+1}+v_{N}^{\prime\prime}\sqrt{Q}}{N+Q+1}\right\}
\left\{m+\frac{\eta_{N}^{\prime}\sqrt{N+1}
+\eta_{N}^{\prime\prime}\sqrt{P}}{N+P+1}\right\}^{-1}-1\\
\hspace{0cm}&=\left\{m+\frac{v_{N}^{\prime}\sqrt{N+1}+v_{N}^{\prime\prime}\sqrt{Q}}{N+Q+1}
\right\}\left\{\frac{1}{m}-\frac{1}{m^2}\frac{\eta_{N}^{\prime}\sqrt{N+1}
+\eta_{N}^{\prime\prime}\sqrt{P}}{N+P+1}\right\}-1\\
\hspace{0cm}&=1-\frac{1}{m}\frac{\eta_{N}^{\prime}\sqrt{N+1}
+\eta_{N}^{\prime\prime}\sqrt{P}}{N+P+1}+\frac{1}{m}\frac{v_{N}^{\prime}\sqrt{N+1}
+v_{N}^{\prime\prime}\sqrt{Q}}{N+Q+1}\\
\hspace{0cm}&\quad-\frac{1}{m^2}\left(\frac{v_{N}^{\prime}\sqrt{N+1}
+v_{N}^{\prime\prime}\sqrt{Q}}{N+Q+1}\right)\left(\frac{\eta_{N}^{\prime}\sqrt{N+1}
+\eta_{N}^{\prime\prime}\sqrt{P}}{N+P+1}\right)-1\\
\hspace{0cm}&=-\frac{1}{m}\frac{\eta_{N}^{\prime}\sqrt{N+1}
+\eta_{N}^{\prime\prime}\sqrt{P}}{N+P+1}+\frac{1}{m}\frac{v_{N}^{\prime}\sqrt{N+1}
+v_{N}^{\prime\prime}\sqrt{Q}}{N+Q+1}\\
\hspace{0cm}&\quad-\frac{1}{m^2}\frac{v_{N}^{\prime}\eta_{N}^{\prime}(N+1)}{(N+Q+1)(N+P+1)}-\frac{1}{m^2}
\frac{v_{N}^{\prime}\eta_{N}^{\prime\prime}\sqrt{N+1}\sqrt{P}}{(N+Q+1)(N+P+1)}\\
\hspace{0cm}&\quad-\frac{1}{m^2}\frac{v_{N}^{\prime\prime}\eta_{N}^{\prime}\sqrt{Q}\sqrt{N+1}}{(N+Q+1)(N+P+1)}
-\frac{1}{m^2}\frac{v_{N}^{\prime\prime}\eta_{N}^{\prime\prime}\sqrt{Q}\sqrt{P}}{(N+Q+1)(N+P+1)}.
\end{flalign*}D'après~(\ref{eq:iaaaaaaaay1}), on a, lorsque $\ \displaystyle N\to\infty,\ $ $\ \displaystyle
\eta_{N}^{\prime}\stackrel{d}{\longrightarrow}\eta^{\prime},\ $ $\ \displaystyle
\eta_{N}^{\prime\prime}\stackrel{d}{\longrightarrow}\eta^{\prime\prime},\ $ $\ \displaystyle
v_{N}^{\prime}\stackrel{d}{\longrightarrow}v^{\prime},\ $ $\ \displaystyle
v_{N}^{\prime\prime}\stackrel{d}{\longrightarrow}v^{\prime\prime},\ $ et donc, on a, lorsque $\ \displaystyle N\to\infty,\ $
\begin{equation}
\eta_{N}^{\prime}=\eta_{N}^{\prime\prime}=v_{N}^{\prime}=v_{N}^{\prime\prime}=O_{\mathbb{P}}(1).\label{eq:iaaaaaaaay2}
\end{equation}D'où, on a, lorsque $\ \displaystyle
N\to\infty,\ $
\begin{flalign*}
\hspace{0cm}\mathcal{Q}_{N}&=-\frac{1}{m}\frac{\eta_{N}^{\prime}\sqrt{N+1}
+\eta_{N}^{\prime\prime}\sqrt{P}}{N+P+1}+\frac{1}{m}\frac{v_{N}^{\prime}\sqrt{N+1}
+v_{N}^{\prime\prime}\sqrt{Q}}{N+Q+1}\\
\hspace{0cm}&\quad+O_{\mathbb{P}}\left(\frac{N+1}{(N+Q+1)(N+P+1)}\right)
+O_{\mathbb{P}}\left(\frac{\sqrt{N+1}\sqrt{P}}{(N+Q+1)(N+P+1)}\right)\\
\hspace{0cm}&\quad+O_{\mathbb{P}}\left(\frac{\sqrt{Q}\sqrt{N+1}}{(N+Q+1)(N+P+1)}\right)
+O_{\mathbb{P}}\left(\frac{\sqrt{Q}\sqrt{P}}{(N+Q+1)(N+P+1)}\right).
\end{flalign*}Or, d'après~(\ref{eq:iaaaaaaaag}), on a, $\ \displaystyle P=P(N)=N_{1}-N\geq 0\ $
et $\ \displaystyle Q=Q(N)=N_{2}-N\geq 0,\ $ donc, on a, lorsque $\ \displaystyle N\to\infty,\ $
\begin{flalign*}
\hspace{0cm}\mathcal{Q}_{N}&=-\frac{1}{m}\frac{\eta_{N}^{\prime}\sqrt{N+1}
+\eta_{N}^{\prime\prime}\sqrt{P}}{N+P+1}+\frac{1}{m}\frac{v_{N}^{\prime}\sqrt{N+1}
+v_{N}^{\prime\prime}\sqrt{Q}}{N+Q+1}\\
\hspace{0cm}&\quad+O_{\mathbb{P}}\left(\frac{N_{1}-P+1}{(N_{1}-P+Q+1)(N_{1}+1)}\right)
+O_{\mathbb{P}}\left(\frac{\sqrt{N_{1}-P+1}\sqrt{N_{1}-N}}{(N_{1}+Q-P+1)(N_{1}+1)}\right)\\
\hspace{0cm}&\quad+O_{\mathbb{P}}\left(\frac{\sqrt{N_{2}-N}\sqrt{N_{2}-Q+1}}{(N_{2}+1)(N_{2}+P-Q+1)}\right)
+O_{\mathbb{P}}\left(\frac{\sqrt{N_{2}-N}\sqrt{N_{1}-N}}{(N_{1}+1)(N_{2}+1)}\right)\\
\hspace{0cm}&=-\frac{1}{m}\frac{\eta_{N}^{\prime}\sqrt{N+1}
+\eta_{N}^{\prime\prime}\sqrt{P}}{N+P+1}+\frac{1}{m}\frac{v_{N}^{\prime}\sqrt{N+1}
+v_{N}^{\prime\prime}\sqrt{Q}}{N+Q+1}\\
\hspace{0cm}&\quad+O_{\mathbb{P}}\left(\frac{1}{N_{1}+1}\right)
+O_{\mathbb{P}}\left(\frac{1}{N_{2}+1}\right)\\
\hspace{0cm}&=-\frac{1}{m}\frac{\eta_{N}^{\prime}\sqrt{N+1}
+\eta_{N}^{\prime\prime}\sqrt{P}}{N+P+1}+\frac{1}{m}\frac{v_{N}^{\prime}\sqrt{N+1}
+v_{N}^{\prime\prime}\sqrt{Q}}{N+Q+1}\\
\hspace{0cm}&\quad+O_{\mathbb{P}}\left(\frac{1}{N+P+1}\right)
+O_{\mathbb{P}}\left(\frac{1}{N+Q+1}\right).
\end{flalign*}Or, $\ \displaystyle \forall\ N\ \geq 1,\ $ on a, les inégalités suivantes,
\begin{flalign}
\hspace{0cm}&\frac{1}{\sqrt{N+P+1}}\leq\frac{\sqrt{N+1}+\sqrt{P}}{N+P+1}\leq \frac{\sqrt{3}}{\sqrt{N+P+1}}\label{eq:iaaaaaaaay3}\\
\hspace{0cm}&\mbox{et }\ \frac{1}{\sqrt{N+Q+1}}\leq\frac{\sqrt{N+1}+\sqrt{Q}}{N+Q+1}\leq \frac{\sqrt{3}}{\sqrt{N+Q+1}},\nonumber
\end{flalign}donc, d'après~(\ref{eq:iaaaaaaaay2}) et~(\ref{eq:iaaaaaaaay3}), on en déduit que, lorsque $\ \displaystyle N\to\infty,\ $
{\footnotesize
\begin{flalign}
\hspace{0cm}\mathcal{Q}_{N}&=\frac{1}{m}\frac{v_{N}^{\prime}\sqrt{N+1}+v_{N}^{\prime\prime}\sqrt{Q}}{N+Q+1}
-\frac{1}{m}\frac{\eta_{N}^{\prime}\sqrt{N+1}
+\eta_{N}^{\prime\prime}\sqrt{P}}{N+P+1}\label{eq:iaaaaaaaaz}\\
\hspace{0cm}&\quad+O_{\mathbb{P}}\left(\frac{1}{N+P+1}\right)
+O_{\mathbb{P}}\left(\frac{1}{N+Q+1}\right)\nonumber\\
\hspace{0cm}&=O_{\mathbb{P}}\left(\frac{1}{\sqrt{N+P+1}}\right)
+O_{\mathbb{P}}\left(\frac{1}{\sqrt{N+Q+1}}\right)
=O_{\mathbb{P}}\left(\frac{1}{\sqrt{N+1}}\right)\stackrel{\mathbb{P}}{\longrightarrow}0.\nonumber
\end{flalign}}Or, d'après~(\ref{eq:iaaaaaaaay}), on a, $\ \displaystyle \forall\ w_{1},\ w_{2},\ w_{3},\ w_{4}\ \in\mathbb{R},\ $ et lorsque $\
\displaystyle N\to\infty,\ $
{\small
\begin{flalign*}
\hspace{0cm}&\displaystyle\varphi_{\eta_{N}^{\prime},\ \eta_{N}^{\prime\prime},\ v_{N}^{\prime},\ v_{N}^{\prime\prime}}\left(w_1,\ w_2,\ w_3,\
w_4\right)=\mathbb{E}\left(e^{\displaystyle i\ w_1\ \eta_{N}^{\prime}+i\ w_2\ \eta_{N}^{\prime\prime}+i\ w_3\ v_{N}^{\prime}+i\ w_4\
v_{N}^{\prime\prime}}\right)\\
\hspace{0cm}&\longrightarrow\varphi(w_1,\ w_2,\ w_3,\ w_{4})=\exp\left[\displaystyle
-\frac{1}{2}\left\{m\left(w_{1}^{2}+w_{2}^{2}+w_{3}^{2}+w_{4}^{2}\right)\right\}\right].
\end{flalign*}}Donc, en posant $\ \displaystyle w_{1}=-\frac{1}{m}\frac{N+1}{N+P+1}\ w_{1},\ w_2=-\frac{1}{m}\frac{\sqrt{(N+1)P}}{N+P+1}\ w_{1},\ $\\*
\\*
$\displaystyle w_3=\frac{1}{m}\ \frac{N+1}{N+Q+1}\ w_{1}\ $ et $\ \displaystyle w_4=\frac{1}{m}\ \frac{\sqrt{(N+1)Q}}{N+Q+1}\ w_{1},\ $ on a, pour tout
$\ w_{1}\in\mathbb{R},\ $ et lorsque $\ N\to\infty,\ $
{\footnotesize
\begin{flalign*}
\hspace{0cm}&\displaystyle\varphi_{\displaystyle\left\{\frac{1}{m}\ \sqrt{N+1}\left(\frac{v_{N}^{\prime}\sqrt{N+1}+v_{N}^{\prime\prime}\sqrt{Q}}{N+Q+1}
-\frac{\eta_{N}^{\prime}\sqrt{N+1}
+\eta_{N}^{\prime\prime}\sqrt{P}}{N+P+1}\right)\right\}}(w_{1})\\
\hspace{0cm}&\longrightarrow\varphi\left(-\frac{1}{m}\frac{N+1}{N+P+1}\ w_{1},\ -\frac{1}{m}\frac{\sqrt{(N+1)P}}{N+P+1}\ w_1,\ \frac{1}{m}\frac{N+1}{N+Q+1}\
w_1,\ \ \frac{1}{m}\frac{\sqrt{(N+1)Q}}{N+Q+1}\ w_1\right)\\
\hspace{0cm}&=\exp\left(\displaystyle-\frac{1}{2}\left[m.\frac{1}{m^2}\
\lim_{N\to\infty}\left\{\frac{(N+1)^2+(N+1)P}{(N+P+1)^2}+\frac{(N+1)^2+(N+1)Q}{(N+Q+1)^2}\right\}\right]\ w_{1}^{2}\right)\\
\hspace{0cm}&=\exp\left(-\frac{1}{2}\
\sigma_{2}^{2}(c,d)\ w_{1}^{2}\right).
\end{flalign*}}Donc, on en déduit que, lorsque $\ \displaystyle N\to\infty,\ $
\begin{equation}
\sqrt{N+1}\ \mathcal{Q}_{N}\stackrel{d}{\longrightarrow}N(0,\sigma^{2}_{2}(c,d)),\label{eq:iaaaaaaaaaa}
\end{equation}où nous faisons usage de~(\ref{eq:iaaaaaaas}) pour écrire,
{\small
\begin{flalign*}
\hspace{0cm}\sigma^{2}_{2}(c,d)&=\frac{1}{m}\ \lim_{N\to\infty}\left\{\frac{(N+1)^2+(N+1)P}{(N+1+P)^2}+\frac{(N+1)^2+(N+1)Q}{(N+1+Q)^2}\right\}\\
\hspace{0cm}&=\frac{1}{m}\ \lim_{N\to\infty}\left\{\frac{\displaystyle 1+\frac{P}{N+1}}{\displaystyle 1+2\
\frac{P}{N+1}+\left(\frac{P}{N+1}\right)^2}+\frac{\displaystyle 1+\frac{Q}{N+1}}{\displaystyle 1+2\ \frac{Q}{N+1}
+\left(\frac{Q}{N+1}\right)^2}\right\}\\
\hspace{0cm}&=\frac{1}{m}\left\{\frac{1}{1+c}+\frac{1}{1+d}\right\}.
\end{flalign*}}Nous obtenons ainsi~(\ref{eq:iaaaaaaaam}) et~(\ref{eq:iaaaaaaaan}).\vskip5pt
\noindent Pour compléter notre preuve, nous combinons~(\ref{eq:iaaaaaaaaj}) avec~(\ref{eq:iaaaaaaaar}) et~(\ref{eq:iaaaaaaaaz}). Nous obtenons
que,

{\footnotesize
\begin{flalign}
\hspace{0cm}&(N+1)^{1/2}\mathcal{D}_{N}\label{eq:iaaaaaaaaab}\\
\hspace{0cm}&=(N+1)^{1/2}\mathcal{Q}_{N}+\frac{1}{2m}
\left\{\frac{N+1}{N+P+1}+\frac{N+1}{N+Q+1}\right\}(N+1)^{1/2}\Delta_{N}\nonumber\\
\hspace{0cm}&\quad+\frac{1}{2m}\left\{\frac{N+1}{N+P+1}-\frac{N+1}{N+Q+1}\right\}(N+1)^{1/2}\Theta_{N}\nonumber\\
\hspace{0cm}&=(N+1)^{1/2}\mathcal{Q}_{N}+\frac{1}{2m}\left\{\frac{N+1}{N+P+1}
+\frac{N+1}{N+Q+1}\right\}\left(\eta_{N}^{\prime}-v_{N}^{\prime}\right)\nonumber\\
\hspace{0cm}&\quad+\frac{1}{2m}\left\{\frac{N+1}{N+P+1}-\frac{N+1}{N+Q+1}\right\}\left(\eta_{N}^{\prime}+v_{N}^{\prime}\right)\nonumber\\
\hspace{0cm}&=\frac{1}{m}\ \frac{v_{N}^{\prime\prime}\sqrt{(N+1)Q}}{N+Q+1}-\frac{1}{m}\ \frac{\eta_{N}^{\prime\prime}\sqrt{(N+1)P}}{N+P+1}\nonumber\\
\hspace{0cm}&\quad+O_{\mathbb{P}}\left(\frac{1}{\sqrt{N+P+1}}\right)
+O_{\mathbb{P}}\left(\frac{1}{\sqrt{N+Q+1}}\right)\stackrel{d}{\longrightarrow}
N\left(0,\sigma_{1}^{2}(c,d)\right),\nonumber
\end{flalign}}où compte tenu de~(\ref{eq:iaaaaaaaay}),
{\footnotesize
\begin{flalign}
\hspace{0cm}\sigma_{1}^{2}(c,d)&=\lim_{N\to\infty}\left\{m.\frac{1}{m^2}\left(\frac{(N+1)Q}{(N+1+Q)^2}
+\frac{(N+1)P}{(N+1+P)^2}\right)\right\}\label{eq:iaaaaaaaaac}\\
\hspace{0cm}&=\frac{1}{m}\lim_{N\to\infty}\left\{\frac{(N+1)Q}{(N+1+Q)^2}+\frac{(N+1)P}{(N+1+P)^2}\right\}\nonumber\\
\hspace{0cm}&=\frac{1}{m}\lim_{N\to\infty}\left\{\frac{\displaystyle\frac{Q}{N+1}}{\displaystyle1+2\
\frac{Q}{N+1}+\left(\frac{Q}{N+1}\right)^2}+\frac{\displaystyle \frac{P}{N+1}}{\displaystyle1+2\
\frac{P}{N+1}+\left(\frac{P}{N+1}\right)^2}\right\}\nonumber\\
\hspace{0cm}&=\frac{1}{m}\left\{\frac{d}{(1+d)^2}+\frac{c}{(1+c)^2}\right\},\nonumber
\end{flalign}}ce qui donne~(\ref{eq:iaaaaaaaal}).$\Box$\vskip5pt
\noindent En faisant usage du fait~\ref{fai:fa202}, nous pouvons définir, pour chaque $\ \displaystyle N\geq 1,\ $ une v.a. $\
\displaystyle\theta_{N}^{\prime}\stackrel{d}{=}N(0,1),\ $ indépendante de $\ \displaystyle \left\{\zeta_{\ell,N}:\ell\geq 1\right\}\ $ telle que, lorsque $\ \displaystyle N\to\infty,\ $
\begin{equation}
\hspace{-1cm}\frac{1}{m}\frac{\eta_{N}^{\prime\prime}\sqrt{(N+1)P}}
{N+P+1}-\frac{1}{\sqrt{m}}\frac{\sqrt{(N+1)P}}{N+P+1}\theta_{N}^{\prime}
=O_{\mathbb{P}}\left(\frac{\sqrt{N+1}\log P}{N+P+1}\right).\label{eq:iaaaaaaaaaf1}
\end{equation}De même, en faisant usage du fait~\ref{fai:fa202}, nous pouvons définir, pour chaque $\ \displaystyle N\geq 1,\ $ une v.a. $\
\displaystyle\theta_{N}^{\prime\prime}\stackrel{d}{=}N(0,1),\ $ indépendante de $\ \displaystyle \left\{\xi_{\ell,N}:\ell\geq 1\right\}\ $ telle que, lorsque $\ \displaystyle N\to\infty,\ $
\begin{equation}
\hspace{-1cm}\frac{1}{m}\frac{v_{N}^{\prime\prime}\sqrt{(N+1)Q}}
{N+Q+1}-\frac{1}{\sqrt{m}}\frac{\sqrt{(N+1)Q}}{N+Q+1}\theta_{N}^{\prime\prime}
=O_{\mathbb{P}}\left(\frac{\sqrt{N+1}\log Q}{N+Q+1}\right).\label{eq:iaaaaaaaaaf2}
\end{equation}Donc d'après~(\ref{eq:iaaaaaaaaaf1}) et~(\ref{eq:iaaaaaaaaaf2}), nous pouvons définir pour chaque $\ \displaystyle N\geq 1,\ $ une v.a. $\
\displaystyle\theta_{N}\stackrel{d}{=}N(0,1),\ $ indépendante de $\ \displaystyle \left\{\zeta_{\ell,N}:\ell\geq 1\right\}\ $ et $\
\displaystyle \left\{\xi_{\ell,N}:\ell\geq 1\right\}\ $ (et donc, de $\ \displaystyle\Delta_{N},\ \Theta_{N},\ \phi_{N}\ \mbox{ et }\ \psi_{N},\
$ comme défini dans la proposition~\ref{propos:propo256}), telle que, lorsque $\ \displaystyle N\to\infty,\ $

{\small
\begin{flalign}
&\hspace{0cm}\frac{1}{m}\frac{v_{N}^{\prime\prime}\sqrt{(N+1)Q}}{N+Q+1}-\frac{1}{m}\frac{\eta_{N}^{\prime\prime}\sqrt{(N+1)P}}
{N+P+1}-\left\{\frac{1}{m}\left(\frac{(N+1)P}{(N+P+1)^2}\right.\right.\label{eq:iaaaaaaaaaf}\\
&\hspace{0cm}\hspace{7cm}\left.\left.+\frac{(N+1)Q}{(N+Q+1)^2}\right)\right\}^{1/2}\theta_{N}
\nonumber\\
&\hspace{0cm}=O_{\mathbb{P}}\left(\frac{\sqrt{N+1}\log P}{N+P+1}\right)+O_{\mathbb{P}}\left(\frac{\sqrt{N+1}\log
Q}{N+Q+1}\right)\nonumber\\
&\hspace{0cm}=O_{\mathbb{P}}\left(\frac{\sqrt{N_1-P+1}\log (N_1-N)}{N_1+1}\right)+O_{\mathbb{P}}\left(\frac{\sqrt{N_2-Q+1}\log(N_2-N)}{N_2+1}\right)\nonumber\\
&\hspace{0cm}=O_{\mathbb{P}}\left(\frac{\log(N_1+1)}{\sqrt{N_1+1}}\right)+O_{\mathbb{P}}\left(\frac{\log
(N_2+1)}{\sqrt{N_2+1}}\right)\nonumber\\
&\hspace{0cm}=O_{\mathbb{P}}\left(\frac{\log (N+P+1)}{\sqrt{N+P+1}}\right)+O_{\mathbb{P}}\left(\frac{\log
(N+Q+1)}{\sqrt{N+Q+1}}\right)=O_{\mathbb{P}}\left(\frac{\log N}{\sqrt{N}}\right).\nonumber
\end{flalign}}{\bf Démonstration du lemme~\ref{lem:le257}.} Faisons un développement limité de Taylor à l'ordre 2, de l'expression $\ \displaystyle G_{3}^{(m)}(\tau_{N}(t))-G_{3}^{(m)}(t)\ $ dans~(\ref{eq:iaaaaaaad}). On a, uniformément en $\ \displaystyle t\in[0,1],\ $ et lorsque $\ N\to\infty,\ $
\begin{flalign}
&G_{3}^{(m)}(\tau_{N}(t))-G_{3}^{(m)}(t)\label{eq:iaaaaaaae}\\
&=\left(\tau_{N}(t)-t\right)G_{3}^{\prime(m)}(t)+\frac{\left(\tau_{N}(t)-t\right)^2}{2}\
G_{3}^{\prime\prime(m)}(t)+O\left[\left(\tau_{N}(t)-t\right)^2\right]\nonumber\\
&=\frac{(2m-1)!}{\left((m-1)!\right)^2}\ \left(t(1-t)\right)^{m-1}\left(\tau_{N}(t)-t\right)\nonumber\\
&\quad+\frac{(2m-1)!}{2(m-2)!(m-1)!}\ (1-2t)\ (t(1-t))^{m-2}(\tau_{N}(t)-t)^2\nonumber\\
&\quad+O\left[\left(\tau_{N}(t)-t\right)^2\right].\nonumber
\end{flalign}Or, en se rappelant de~(\ref{eq:iaaaaaaai}) et~(\ref{eq:iaaaaaaae}), on a, lorsque $\ \displaystyle N\to\infty,\ $
\begin{displaymath}
\sup_{0\leq t\leq 1}\left\vert\tau_{N}(t)-t\right\vert\leq
\frac{1}{4}\left\vert\mathcal{Q}_{N}\right\vert\left\{1-\left\vert\mathcal{Q}_{N}\right\vert\right\}^{-1}
=O_{\mathbb{P}}\left(\frac{1}{\sqrt{N+1}}\right).
\end{displaymath}D'où, on a, uniformément en $\ \displaystyle t\in[0,1],\ $ et lorsque $\ N\to\infty,\ $
\begin{equation}
\tau_{N}(t)-t=O_{\mathbb{P}}\left(\frac{1}{\sqrt{N+1}}\right).\label{eq:iaaaaaaaf}
\end{equation}Donc, de~(\ref{eq:iaaaaaaae}),~(\ref{eq:iaaaaaaaf}) et comme la fonction $\ \displaystyle t\mapsto (1-2t)\ (t(1-t))^{m-2}\ $ est
bornée sur $\ \displaystyle [0,1],\ $ on en déduit que, l'on a, uniformément en $\ \displaystyle t\in[0,1],\ $ et lorsque $\ \displaystyle N\to\infty,\ $
{\footnotesize
\begin{equation}
G_{3}^{(m)}(\tau_{N}(t))-G_{3}^{(m)}(t)=\frac{(2m-1)!}{\left((m-1)!\right)^2}\
\left(t(1-t)\right)^{m-1}\left(\tau_{N}(t)-t\right)+O_{\mathbb{P}}\left(\frac{1}{N+1}\right)\label{eq:iaaaaaaag}.
\end{equation}}D'où, de~(\ref{eq:iaaaaaaad}) et~(\ref{eq:iaaaaaaag}), on en déduit~(\ref{eq:iaaaaaaah}), ce qui achève la démonstration du
lemme~\ref{lem:le257}.$\Box$\vskip5pt
\begin{lemme}
\label{lem:le258}
\noindent Nous avons, lorsque $\ \displaystyle N\to\infty,\ $
{\small
\begin{equation}
\sup_{\displaystyle0\leq t\leq 1}\left\vert \left(N+1\right)^{1/2}\left(\tau_{N}(t)-t\right)+t(1-t)\
(N+1)^{1/2}\mathcal{Q}_{N}\right\vert=O_{\mathbb{P}}\left(\frac{1}{\sqrt{N}}\right).\label{eq:iaaaaaaaj}
\end{equation}}
\end{lemme}
\noindent{\bf Démonstration.} En se rappelant de l'expression~(\ref{eq:iaaaaaaai}) de $\ \displaystyle \tau_{N}(t)-t\ $ et en faisant un D.L. à
l'ordre 2 de $\ \displaystyle \tau_{N}(t)-t,\ $ on a, uniformément en $\ t\in[0,1],\ $ et lorsque $\ \displaystyle N\to\infty,\ $
\begin{flalign*}
&\left(N+1\right)^{1/2}\left(\tau_{N}(t)-t\right)\\
&=-\left(N+1\right)^{1/2}t(1-t)\mathcal{Q}_{N}.\frac{1}{1+(1-t)\mathcal{Q}_{N}}\\
&=-t(1-t)(N+1)^{1/2}\mathcal{Q}_{N}\left(1-(1-t)\mathcal{Q}_{N}\right)\\
&=-t(1-t)(N+1)^{1/2}\mathcal{Q}_{N}+t(1-t)^2\ (N+1)^{1/2}\mathcal{Q}_{N}^{2}.
\end{flalign*}Or, d'après~(\ref{eq:iaaaaaaaam}), on a, lorsque $\ \displaystyle N\to\infty,\ $ $\ \displaystyle
(N+1)^{1/2}\mathcal{Q}_{N}=O_{\mathbb{P}}(1),\ $ de plus\\*
$\displaystyle t\mapsto t(1-t)^2\ $ est bornée sur $\ [0,1],\ $ donc, on a, lorsque $\ \displaystyle N\to\infty,\ $
\begin{displaymath}
\left(N+1\right)^{1/2}\left(\tau_{N}(t)-t\right)=-t(1-t)(N+1)^{1/2}\mathcal{Q}_{N}+O_{\mathbb{P}}\left((N+1)^{1/2}\mathcal{Q}_{N}^{2}\right).
\end{displaymath}Et donc, on en déduit que, lorsque $\ \displaystyle N\to\infty,\ $
\begin{flalign*}
&\sup_{0\leq t\leq 1}\left\vert \left(N+1\right)^{1/2}\left(\tau_{N}(t)-t\right)+t(1-t)(N+1)^{1/2}\mathcal{Q}_{N}\right\vert\\
&=O_{\mathbb{P}}\left((N+1)^{1/2}\mathcal{Q}_{N}^{2}\right)=O_{\mathbb{P}}\left(\frac{1}{\sqrt{N}}\right),
\end{flalign*}ce qui achève la démonstration du lemme~\ref{lem:le258}.$\Box$\vskip5pt
\noindent Le fait suivant est une conséquence directe du théorème 1.2 de Deheuvels et Einmahl\cite{MR1797314} (voir aussi Stute\cite{MR637378} et le théorème 3.1 de Deheuvels et Mason\cite{MR1175262}, pour les versions antérieures et les variantes de ce résultat).\vskip5pt
\begin{fait}\label{fai:fa204} Soit $\ \displaystyle\left\{h_{N}:N\geq 1\right\}\ $ une suite de constantes positives telle que, lorsque $\
\displaystyle N\to\infty,\ $
\begin{displaymath}
h_{N}\to\infty,\ \frac{N\ h_{N}}{\log N}\to\infty\mbox{ et }\ \frac{\log(1/h_{N})}{\log\log N}\to\infty.\
\end{displaymath}Alors, nous avons, lorsque $\ \displaystyle N\to\infty,\ $
\begin{equation}
\left\{2\ h_{N}\ \log(1/h_{N})\right\}^{-1/2}\sup_{\displaystyle\substack{\displaystyle0\leq t,s\leq 1 \\ \displaystyle\vert t-s\vert\leq
h_{N}}}\left\vert\alpha_{N;1}(t)-\alpha_{N;1}(s)\right\vert\stackrel{\mathbb{P}}{\longrightarrow}1.\label{eq:iaaaaaaap}
\end{equation}
\end{fait}
\begin{lemme}
\label{lem:le259}
\noindent Nous avons, lorsque $\ \displaystyle N\to\infty,\ $
{\footnotesize
\begin{equation}
\hspace{-0,5cm}\sup_{\displaystyle0\leq t\leq
1}\left\vert\alpha_{N;1}\left(G_{3}^{(m)}\left(\tau_{N}(t)\right)\right)-\alpha_{N;1}\left(G_{3}^{(m)}(t)\right)\right\vert
=O_{\mathbb{P}}\left(N^{-1/4}(\log(N))^{1/2}\right).\label{eq:iaaaaaaaq}
\end{equation}}
\end{lemme}
\noindent{\bf Démonstration.} D'après~(\ref{eq:iaaaaaaaf}) et~(\ref{eq:iaaaaaaag}), et comme $\ \displaystyle
t\mapsto\frac{(2m-1)!}{((m-1)!)^2}\left(t(1-t)\right)^{m-1}\ $ est bornée sur $\ \displaystyle [0,1],\ $ on a, uniformément en
$\ \displaystyle t\in[0,1],\ $ et lorsque $\ \displaystyle N\to\infty,\ $
\begin{displaymath}
G_{3}^{(m)}(\tau_{N}(t))-G_{3}^{(m)}(t)=O_{\mathbb{P}}\left(\displaystyle\frac{1}{\sqrt{N}}\right).
\end{displaymath}D'où, on en déduit, que lorsque $\ \displaystyle N\to\infty,\ $ on a,
\begin{displaymath}
\sup_{\displaystyle0\leq t\leq 1}\left\vert
G_{3}^{(m)}(\tau_{N}(t))-G_{3}^{(m)}(t)\right\vert=O_{\mathbb{P}}\left(\displaystyle\frac{1}{\sqrt{N}}\right).
\end{displaymath}Par conséquent, $\ \displaystyle \forall\ \varepsilon>0,\ \exists\ C_{\varepsilon}<\infty,\ $ tel que
\begin{displaymath}
\limsup_{\displaystyle N\to\infty}\mathbb{P}\left[\sup_{\displaystyle0\leq t\leq 1}\left\vert
G_{3}^{(m)}(\tau_{N}(t))-G_{3}^{(m)}(t)\right\vert\geq \frac{C_{\varepsilon}}{\sqrt{N}}\right]<\varepsilon.
\end{displaymath}Donc, avec une probabilité supérieure ou égale à $\ \displaystyle1-\varepsilon,\ $ on a, à partir d'un certain rang,
\begin{displaymath}
\sup_{\displaystyle0\leq t\leq 1}\left\vert G_{3}^{(m)}(\tau_{N}(t))-G_{3}^{(m)}(t)\right\vert\leq\frac{C_{\varepsilon}}{\sqrt{N}}.
\end{displaymath}D'où, pour tout $\ \displaystyle t\in[0,1],\ $ on a,
\begin{displaymath}
\left\vert G_{3}^{(m)}\left(\tau_{N}(t)\right)-G_{3}^{(m)}(t)\right\vert\leq\frac{C_{\varepsilon}}{\sqrt{N}}.
\end{displaymath}Posons, pour tout $\ \displaystyle N\geq 1,\ h_{N}=\frac{C_{\varepsilon}}{\sqrt{N}}.\ $ On constate que, lorsque $\ \displaystyle
N\to\infty,\ $
\begin{displaymath}
h_{N}\downarrow0,\ \frac{N\ h_{N}}{\log N}\to\infty\ \mbox{ et }\frac{\log(1/h_{N})}{\log\log N}\to\infty.
\end{displaymath}Donc, $\ \displaystyle h_{N}\ $ vérifie les conditions du fait~\ref{fai:fa204} appliqué au processus empirique uniforme $\ \displaystyle \alpha_{N;1}.\ $ On a, donc, lorsque $\ \displaystyle N\to\infty,\ $
\begin{flalign*}
&\left\{2\ C_{\varepsilon}\ N^{-1/2}\ \log\left(N^{1/2}C_{\varepsilon}^{-1}\right)\right\}^{-1/2}\times\\
&\sup_{\substack{\displaystyle0\leq t\leq 1 \\ \displaystyle G_{3}^{(m)}(\tau_{N}(t))-G_{3}^{(m)}(t)\leq
\frac{C_{\varepsilon}}{\sqrt{N}}}}\left\vert\alpha_{N;1}\left(G_{3}^{(m)}(\tau_{N}(t))\right)-\alpha_{N;1}\left(G_{3}^{(m)}(t)
\right)\right\vert=O_{\mathbb{P}}(1),
\end{flalign*}i.e,
\begin{flalign*}
&\left\{\sqrt{C_{\varepsilon}}\ N^{-1/4}\ \left(\log N\right)^{1/2}\right\}^{-1}\times\nonumber\\
&\sup_{\substack{\displaystyle0\leq t\leq 1 \\ \displaystyle G_{3}^{(m)}(\tau_{N}(t))-G_{3}^{(m)}(t)\leq
\frac{C_{\varepsilon}}{\sqrt{N}}}}\left\vert\alpha_{N;1}\left(G_{3}^{(m)}(\tau_{N}(t))\right)-\alpha_{N;1}
\left(G_{3}^{(m)}(t)\right)\right\vert=O_{\mathbb{P}}(1).
\end{flalign*}D'où, on en déduit~(\ref{eq:iaaaaaaaq}), ce qui achève la démonstration du lemme~\ref{lem:le259}.$\Box$\vskip5pt
\begin{prop}\label{propos:propo258} Sur un espace $(\Omega,\mathcal{A},\mathbb{P})$ de probabilité
convenable, il existe une suite $\{B_N(.):N\geq 1\}$ de ponts Browniens, et des suites $\{\phi_N,\psi_N,\theta_N:N\geq 1\}$ de variables aléatoires normales $N(0,1)$, tels que, pour chaque $N\geq 1$,  $\left\{B_{N}(v):0\leq v\leq 1\right\}$, $\phi_{N}$, $\psi_{N}$ et $\theta_{N}$ soient indépendants, et vérifient, lorsque $N\rightarrow\infty$,
{\footnotesize
\begin{flalign}
&\sup_{0\leq t\leq 1}\left\vert \gamma_{N;n_1,n_2}(t)-B_{N}\left(H_{m}(t)\right)+\frac{(2m-1)!}{\left((m-1)!\right)^2}\
\left(t(1-t)\right)^m\times\right.\label{eq:iaaaaaaaaae}\\
&\hspace{0,5cm}\bigg[\left\{\frac{1}{m}\left(\frac{(N+1)P}{(N+1+P)^2}+\frac{(N+1)Q}{(N+1+Q)^2}\right)\right\}^{1/2}\ \theta_{N}\nonumber\\
&\hspace{0,5cm}-\frac{1}{2m}\left\{\frac{N+1}{N+1+P}+\frac{N+1}{N+1+Q}\right\}\left(\phi_{N}\ \sqrt{\frac{2m}{2m+1}}-4m\int_{0}^{1}B_{N}(H_{m}(x))\ \ud
x\right)\nonumber\\
&\hspace{0,5cm}\left.\left.-\frac{1}{2m}\left\{\frac{N+1}{N+1+P}-\frac{N+1}{N+1+Q}\right\}\psi_{N}\
\sqrt{2m}\right]\right\vert=O_{\mathbb{P}}\left(N^{-1/4}\ \left(\log N\right)^{1/2}\right)
.\nonumber
\end{flalign}}
\end{prop}
\noindent{\bf Démonstration.} Nous déduisons facilement~(\ref{eq:iaaaaaaaaae}) de~(\ref{eq:iaaaaaaah}) dans le lemme~\ref{lem:le257}, de~(\ref{eq:iaaaaaaaq}) dans
le lemme~\ref{lem:le259}, de~(\ref{eq:iaaaaaaaj}) dans le lemme~\ref{lem:le258}, des 1\iere{} et 4\ieme{} parties de~(\ref{eq:iaaaaaaaaab}) dans
la démonstration du lemme~\ref{lem:le260} et de~(\ref{eq:iaaaaaaaaaf}), en combinaison avec~(\ref{eq:iaaaaax}),~(\ref{eq:iaaaaay})
et~(\ref{eq:iaaaaaz}) dans le théorème~\ref{thm:teo11111}.$\Box$\vskip5pt

\subsection{Preuve du théorème 1.2}\label{ma-soussectionm}
\noindent D'après~(\ref{eq:iaaaaaaaaae}), nous avons, lorsque $\ \displaystyle N\to\infty,$
{\footnotesize
\begin{flalign*}
&\sup_{0\leq t\leq 1}\left\vert \gamma_{N;n_1,n_2}(t)-B_{N}\left(H_{m}(t)\right)+\frac{(2m-1)!}{\left((m-1)!\right)^2}\
\left(t(1-t)\right)^m\times\right.\\
&\hspace{0,5cm}\bigg[\left\{\frac{1}{m}\left(\frac{(N+1)P}{(N+1+P)^2}+\frac{(N+1)Q}{(N+1+Q)^2}\right)\right\}^{1/2}\ \theta_{N}\nonumber\\
&\hspace{0,5cm}-\frac{1}{2m}\left\{\frac{N+1}{N+1+P}+\frac{N+1}{N+1+Q}\right\}\left(\phi_{N}\ \sqrt{\frac{2m}{2m+1}}-4m\int_{0}^{1}B_{N}(H_{m}(x))\ \ud
x\right)\nonumber\\
&\hspace{0,5cm}\left.\left.-\frac{1}{2m}\left\{\frac{N+1}{N+1+P}-\frac{N+1}{N+1+Q}\right\}\psi_{N}\
\sqrt{2m}\right]\right\vert=O_{\mathbb{P}}\left(N^{-1/4}\ \left(\log N\right)^{1/2}\right)
.\nonumber
\end{flalign*}}Maintenant, nous avons les identités en loi suivantes. En se rappelant de la définition~(\ref{eq:iaaaaaaaw}) de $\ \displaystyle H_{m}\ $ et en posant $\ \displaystyle \left\{B(t):0\leq t\leq 1\right\}\ $ qui désigne un pont Brownien, $\ \displaystyle \phi,\ $ $\ \displaystyle \psi,\ $ et $\ \displaystyle \theta,\ $ qui désignent des variables aléatoires normales $N(0,1)$ tels que,  $\ \displaystyle\left\{B(t):0\leq t\leq 1\right\},\ $ $\ \displaystyle\phi,\ $ $\ \displaystyle\psi,\ $ et $\ \displaystyle\theta\ $ soient indépendants, nous avons,

{\footnotesize
\begin{flalign*}
{\mathfrak{Y}}_{N}(t)&:=B_{N}\left(H_{m}(t)\right)-\frac{(2m-1)!}{\left((m-1)!\right)^2}\
\left(t(1-t)\right)^m\times\\
&\hspace{0,5cm}\bigg[\left\{\frac{1}{m}\left(\frac{(N+1)P}{(N+1+P)^2}+\frac{(N+1)Q}{(N+1+Q)^2}\right)\right\}^{1/2}\ \theta_{N}\nonumber\\
&\hspace{0,5cm}-\frac{1}{2m}\left\{\frac{N+1}{N+1+P}+\frac{N+1}{N+1+Q}\right\}\\
&\hspace{0,5cm}\times\left(\phi_{N}\ \sqrt{\frac{2m}{2m+1}}-4m\int_{0}^{1}B_{N}(H_{m}(x))\ \ud
x\right)\nonumber\\
&\hspace{0,5cm}\left.-\frac{1}{2m}\left\{\frac{N+1}{N+1+P}-\frac{N+1}{N+1+Q}\right\}\psi_{N}\
\sqrt{2m}\right]\\
&\stackrel{d}{=}B\left(H_{m}(t)\right)-4(2m+1)\frac{(2m-1)!}{2m\left((m-1)!\right)^2}\
\left(t(1-t)\right)^m\\
&\quad\times\left\{\frac{m}{2m+1}\ \left[\frac{N+1}{N+1+P}+\frac{N+1}{N+1+Q}\right]\right\}\int_{0}^{1}B(H_{m}(x))\ \ud
x\\
&\quad+\frac{(2m-1)!}{2m\ \left((m-1)!\right)^2}\
\left(t(1-t)\right)^m\left[\left\{\frac{N+1}{N+1+P}+\frac{N+1}{N+1+Q}\right\}\right.\\
&\quad\times\phi\ \sqrt{\frac{2m}{2m+1}}+\left\{\frac{N+1}{N+1+P}-\frac{N+1}{N+1+Q}\right\}\psi\
\sqrt{2m}\\
&\quad-\left\{\frac{(N+1)P}{(N+1+P)^2}+\frac{(N+1)Q}{(N+1+Q)^2}\right\}^{1/2}\ \theta\ 2\sqrt{m}\bigg]\\
&=(B\circ H_{m})_{\left\{\displaystyle\frac{m}{2m+1}\ \left[\frac{N+1}{N+1+P}+\frac{N+1}{N+1+Q}\right]\right\}}(t)+\Psi(t)Y,
\end{flalign*}}où $\ \displaystyle (B\circ H_{m})_{C}(\cdot)\ $ est comme dans~(\ref{eq:mr}), $\ \displaystyle\Psi(t)=\frac{(2m-1)!}{2m\left((m-1)!\right)^2}\
\left(t(1-t)\right)^m\ $ et $\ \displaystyle\sigma^{2}_{\mathfrak{X}}=\sigma^{2}_{B\circ H_{m}}=\frac{1}{4(2m+1)}\ $ sont comme dans~(\ref{eq:mo}) et~(\ref{eq:mn2}), et
{\small
\begin{flalign*}
Y&\stackrel{d}{=}N\left(0,\frac{2m}{2m+1}\left\{\frac{N+1}{N+1+P}+\frac{N+1}{N+1+Q}\right\}^2
+2m\left\{\frac{N+1}{N+1+P}-\frac{N+1}{N+1+Q}\right\}^2\right.\\
&\left.\hspace{2cm}+4m\left\{\frac{(N+1)P}{(N+1+P)^2}+\frac{(N+1)Q}{(N+1+Q)^2}\right\}\right)
\end{flalign*}}est indépendant de $\ \displaystyle (B\circ H_{m})(.).\ $ Maintenant, nous appliquons~(\ref{eq:ml}), pris avec
{\footnotesize
\begin{flalign*}
&C=\frac{m}{2m+1}\left\{\frac{N+1}{N+1+P}+\frac{N+1}{N+1+Q}\right\},\ \sigma^2=\frac{2m}{2m+1}\left\{\frac{N+1}{N+1+P}+\frac{N+1}{N+1+Q}\right\}^2\\
&+2m\left\{\frac{N+1}{N+1+P}-\frac{N+1}{N+1+Q}\right\}^2+4m\left\{\frac{(N+1)P}{(N+1+P)^2}+\frac{(N+1)Q}{(N+1+Q)^2}\right\},\ \mathfrak{X}=B\circ H_{m}\\
&\mbox{et } \sigma^{2}_{\mathfrak{X}}=\frac{1}{4(2m+1)}.
\end{flalign*}}Il est facile de montrer que
{\small
\begin{flalign*}
(1-C)^2+\sigma^2\sigma_{\mathfrak{X}}^2&=\left(1-\frac{m}{2m+1}\left\{\frac{N+1}{N+1+P}+\frac{N+1}{N+1+Q}\right\}\right)^2\\
&\quad+\frac{m}{2(2m+1)^2}\left\{\frac{N+1}{N+1+P}+\frac{N+1}{N+1+Q}\right\}^2\\
&\quad+\frac{m}{2(2m+1)}\left\{\frac{N+1}{N+1+P}-\frac{N+1}{N+1+Q}\right\}^2\\
&\quad+\frac{m}{2m+1}\left\{\frac{(N+1)P}{(N+1+P)^2}+\frac{(N+1)Q}{(N+1+Q)^2}\right\}\\
&=1-\frac{m}{2m+1}\left\{\frac{N+1}{N+1+P}+\frac{N+1}{N+1+Q}\right\}:=R_{N,m}.
\end{flalign*}}En se rappelant de~(\ref{eq:me}) et~(\ref{eq:mh}), nous obtenons aussi que
\begin{flalign*}
\hspace{-1cm}{\mathfrak{Y}}_{N}(\cdot)&=(B\circ H_{m})_{\left\{\displaystyle\frac{m}{2m+1}\ \left[\frac{N+1}{N+1+P}+\frac{N+1}{N+1+Q}\right]\right\}}(\cdot)+\Psi(\cdot)Y\\
\hspace{-1cm}&\stackrel{d}{=}(B\circ H_{m})_{\left\{1\pm\sqrt{R_{N,m}}\right\}}(\cdot)=J_{\left\{1\pm\sqrt{R_{N,m}}\right\}}B\circ H_{m}(\cdot).
\end{flalign*}Ceci montre que les suites de processus $\ \displaystyle \left\{{\mathcal{B}}_{N}^{+}(t):0\leq t\leq 1\right\}\ $ et\\*
$\displaystyle \left\{{\mathcal{B}}_{N}^{-}(t):0\leq t\leq 1\right\},\ $ $\ \displaystyle N=1,2,\ldots,\ $ définies en posant
\begin{flalign}
{\mathcal{B}}_{N}^{+}\circ H_{m}(\cdot)&=J_{\left\{\frac{R_{N,m}+\sqrt{R_{N,m}}}{R_{N,m}}\right\}}\mathfrak{Y}_{N}(\cdot)
=J_{\left\{1+\sqrt{R_{N,m}}\right\}}^{-1}\mathfrak{Y}_{N}(\cdot)\label{eq:maf}\\
&\stackrel{d}{=}J_{\left\{1+\sqrt{R_{N,m}}\right\}}^{-1}\circ J_{\left\{1+\sqrt{R_{N,m}}\right\}}B\circ H_{m}(\cdot)=B\circ H_{m}(\cdot),\nonumber
\end{flalign}
\begin{flalign}
{\mathcal{B}}_{N}^{-}\circ H_{m}(\cdot)&=J_{\left\{\frac{R_{N,m}-\sqrt{R_{N,m}}}{R_{N,m}}\right\}}\mathfrak{Y}_{N}(\cdot)
=J_{\left\{1-\sqrt{R_{N,m}}\right\}}^{-1}\mathfrak{Y}_{N}(\cdot)\label{eq:mag}\\
&\stackrel{d}{=}J_{\left\{1-\sqrt{R_{N,m}}\right\}}^{-1}\circ J_{\left\{1-\sqrt{R_{N,m}}\right\}}B\circ H_{m}(\cdot)=B\circ H_{m}(\cdot),\nonumber
\end{flalign}sont toutes les deux des ponts Browniens qui vérifient les identités
\begin{equation}
{\mathfrak{Y}}_{N}=J_{\left\{1+\sqrt{R_{N,m}}\right\}}
{\mathcal{B}}_{N}^{+}\circ H_{m}=J_{\left\{1-\sqrt{R_{N,m}}\right\}}{\mathcal{B}}_{N}^{-}\circ H_{m},\label{eq:mah}
\end{equation}ensemble avec
\begin{equation}
\hspace{-1cm}{\mathcal{B}}_{N}^{+}\circ H_{m}=J_{\left\{\frac{R_{N,m}+\sqrt{R_{N,m}}}{R_{N,m}}\right\}}\circ J_{\left\{1-\sqrt{R_{N,m}}\right\}}{\mathcal{B}}_{N}^{-}\circ H_{m}=J_{2}\ {\mathcal{B}}_{N}^{-}\circ H_{m},\label{eq:mai}
\end{equation}et
\begin{equation}
\hspace{-1cm}{\mathcal{B}}_{N}^{-}\circ H_{m}=J_{\left\{\frac{R_{N,m}-\sqrt{R_{N,m}}}{R_{N,m}}\right\}}\circ J_{\left\{1+\sqrt{R_{N,m}}\right\}}{\mathcal{B}}_{N}^{+}\circ H_{m}=J_{2}\ {\mathcal{B}}_{N}^{+}\circ H_{m}.\label{eq:maj}
\end{equation}Ceci suffit pour~(\ref{eq:my}),~(\ref{eq:mz}) et~(\ref{eq:maa}) et complète la preuve du théorème~\ref{thm:teo11113}.$\Box$\vskip5pt

\subsection{Preuve du théorème 1.1}\label{ma-soussectionn}
\noindent La preuve est essentiellement similaire à celle du théorème~\ref{thm:teo11113}. Soit~(\ref{eq:iaaaaaaas}) vérifié avec $\ \displaystyle c=d=0.\ $ Une application de~(\ref{eq:iaaaaaaaaae}) montre que, lorsque $\ \displaystyle N\to\infty,\ $
\begin{flalign}
&\hspace{-0,5cm}\sup_{0\leq t\leq 1}\left\vert \gamma_{N;n_1,n_2}(t)-B_{N}\left(H_{m}(t)\right)+\frac{(2m-1)!}{\left((m-1)!\right)^2}\
\left(t(1-t)\right)^m\times\right.\label{eq:iaaaaaaaab1}\\
&\left.\left\{4\int_{0}^{1}B_{N}(H_{m}(x))\ \ud x-\frac{1}{m}\ \phi_{N}\ \sqrt{\frac{2m}{2m+1}}\right\}\right\vert=o_{\mathbb{P}}(1).\nonumber
\end{flalign}Par commodité, nous réécrivons le processus gaussien approchant $\ \displaystyle \gamma_{N;n_1,n_2}(.)\ $ dans~(\ref{eq:iaaaaaaaab1}) en
\begin{flalign*}
&B_{N}\left(H_{m}(t)\right)-4(2m+1)\frac{(2m-1)!}{2m\left((m-1)!\right)^2}\
\left(t(1-t)\right)^m\left\{\frac{2m}{2m+1}\right\}\\
&\times\int_{0}^{1}B_{N}(H_{m}(x))\ \ud x+\frac{(2m-1)!}{2m\left((m-1)!\right)^2}\
\left(t(1-t)\right)^m\left\{2\phi_{N}\ \sqrt{\frac{2m}{2m+1}}\right\}.
\end{flalign*}Nous appliquons donc~(\ref{eq:ml}) pris avec
\begin{displaymath}
C=\frac{2m}{2m+1},\ \sigma^2=\frac{8m}{2m+1},\ \mathfrak{X}=B\circ H_{m}\ \mbox{ et}\ \sigma_{\mathfrak{X}}^{2}=\frac{1}{4(2m+1)}.
\end{displaymath}Des calculs de routine montrent que
\begin{displaymath}
(1-C)^2+\sigma^2\sigma_{\mathfrak{X}}^2=1-\frac{2m}{2m+1}=\frac{1}{2m+1}.
\end{displaymath}\'Etant donné ce fait, nous faisons usage de~(\ref{eq:iaaaaaaaab1}) et complétons la preuve du théorème~\ref{thm:teo11112} en suivant mot à mot les arguments utilisés pour la preuve juste donnée du théorème~\ref{thm:teo11113}.$\Box$\vskip5pt

\subsection{Preuve du théorème 1.3}\label{ma-soussectiono}
\noindent Soient, pour $\ \displaystyle k=0,\ldots,N,\ $ $\ \displaystyle S_{k,N;X}^{(m)}$ et $\displaystyle S_{k,N;Y}^{(m)}$ qui désignent les paires de $m$-espacements uniformes disjoints sur $\ \displaystyle (0,1),\ $  engendrées par des suites indépendantes de v.a. i.i.d. de loi uniforme sur $\ \displaystyle(0,1) $ définies comme dans~(\ref{eq:iaaaaaaar}), et posons $\ \displaystyle P=P(N)\ $ et $\ \displaystyle Q=Q(N)\ $ comme dans~(\ref{eq:iaaaaaaaag}). Il est facile de montrer que, sous les hypothèses du théorème~\ref{thm:teo11114}, nous avons, pour tout $\ \displaystyle 0\leq k\leq N,\ $ et pour tout $\ \displaystyle 0\leq t\leq 1,\ $ l'égalité d'événements suivante,
{\small
\begin{flalign}
&\left\{\frac{e(N+1+P)\ S_{k,N;X}^{(m)}}{e(N+1+P)\ S_{k,N;X}^{(m)}+h(N+1+Q)\ S_{k,N;Y}^{(m)}}\leq t\right\}\label{eq:iaaaaaaaaae2}\\
&=\left\{\frac{(N+1+P)\ S_{k,N;X}^{(m)}}{(N+1+P)\ S_{k,N;X}^{(m)}+(N+1+Q)\ S_{k,N;Y}^{(m)}}\leq \frac{th}{th+(1-t)e}\right\}.\nonumber
\end{flalign}}Compte tenu de~(\ref{eq:iaaaaaaaaae2}) , nous voyons que~(\ref{eq:mab}),~(\ref{eq:mac}) et~(\ref{eq:mad}) sont des conséquences directes de~(\ref{eq:iaaaaaaaaae3}),~(\ref{eq:my}) et~(\ref{eq:mz}).$\Box$\vskip5pt

\nocite{*}

\end{document}